\documentclass[a4paper]{scrreprt}
\usepackage{a4}
\usepackage[francais]{babel}
\usepackage[latin1]{inputenc}
\usepackage{amsfonts}
\usepackage{amssymb}
\usepackage{amsmath}
\usepackage[mathscr]{euscript}
\usepackage{rotating}
\include{xy}
\xyoption{all}
\usepackage{vmargin}
\setmarginsrb{3cm}{2.5cm}{2.5cm}{3.5cm}{1cm}{1cm}{0cm}{0.8cm}
\font\petit=cmr10 at 4pt

\def\leq{\leqslant}
\def\geq{\geqslant}

\def\C{\mathbb{C}}
\def\Nat{\mathbb{N}}
\def\R{\mathbb{R}}
\def\Z{\mathbb{Z}}
\def\Q{\mathbb {Q}}
\def\Fi{\mathbb{F}}
\def\Af{\mathbb{A}_f}
\def\Ad{\mathbb{A}}

\def\O{\mathcal{O}}

\def\Hp{{}^{p}{H}}

\def\Hom{\mathrm{Hom}}
\def\Ext{\mathrm{Ext}}

\newcounter{compteurun}[section]
\newcounter{compteurdeux}[subsection]
\makeatletter
\def\thecompteurun{\thesection.\@arabic\c@compteurun}
\def\thecompteurdeux{\thesubsection.\@arabic\c@compteurdeux}
\makeatother

\newenvironment{enonce}[1]{

\vspace{0.3cm}
\noindent{\scshape #1}\,\,\begin{itshape} }{\end{itshape}
\vspace{0.3cm}
}

\newenvironment{enonce_rom}[1]{

\vspace{0.3cm}
\noindent{\scshape #1}\,\, }
{\vspace{0.3cm}}

\newenvironment{propdef}[1][]{\ifthenelse{\equal{\arabic{subsection}}{0}}{\refstepcounter{compteurun}}{\refstepcounter{compteurdeux}}\begin{enonce}{\textbf{Proposition-définition \ifthenelse{\equal{\arabic{subsection}}{0}}{\thecompteurun}{\thecompteurdeux}}}}{\end{enonce}}

\newenvironment{rappel}[1][]{\ifthenelse{\equal{\arabic{subsection}}{0}}{\refstepcounter{compteurun}}{\refstepcounter{compteurdeux}}\begin{enonce}{\textbf{Rappel \ifthenelse{\equal{\arabic{subsection}}{0}}{\thecompteurun}{\thecompteurdeux}}}}{\end{enonce}}

\newenvironment{theoreme}[1][]{\ifthenelse{\equal{\arabic{subsection}}{0}}{\refstepcounter{compteurun}}{\refstepcounter{compteurdeux}}\begin{enonce}{\textbf{Théorème \ifthenelse{\equal{\arabic{subsection}}{0}}{\thecompteurun}{\thecompteurdeux}}}}{\end{enonce}}

\newenvironment{lemme}[1][]{\ifthenelse{\equal{\arabic{subsection}}{0}}{\refstepcounter{compteurun}}{\refstepcounter{compteurdeux}}\begin{enonce}{\textbf{Lemme \ifthenelse{\equal{\arabic{subsection}}{0}}{\thecompteurun}{\thecompteurdeux}}}}{\end{enonce}}
\newenvironment{definition}[1][]{\ifthenelse{\equal{\arabic{subsection}}{0}}{\refstepcounter{compteurun}}{\refstepcounter{compteurdeux}}\begin{enonce}{\textbf{Définition \ifthenelse{\equal{\arabic{subsection}}{0}}{\thecompteurun}{\thecompteurdeux}}}}{\end{enonce}}
\newenvironment{exemple}[1][]{\ifthenelse{\equal{\arabic{subsection}}{0}}{\refstepcounter{compteurun}}{\refstepcounter{compteurdeux}}\begin{enonce}{\textbf{Exemple \ifthenelse{\equal{\arabic{subsection}}{0}}{\thecompteurun}{\thecompteurdeux}}}}{\end{enonce}}
\newenvironment{corollaire}[1][]{\ifthenelse{\equal{\arabic{subsection}}{0}}{\refstepcounter{compteurun}}{\refstepcounter{compteurdeux}}\begin{enonce}{\textbf{Corollaire \ifthenelse{\equal{\arabic{subsection}}{0}}{\thecompteurun}{\thecompteurdeux}}}}{\end{enonce}}
\newenvironment{remarque}[1][]{\ifthenelse{\equal{\arabic{subsection}}{0}}{\refstepcounter{compteurun}}{\refstepcounter{compteurdeux}}\begin{enonce_rom}{\textbf{Remarque \ifthenelse{\equal{\arabic{subsection}}{0}}{\thecompteurun}{\thecompteurdeux}}}}{\end{enonce_rom}}
\newenvironment{remarques}[1][]{\ifthenelse{\equal{\arabic{subsection}}{0}}{\refstepcounter{compteurun}}{\refstepcounter{compteurdeux}}\begin{enonce_rom}{\textbf{Remarques \ifthenelse{\equal{\arabic{subsection}}{0}}{\thecompteurun}{\thecompteurdeux}}}}{\end{enonce_rom}}
\newenvironment{proposition}[1][]{\ifthenelse{\equal{\arabic{subsection}}{0}}{\refstepcounter{compteurun}}{\refstepcounter{compteurdeux}}\begin{enonce}{\textbf{Proposition \ifthenelse{\equal{\arabic{subsection}}{0}}{\thecompteurun}{\thecompteurdeux}}}}{\end{enonce}}

\newenvironment{notation}[1][]{\ifthenelse{\equal{\arabic{subsection}}{0}}{\refstepcounter{compteurun}}{\refstepcounter{compteurdeux}}\begin{enonce}{\textbf{Notation \ifthenelse{\equal{\arabic{subsection}}{0}}{\thecompteurun}{\thecompteurdeux}}}}{\end{enonce}}

\newenvironment{fait}[1][]{\ifthenelse{\equal{\arabic{subsection}}{0}}{\refstepcounter{compteurun}}{\refstepcounter{compteurdeux}}\begin{enonce}{\textbf{Fait \ifthenelse{\equal{\arabic{subsection}}{0}}{\thecompteurun}{\thecompteurdeux}}}}{\end{enonce}}

\newenvironment{souslemme}[1][]{\ifthenelse{\equal{\arabic{subsection}}{0}}{\refstepcounter{compteurun}}{\refstepcounter{compteurdeux}}\begin{enonce}{\textbf{Sous-lemme \ifthenelse{\equal{\arabic{subsection}}{0}}{\thecompteurun}{\thecompteurdeux}}}}{\end{enonce}}

\newenvironment{preuve}{\vspace{0.3cm}\emph{Démonstration.}}{\hfill $\square$ \vspace{0.3cm}}
\newenvironment{preuvet}{\vspace{0.3cm}\emph{Démonstration du théorème.}}{\hfill $\square$ \vspace{0.3cm}}

\newenvironment{preuvepn}[1]{\vspace{0.3cm}\emph{Démonstration de la proposition #1.}} {\hfill $\square$ \vspace{0.3cm}}

\def\nopreuve{\hfill $\square$}

\def\Gm{\mathbb{G}_m}

\def\SD{\mathbb{S}}

\def\T{{\bf T}}
\def\H{{\bf H}}
\def\G{{\bf G}}
\def\GU{{\bf GU}}
\def\SU{{\bf SU}}

\def\GL{{\bf GL}}
\def\SL{{\bf SL}}

\def\P{{\bf P}}
\def\QP{{\bf Q}}
\def\RP{{\bf R}}
\def\B{{\bf B}}

\def\Gr{\mathbb{G}}

\def\N{{\bf N}}
\def\U{{\bf U}}

\def\L{{\bf L}}
\def\M{{\bf M}}
\def\S{{\bf S}}
\def\Ar{\mathrm{A}}

\def\K{\mathrm {K}}

\def\Ze{\mathcal{Z}}
\def\Hr{\mathrm{H}}

\def\F{{\mathcal F}}

\def\DP{{}^{w}D}
\def\as{\underline{a}}
\def\lin{\ell}
\def\Cat{{\mathcal C}}

\def\Mod{{\mathcal M}}
\def\O{{\mathcal O}}
\def\sous{\setminus}
 \def\til{\widetilde}

\def\X{{\mathcal X}}
\def\Y{{\mathcal Y}}
\def\fl{\longrightarrow}
\def\fle{\longmapsto}
\def\flecheinj{\ar@{^{(}->}[r]}
\def\flechesurj{\ar@{->>}[r]}
\def\flinj #1 #2 {\xymatrix{#1\flecheinj & #2}}
\def\flsur #1 #2 {\xymatrix{#1\flechesurj & #2}}
\def\flinja #1 #2 #3{\xymatrix{#1\flecheinj^{#3} & #2}}
\def\flsura #1 #2 #3{\xymatrix{#1\flechesurj^{#3} & #2}}
\def\restr #1 #2 {{#1}_{\mid {#2}}}
\def\flnom #1 {\stackrel {#1} {\longrightarrow}}
\def\doublef #1 #2 #3 #4 {\xymatrix{#1 \ar@<2pt>[r]^{#3}
    \ar@<-2pt>[r]_{#4} & #2}}
\def\dbf #1 #2{\ar@<2pt>[r]^{#1}\ar@<-2pt>[r]_{#2}}
\def\iso{\stackrel {\sim} {\fl}}

\def\limia #1 {\displaystyle{\lim_{\overrightarrow{\scriptstyle #1}}}}
\def\limpa #1 {\displaystyle{\lim_{\overleftarrow{\scriptstyle #1}}}}
\def\limpd {\displaystyle{\lim_{\overleftarrow{\scriptstyle \,\,\,\Delta\,\,\,}}}}
\def\limid {\displaystyle{\lim_{\overrightarrow{\scriptstyle \,\,\,\Delta\,\,\,}}}}

\def\wq{{\circ}\raise 3pt
  \vbox{$\!\!\!\!\!\!\!\!\!\!\!\scriptscriptstyle{\smile}\lower
    6pt\vbox{$\!\!\!\!\!\!\!\!\!\!\!\mkern -1.8mu ($}$}\mkern -917mu}
\def\wa{{\circ}\raise 3pt
\vbox{$\!\!\!\!\!\!\!\!\!\!\!\scriptscriptstyle{\smile}\lower
  2pt\vbox{$\!\!\!\!\!\!\!\!\!\!\mkern -1.8mu${\petit (}}$}\mkern -917mu}

\def\ungras{1\!\!\mkern -1mu1}

\newdimen\lengtharrow            \newbox\uppertext
 \newbox\lowertext               \newbox\arrowbox
   \def\dimmax #1#2{\ifdim #1<#2 #2\else #1\fi}

\def\arrow #1#2%
    {\setbox\uppertext=\hbox{$\scriptstyle #1$}
      \setbox\lowertext=\hbox{$\scriptstyle #2$}
       \lengtharrow=\dimmax{\wd \lowertext}{\wd \uppertext}
        \ifdim \lengtharrow<1pt \setbox\arrowbox=\hbox{$\rightarrow$}
         \else \advance\lengtharrow by 8pt
          \setbox\arrowbox=\hbox to\lengtharrow{\rightarrowfill}
           \ht\arrowbox=3.66875pt \dp\arrowbox=-1.5pt
            \setbox\arrowbox=\hbox{$\mathop{\box\arrowbox}
             \limits_{\box\lowertext}^{\box\uppertext}$}
              \fi
               \mathrel{\box\arrowbox}}

\def\ssi{si et seulement si }

\title{Complexes d'intersection des compactifications de Baily-Borel}
\author{Le cas des groupes unitaires}
\date{}

\begin{document}

{\LARGE\textbf{Complexes d'intersection des compactifications de Baily-Borel}}

\bigskip

{\LARGE Le cas des groupes unitaires sur $\Q$}

\bigskip

\textbf{Sophie Morel}

Laboratoire de mathématique, Université Paris Sud, bâtiment 425, 91405 Orsay Cedex, France

e-mail : sophie.morel@math.u-psud.fr

\vspace{2cm}

\tableofcontents

\setcounter{chapter}{-1}

\chapter{Introduction}

Les variétés de Shimura les plus étudiées sont celles associées au groupe $\GL_2$,
autrement dit les courbes modulaires.
Si $Y$ est une courbe modulaire, on obtient sa compactification de Baily-Borel $j:Y\fl Y^*$
en ajoutant un nombre fini de pointes
(car $\GL_2$ est de rang semi-simple $1$).
Comme $Y^*$ est lisse, le complexe d'intersection associé à un système de coefficients
$\F$ sur $Y$ est $j_*\F$.
Il est possible dans ce cas de calculer la fonction $L$ de $j_*\F$,
et il a été prouvé dans des travaux d'Eichler, Shimura, Deligne et Ihara (entre autres)
qu'elle s'écrit comme un produit alterné de fonctions $L$ de formes modulaires 
et de fonctions zêta.

La fonction $L$ est un produit de facteurs locaux $L_p$, où $p$ parcourt l'ensemble des 
nombres premiers, et ce sont les $L_p$ que l'on calcule. Le cas essentiel est celui où
$Y^*$ et $j_*\F$ ont bonne réduction en $p$. Le facteur local $L_p$ ne dépend alors
que des réductions modulo $p$ de $Y^*$ et $j_*\F$.

Pour calculer $L_p$ dans le cas des courbes modulaires,
il existe deux méthodes : 
la méthode des congruences, qui ne se généralise pas en dimension supérieure,
et la comparaison de la formule des traces d'Arthur-Selberg 
et de la formule des points fixes de Grothendieck-Lefschetz.
C'est cette deuxième méthode que l'on cherche à généraliser.

Pour une variété de Shimura générale $M^{\K}(\G,\X)$,
l'application de cette méthode est plus délicate.
Un premier pas 
est le calcul de la trace d'une puissance du morphisme de Frobenius
sur la cohomologie du complexe d'intersection 
de la compactification de Baily-Borel, ou,
ce qui suffit grâce à la formule des traces de Grothendieck-Lefschetz (cf SGA 4 1/2 Rapport),
de la fonction trace de Frobenius
de ce complexe.

Brylinski et Labesse ([BL]) ont effectué ce calcul 
pour $\G=R_{E/\Q}\GL_2$, avec $E$ une extension totalement réelle 
de $\Q$ de degré supérieur ou égal à $2$, et ils en ont déduit
que la fonction $L$ du complexe d'intersection
était bien de la forme attendue.

Pour $\G=\GU(2,1)$, le calcul a été fait par Kottwitz et Rapoport dans l'article [KR] du livre [LR] (dans ce cas, on peut aussi montrer que la fonction $L$ à coefficients dans le
complexe d'intersection est un produit alterné de fonctions $L$ automorphes, voir
l'article de Langlands et Ramanujan dans le même livre [LR]).
Le cas général d'un groupe de rang $1$ a été traité par Rapoport dans son article [R], paru dans le même livre.

Dans cette thèse, nous traitons le cas des groupes unitaires sur $\Q$ de rang arbitraire,
avec la restriction suivante :
comme on ne dispose pas de compactifications toroïdales des modèles entiers,
on doit exclure un ensemble fini de nombres premiers.

\vspace{1cm}

Du point de vue topologique, 
le complexe d'intersection peut être calculé en utilisant la compactification de Borel-Serre réductive.
Plus précisément, Goresky, Harder et MacPherson construisent une famille de complexes pondérés
sur la compactification de Borel-Serre réductive
et montrent que le complexe d'intersection sur la compactification de Baily-Borel est l'image 
directe de deux de ces complexes pondérés
par le morphisme naturel 
de la compactification de Borel-Serre réductive sur la compactification de Baily-Borel (cf [GHM]).

On peut définir en caractéristique finie des complexes pondérés analogues
aux images directes des complexes pondérés de Goresky, Harder et MacPherson.

La méthode consiste à remplacer les troncatures par les poids pour les actions
des tores centraux des groupes associés aux composantes de bord
par les troncatures par les poids de l'endomorphisme de Frobenius sur les strates de bord. 
Les premiers, à notre connaissance, à utiliser ce type de méthode sont 
Looijenga et
Rapoport dans [LR2].

\vspace{1cm}

Les deux points clé de notre méthode sont le théorème de Pink calculant les restrictions aux strates
de la compactification de Baily-Borel du prolongement d'un système de coefficients 
sur la variété de Shimura (cf [P2])
et la formule suivante :
Si $j:U\fl X$ est l'inclusion d'un ouvert non vide dans un schéma $X$ séparé de type fini sur un corps fini
et $K$ un faisceau pervers pur de poids $a$ sur $U$, alors
$$j_{!*}K=w_{\leq a}Rj_*K.$$
Dans cette formule, $w_{\leq a}$ est le tronqué pour la t-structure $(\DP^{\leq a}(X),\DP^{>a}(X))$,
où $\DP^{\leq a}(X)$ (resp. $\DP^{>a}(X)$) est la sous-catégorie pleine de la catégorie des complexes mixtes sur $X$
dont les objets sont les complexes qui ont tous leurs faisceaux de cohomologie perverse
de poids $\leq a$ (resp. $>a$).
Cette t-structure est assez inhabituelle; par exemple, elle est dégénérée et de coeur nul.

\vspace{1cm}

Pour les variétés de Shimura $M^{\K}(\G,\X)$ considérées dans cette thèse,
le théorème de Pink s'écrit (cf le corollaire \ref{restriction_bord})
$$i^*Rj_*\F^{\K}V\simeq\F^{\K'}R\Gamma(\Gamma_\lin,R\Gamma(Lie(\N),V)).$$
$j$ est l'inclusion de la variété de Shimura dans sa compactification de Baily-Borel $M^{\K}(\G,\X)^*$,
$V$ est une représentation algébrique du groupe $\G$,
$\F^{\K}V$ est le système de coefficients associé à $V$,
$i$ est l'inclusion d'une strate associée à un parabolique maximal $\P$,
$\N$ est le radical unipotent de $\P$
et $\Gamma_\lin$ est un sous-groupe arithmétique de la parte linéaire du quotient de Levi de $\P$.

On peut aussi appliquer ce théorème à des systèmes locaux sur les strates de bord
de la compactification de Baily-Borel,
ce qui permet, par exemple, si $M_2\subset\overline{M}_1\subset M^{\K}(\G,X)^*$ sont deux strates,
de calculer la restriction à $M_2$ de $Ri_{1*}i_1^*Rj_*\F^{\K}V$,
où $i_1$ est l'inclusion de $M_1$ dans $M^{\K}(\G,\X)$.

Rappelons que les poids du faisceau lisse $\F^{\K}V$ sont déterminés par le caractère central de $V$ (cf 3.4).

On a vu ci-dessus qu'on disposait d'une formule pour calculer le prolongement intermédiaire.
Supposons par exemple que le caractère central de $V$ est trivial, donc que $\F^{\K}V$ est pur de poids $0$,
et notons $d$ la dimension de $M^{\K}(\G,\X)$.
$\F^{\K}V[d]$ est alors un faisceau pervers pur de poids $d$ sur $M^{\K}(\G,\X)$,
et le complexe d'intersection $IC_V$ à coefficients dans $V$, 
égal par définition à $(j_{!*}(\F^{\K}V[d]))[-d]$, est donné par la formule suivante :
$$IC_V=w_{\leq d}Rj_*\F^{\K}V.$$

Il s'agit ensuite de calculer le complexe $w_{\leq d}Rj_*\F^{\K}V$.
La t-structure $(\DP^{\leq d},\DP^{>d})$ sur la compactification de Baily-Borel
s'obtient en recollant les t-structures analogues sur la variété de Shimura et les strates du bord.
En utilisant cette remarque et le fait que 
les foncteurs de troncature $w_{\leq d}$ et $w_{>d}$ sont triangulés, 
on obtient la formule suivante dans le groupe de Grothendieck
(qui est un cas particulier du théorème \ref{th:simplification_formule_traces}) :
$$[IC_V]=[w_{\leq d}Rj_*\F^{\K}V]=\sum_{1\leq n_1<\dots<n_r\leq n}(-1)^r[Ri_{n_r*}w_{>d}i_{n_r}^*\dots Ri_{n_1*}w_{>d}i_{n_1}^*Rj_*\F^{\K}V],$$
où $i_k$ est l'inclusion de la $k$-ième strate dans $M^{\K}(\G,\X)^*$
(on a choisi un ordre total sur les strates tel que la dimension soit décroissante).

On peut calculer explicitement les complexes qui apparaissent dans cette somme alternée 
grâce à des applications successives du théorème de Pink.
On obtient finalement le théorème \ref{th:restriction_W_bord} :
\begin{flushleft}$\displaystyle{[i^*IC_V]=}$\end{flushleft}
\begin{flushright}$\displaystyle{\F^{\K_{g,\{r\}}}\left(\sum_{S\subset\{1,\dots,r-1\}}\sum_{i\in I_S}\left[Ind_{\K_{p_ig,S\cup\{r\}}}^{\K_{g,\{r\}}}(-1)^{card(S)}R\Gamma\left(\Gamma_{p_ig,S\cup\{r\}},R\Gamma(Lie(\N_{S\cup\{r\}}),V)_{\geq t_r,< t_s,s\in S}\right)\right]\right).}$\end{flushright}
Les notations précises sont expliquées dans 5.2.
Disons seulement que $i$ est comme plus haut l'inclusion d'une strate de bord de la compactification de Baily-Borel,
que les $t_i$ sont des entiers qui ne dépendent que des dimensions des différentes strates
(on peut prendre par exemple $t_i$ égal à l'opposé
de la codimension de la $i$-ième strate),
que $S\cup\{r\}$ parcourt un système d'indices des sous-groupes paraboliques standard
dont les strates de bord associées dans la compactification de Borel-Serre réductive
s'envoient sur la strate de la compactification de Baily-Borel considérée,
que, si $\P_{S\cup\{r\}}$ est le sous-groupe parabolique correspondant à $S\cup\{r\}$,
alors $\N_{S\cup\{r\}}$ est le radical unipotent de $\P_{S\cup\{r\}}$ 
et les $\Gamma_{p_ig,S\cup\{r\}}$ sont des sous-groupes arithmétique de la partie linéaire du
quotient de Levi $\P_{S\cup\{r\}}/\N_{S\cup\{r\}}$ de $\P_{S\cup\{r\}}$,
et enfin que $R\Gamma(Lie(\N_{S\cup\{r\}}),V)_{\geq t_r,<t_s,s\in S}$ est un tronqué
pour les poids de certains tores centraux de $\P_{S\cup\{r\}}/\N_{S\cup\{r\}}$.

\vspace{1cm}

Passons en revue les différentes parties.

Les trois premières parties contiennent essentiellement des rappels.

Dans la partie 1, 
nous introduisons les groupes unitaires $\GU(p,q)$ que nous comptons étudier
et leurs données de Shimura,
puis nous rappelons la définition de l'ensemble des points complexes des variétés de Shimura associées 
et celle de l'ensemble des points complexes des compactifications de Baily-Borel.
En suivant Pink ([P2] 3.7), nous définissons une stratification du bord
des compactifications de Baily-Borel par des variétés de Shimura associées à des groupes unitaires plus petits.
Enfin, nous rappelons les théorèmes d'algébricité (sur $\C$) pour les variétés ci-dessus
et les théorèmes d'existence des modèles canoniques sur le corps reflex.
Signalons aussi que la section 1.3 contient une liste des
sous-groupes paraboliques standard de $\GU(p,q)$,
et la section 1.5 un calcul de nature combinatoire sur la compactification de Baily-Borel,
qui sert dans la section 5.2, 
et pour lequel nous n'avons pas trouvé de référence.

La partie 2 présente la construction des systèmes de coefficients sur la variété de Shimura
provenant de représentations du groupe.
Nous rappelons d'abord la construction de ces systèmes de coefficients sur les points complexes,
puis nous expliquons une méthode, due à Pink ([P2] 1),
pour montrer que ces systèmes de coefficients proviennent de faisceaux étales sur les modèles canoniques.
Ensuite, nous énonçons le théorème de Pink sur le prolongement de ces faisceaux étales 
à la compactification de Baily-Borel.

La partie 3 est consacrée aux modèles entiers.
Le but est de montrer que, quitte à inverser un ensemble fini de nombres premiers,
la situation des deux premières parties sur le corps reflex 
s'étend à l'anneau des entiers de ce corps.
Comme les variétés de Shimura considérées sont PEL, 
elles ont des modèles sur l'anneau des entiers du corps reflex où
on a inversé le discriminant et certains nombres premiers qui dépendent du niveau.
Une fois qu'on a ces modèles, la construction de Pink des systèmes de coefficients étales
s'étend automatiquement (cf 3.3).
En revanche, pour s'assurer que ces systèmes de coefficients sont bien mixtes (3.4),
pour avoir une compactification de Baily-Borel avec des propriétés convenables (3.2)
et pour que le théorème de Pink reste vrai (3.5),
on doit inverser d'autres nombres premiers, sur lesquels on n'a que très peu d'informations.

La partie 4 est indépendante des autres.
Dans un premier temps, nous y étudions la t-structure $(\DP^{\leq a}(X),\DP^{>a}(X))$ définie plus haut
et y montrons la formule pour le prolongement intermédiaire 
d'un faisceau pervers pur $K$ de poids $a$ sur un ouvert non vide $j:U\fl X$.
Dans un deuxième temps, nous considérons des t-structures sur un schéma stratifié $X$ 
qui s'obtiennent en recollant des t-structures $(\DP^{\leq a'},\DP^{>a'})$ sur les strates (avec un $a'$ qui dépend de la strate).
Grâce aux propriétés de ces t-structures,
nous obtenons si l'ouvert $U$ est réunion de strates une égalité
entre les classes dans le groupe de Grothendieck d'un $w_{\leq a}Rj_*K$ 
et d'une somme alternée de tronqués par le poids qui se calculent sur les strates contenues dans $X-U$ (théorème \ref{th:simplification_formule_traces}).

Dans la partie 5, nous appliquons les résultats de la partie 4 aux variétés de Shimura.
Dans 5.1, nous définissons les complexes pondérés
et montrons que deux de ces complexes sont isomorphes au complexe d'intersection.
En particulier, à l'aide du théorème de Pink et du théorème \ref{th:simplification_formule_traces},
nous obtenons une formule explicite, dans le groupe de Grothendieck,
pour la restriction à une strate de bord de la compactification de Baily-Borel
du complexe d'intersection (plus généralement, d'un complexe pondéré).
Grâce à cette formule et à un résultat de Kottwitz ([K2]), 
nous calculons dans 5.3 la trace d'une puissance de l'endomorphisme de Frobenius
sur la cohomologie d'intersection.

\vspace{1cm}

C'est un plaisir de remercier G. Laumon, qui a passé beaucoup de temps à discuter
avec moi, et les rapporteurs de cette thèse, M. Harris et R. Kottwitz,
qui a corrigé ou simplifié les démonstrations de certains résultats
de la partie 4.

\setcounter{compteurun}{0}
\setcounter{compteurdeux}{0}

\chapter{Variétés de Shimura en caractéristique $0$}

\section{Quelques notations}

Soit $\G$ un groupe algébrique linéaire sur $\Q$.
On note $\S$ son centre déployé.
On appelle \emph{espace symétrique} de $\G$ un espace homogène $X$ sous $\G(\R)$ tel que pour un $x\in X$,
le stabilisateur de $x$ dans $\G(\R)$ soit de la forme $\Ar.\K_\infty$,
où $\Ar=\S(\R)^\circ$ et $\K_\infty$ est un sous-groupe compact maximal de $\G(\R)$
(cette propriété est alors vraie pour tout $x\in X$).
Si $\Gamma$ est un sous-groupe arithmétique de $\G(\Q)$, on dit que $\Gamma\sous X$ est un \emph{espace localement symétrique} associé à $\G$.
 

Pour tout sous-groupe ouvert compact $\K$ de $\G(\Af)$, on pose
$$M^{\K}(\G,X)(\C)=\G(\Q)\sous (X\times\G(\Af)/\K).$$
Si $(g_i)_{i\in I}$ est un système de représentants du double quotient $\G(\Q)\sous\G(\Af)/\K$, on a
$$M^{\K}(\G,X)(\C)=\coprod_{i\in I}\Gamma_i\sous X,$$
où $\Gamma_i=\G(\Q)\cap g_i\K g_i^{-1}$ est un sous-groupe arithmétique de $\G(\Q)$, net si $\K$ est net.

En particulier, si $\K$ est net, $M^{\K}(\G,\X)(\C)$ est une variété analytique réelle.

Si $\K,\K'$ sont deux sous-groupes ouverts compacts de $\G(\Af)$ et $g\in\G(\Af)$ est tel que $\K'\subset g\K g^{-1}$, 
on définit une application analytique finie étale
$$T_g:M^{\K'}(\G,X)(\C)\fl M^{\K}(\G,X)(\C)$$
par : pour tous $x\in X$ et $h\in\G(\Af)$,
$$T_g(\G(\Q).(x,h\K'))=\G(\Q).(x,hg\K).$$

Si $\Gamma,\Gamma'$ sont deux sous-groupes arithmétiques de $\G(\Q)$ et $\gamma\in\G(\Q)$ est tel que $\Gamma'\subset\gamma\Gamma\gamma^{-1}$, 
on définit une application analytique finie étale
$$T_\gamma:\Gamma'\sous X\fl \Gamma\sous X$$
par : pour tout $x\in X$,
$$T_\gamma(\Gamma'.x)=\Gamma.\gamma^{-1}x.$$

\section{Données de Shimura pour certains groupes unitaires}

\begin{notation} Pour $n\in\Nat^*$, on note
$$I=I_n=\left(\begin{array}{ccc}1 & & 0 \\ & \ddots & \\ 0 & & 1\end{array}\right)\in\GL_n(\Z)$$
et
$$J_n=\left(\begin{array}{ccc}0 & & 1 \\ & \begin{turn}{45}\large\ldots\end{turn}& \\ 1 & & 0\end{array}\right)\in\GL_n(\Z).$$

On note $\SD=Res_{\C/\R}\Gm$ le tore de Serre.

\end{notation}

On fixe $n\in\Nat^*$, $p,q\in\Nat$ avec $p+q=n$ et $p\geq q$. 
On note $J_{p,q}$ ou simplement $J$ la matrice
$$J=\left(\begin{array}{ccc}0 & 0 & J_q \\ 0 & I_{p-q} & 0 \\ J_q & 0 & 0\end{array}\right)\in\GL_n(\Z).$$

Soit $E=\Q[i\sqrt{d}]$ ($d\in\Nat^*$ sans facteur carré) une extension quadratique imaginaire 
de $\Q$, dont l'automorphisme non trivial, qui est la conjugaison complexe, 
sera noté \raisebox{7pt}{$\overline{\hspace{0.2cm}}$}; 
on fixe une fois pour toutes des injections 
$E\subset\overline{\Q}\subset\C$ et $\overline{\Q}\subset\overline{\Q}_p$ 
pour tout nombre premier $p$. 
On s'intéresse au groupe algébrique sur $\Q$, $\G=\GU(p,q)$, 
dont les points à valeurs dans une $\Q$-algèbre $A$ sont donnés par
$$\G(A)=\{g\in\GL_n(E\otimes_\Q A),g^* Jg=c(g)J, c(g)\in A^\times\},$$ 
avec, si $g\in \GL_n(E\otimes_\Q A)$, $g^*={}^t \overline{g}$.
Si $q=0$, on note aussi $\GU(p)=\GU(p,q)$.

On a des morphismes de groupes algébriques sur $\Q$:
$$c:\G\fl\Gm\mbox{ et }det:\G\fl Res_{E/\Q}\Gm.$$

On utilisera aussi les groupes $\U(p,q)=Ker(c)$ et $\SU(p,q)=Ker(c)\cap Ker(det)$.

Les groupes $\GU(p,q)$ et $\U(p,q)$ sont réductifs, et $\SU(p,q)$ est semi-simple. 

\begin{fait}
\begin{itemize}
\item[(1)] Le groupe $\G$ est connexe. Le $\Q$-rang semi-simple et le $\R$-rang semi-simple de $\G$ sont tous les deux égaux à $q$.
% \item[(2)] Si $p>q$, on a $c(g)>0$ pour tout $g\in\G(\R)$. Si $p=q$, on a $c(\G(A))=A^\times$ pour toute $\Q$-algèbre $A$. Dans les deux cas, $c(\G(\C))=\C^\times$.
\item[(2)] Si $p>q$, $\G(\R)$ est connexe. Si $p=q$, $\G(\R)$ a deux composantes connexes, $\G(\R)^+=c^{-1}(\R^{+\times})$ et $\G(\R)^-=c^{-1}(\R^{-\times})$.
\item[(3)] Le groupe dérivé de $\G$ est $\G^{der}=\SU(p,q)$. Il est simplement connexe.
\item[(4)] Le quotient $\G/\G^{der}$ est isomorphe au tore
$$\H=\{(x,\lambda),x\overline{x}=\lambda^n\}\subset Res_{E/\Q}\Gm\times\Gm.$$
% Si $n$ est impair, $n=2k+1$, l'application $(x,\lambda)\fle \lambda^k/x$ identifie $\H$ à $Res_{E/\Q}\Gm$, c'est-à-dire qu'on a une suite exacte
% $$1\fl\G^{der}\fl\G\flnom {\kappa} Res_{E/\Q}\Gm\fl 1,$$
% avec $\kappa=c^k/det$.

\end{itemize}
\end{fait}

Introduisons l'espace symétrique $\X$ de $\G$. 
Pour cela, on a besoin d'un certain sous-groupe compact maximal de $\G(\R)$,
 qu'il est plus facile d'écrire pour un sous-groupe de $\G(\C)$ conjugué à $\G(\R).$ Soit donc
$$\til{J}=\left(\begin{array}{cc}I_p & 0 \\ 0 & -I_q \end{array}\right),$$

et
$$\til{G}=\mathrm{GU}(\til{J})=\{g\in\GL_n(\R),g^*\til{J}g=c(g)\til{J},c(g)\in\R^\times\}.$$

Soit
$$u=\left(\begin{array}{ccc} 1/\sqrt{2}I_q & 0       & 1/\sqrt{2}J_q \\
                             0          & I_{p-q} & 0 \\
                             1/\sqrt{2}J_q & 0       & -1/\sqrt{2}I_q \\\end{array}\right).$$
Alors $u=u^*=u^{-1}$ et
et $\til{J}=uJu^{-1},$ donc $\G(\R)\simeq \til{\Gr}$ par $\varphi:g\fle ugu.$

On note $\til{\K}_\infty$ le sous-groupe compact maximal de $\til{\Gr}$ défini par
$$\til{\K}_\infty= \U(p,0)(\R) \times \U(0,q)(\R)=\left\{
\left(\begin{array}{cc} A & 0 \\ 0 & B\end{array}\right),A\in \U(p)(\R),B\in \U(q)(\R)
\right\}.$$

$K_\infty=\varphi^{-1}(\til{\K}_\infty)$ est un sous-groupe compact maximal de $\G(\R).$ 
Soient $\S=\Gm I_n$ le centre déployé de $\G$ et $\Ar=\S(\R)^\circ=\{\lambda I_n,\lambda\in\R^{+\times}\}.$ On pose
$$\X=\G(\R)/\Ar\K_\infty\mbox{ et }x_0=\Ar\K_\infty\in\X.$$

Comme $\Ar\K_\infty\subset\G(\R)^\circ$, on a $\pi_0(\X)=\pi_0(\G(\R))$, 
c'est-à-dire que $\X$ est connexe si $p>q$ et que, si $p=q$, $\X$ a deux composantes connexes (isomorphes) 
$\X^+=\G(\R)^+/\Ar\K_\infty$ et $\X^-=\G(\R)^-/\Ar\K_\infty$. 

On veut définir une application $\G(\R)$-équivariante $h:\X\fl Hom(\SD,\G_\R)$. 
Comme $\G(\R)$ agit transitivement sur $\X$, il suffit de se donner $h_0=h(x_0):\SD\fl\G_\R$ 
tel que $h_0(\S(\R))$ soit dans le centralisateur de $\Ar\K_\infty.$ On prend
$$ h_0:\left\{\begin{array}{ccl} \SD & \fl & \quad\quad\G_\R \\
                             z=a+ib & \fle  &  
\left({\begin{array}{ccc} aI_q & 0 & ibJ_q \\ 0 & zI_{p-q} & 0 \\ ibJ_q & 0 & aI_q\end{array}}\right)\end{array}\right..$$

Si $\til{h}_0=\varphi\circ h_0$, on a
$$\til{h}_0(z)=\left(\begin{array}{cc}zI_p & 0 \\ 0 & \overline{z}I_q\end{array}\right).$$

On suppose désormais que $q\geq 1$ ou $q=0$ et $p=1$
(si $p\geq 2$, $\GU(p)$ ne peut pas être le groupe d'une donnée de Shimura, car il est de type compact modulo son centre).

\begin{fait} 
$(\G,\X,h)$ est une donnée de Shimura pure (cf [P1] 2.1).
\end{fait}

Si $\K$ est un sous-groupe compact ouvert de $\G(\Af)$, l'ensemble des points complexes de la variété de Shimura associée est
$$M^\K(\G,\X)(\C)=\G(\Q)\sous (\X\times\G(\Af)/\K).$$

\begin{proposition} ([P1] 3.3)
$M^{\K}(\G,\X)(\C)$ a une structure canonique d'espace complexe normal. 
Si de plus $\K$ est net (cf [P1] 0.6), $M^{\K}(\G,\X)(\C)$ est une variété analytique complexe.

\end{proposition}

\section{Sous-groupes paraboliques de $\GU(p,q)$}

On note $E_{ij}\in M_n(\Z)$, $1\leq i,j\leq n$, les matrices élémentaires : 
pour tous  $k,l\in\{1,\dots,n\}$, $(E_{ij})_{kl}=\delta_{ik}\delta_{jl}$.
Un tore maximal de $\G=\GU(p,q)$ est le tore $\T$ tel que
$$\T(\Q)=\left\{\left(\begin{array}{ccc}\lambda_1 & 0 & 0 \\ 0 & \ddots & 0 \\ 0 & 0 & \lambda_n\end{array}\right),\lambda_1,\dots,\lambda_n\in E^\times,\lambda_1\overline{\lambda_n}=\dots=\lambda_q\overline{\lambda_{p+1}}=\lambda_{q+1}\overline{\lambda_{q+1}}=\dots=\lambda_p\overline{\lambda_p}\in\Q^\times\right\}.$$
Le sous-tore déployé maximal $\S$ de $\T$ vérifie
$$\S(\Q)=\left\{\lambda\left(\begin{array}{ccccccc}\lambda_1 & 0 & 0 & 0 & 0 & 0 & 0 \\
0 & \ddots & 0 & 0 & 0 & 0 & 0 \\
0 & 0 & \lambda_q & 0 & 0 & 0 & 0 \\
0 & 0 & 0 &  I_{p-q} & 0 & 0 & 0 \\
0 & 0 & 0 & 0 & \lambda_q^{-1} & 0 & 0 \\
0 & 0 & 0 & 0 & 0 & \ddots & 0 \\
0 & 0 & 0 & 0 & 0 & 0 & \lambda_1^{-1} \end{array}\right),\lambda,\lambda_1,\dots,\lambda_q\in\Q^\times\right\}.$$
Un parabolique minimal contenant $\S$ est le groupe $\P$ dont les $\Q$-points sont
$$\P(\Q)=\left\{\left(\begin{array}{ccc}A & * & * \\
0 & B & * \\
0 & 0 & C\end{array}\right),A,C\in \B_q(E),B\in\GL_{p-q}(E)\right\}\cap\GU(p,q)(\Q),$$
où $\B_q\subset\GL_q$ est le groupe des matrices triangulaires supérieures.

On dit qu'un sous-groupe parabolique de $\G$ est \emph{standard} s'il contient $\P$.
On sait qu'un sous-groupe parabolique de $\G$ est conjugué par $\G(\Q)$ à un unique
sous-groupe parabolique standard,
donc il suffit de décrire les sous-groupes paraboliques standard.
Il sont indexés par les sous-ensembles $I\subset\{1,\dots,q\}$ de la manière suivante.

% Pour tout $i\in\{1,\dots,q\}$, on note
% $$\epsilon_i:\left(\begin{array}{ccccccc}\lambda\lambda_1 & 0 & 0 & 0 & 0 & 0 & 0 \\
% 0 & \ddots & 0 & 0 & 0 & 0 & 0 \\
% 0 & 0 & \lambda\lambda_q & 0 & 0 & 0 & 0 \\
% 0 & 0 & 0 & \lambda I_{p-q} & 0 & 0 & 0 \\
% 0 & 0 & 0 & 0 & \lambda\lambda_q^{-1} & 0 & 0 \\
% 0 & 0 & 0 & 0 & 0 & \ddots & 0 \\
% 0 & 0 & 0 & 0 & 0 & 0 & \lambda\lambda_1^{-1}\end{array}\right)\fle\lambda_i.$$
% On a $2\epsilon_i\in X^*(\S)$, et $\epsilon_i\in X^*(\S)$ \ssi $p>q$. Supposons d'abord $p=q$. Alors
% $$\Phi(\S,\G)=\Phi_{nd}=\{\epsilon_i-\epsilon_j,i,j\in\{1,\dots,q\},i\not =j\}\cup\{\pm (\epsilon_i+\epsilon_j),i,j\in\{1,\dots,q\}\}$$
% $$\Phi^+=\{\epsilon_i-\epsilon_j,i<j\}\cup\{\epsilon_i+\epsilon_j\}$$
% $$\Delta=\{\alpha_i=\epsilon_i-\epsilon_{i+1},1\leq i\leq q-1\}\cup\{2\epsilon_q\}.$$
% Soient $i,j\in\{1,\dots,q\}$. Si $i\not =j$, alors
% $$\Lg_{\epsilon_i-\epsilon_j}=\{\lambda E_{ij}-\overline{\lambda}E_{n+1-j,n+1-i},\lambda\in E\}$$
% $$\Lg_{\epsilon_i+\epsilon_j}=\{\lambda E_{i,n+1-j}-\overline{\lambda}E_{j,n+1-i},\lambda\in E\}$$
% $$\Lg_{-\epsilon_i-\epsilon_j}=\{\lambda E_{n+1-i,j}-\overline{\lambda}E_{n+1-j,i},\lambda\in E\}.$$
% Enfin,
% $$\Lg_{2\epsilon_i}=\{\lambda E_{i,n+1-i},\lambda\in E,\lambda+\overline{\lambda}=0\}$$
% $$\Lg_{-2\epsilon_i}=\{\lambda E_{n+1-i,i},\lambda\in E,\lambda+\overline{\lambda}=0\}.$$
Soit $I\subset\{1,\dots,q\}$. On écrit
$$I-\{q\}=\{q_1,\dots,q_1+\dots q_{r-1}\}$$
et on note $q_r=q-(q_1+\dots+q_{r-1})$. 
On pose
$$\P_I=\left(\begin{array}{ccccc}R_{E/\Q}\GL_{q_1} & * & * & * & * \\
0 & \ddots & * & * & * \\
0 & 0 & \GU(q_r,q_r) & * & * \\
0 & 0 & 0 & \ddots & * \\
0 & 0 & 0 & 0 & R_{E/\Q}\GL_{q_1}\end{array}\right)\cap\GU(q,q).$$
En particulier, les sous-groupes paraboliques maximaux standard de $\GU(p,q)$ sont les
$$\P_r=\P_{\{r\}}=\left(\begin{array}{ccc}R_{E/\Q}\GL_r & * & * \\
0 & \GU(p-r,q-r) & * \\
0 & 0 & R_{E/\Q}\GL_r\end{array}\right)\cap\GU(p,q)$$
pour $r\in\{1,\dots,q\}$.

\section{Compactification de Baily-Borel}

\subsection{Description}

Soit $(\G,\X)$ la donnée de Shimura pure de la section 1.2.

On obtient la compactification de Baily-Borel des $M^{\K}(\G,\X)(\C)$ en ajoutant à $\X$ des composantes rationnelles de bord, indexées par les sous-groupes paraboliques maximaux de $\G$,
 pour former une compactification partielle $\X^*$, en étendant l'action de $\G(\Q)$ à $\X^*$ et en faisant le double quotient
$$\G(\Q)\sous (\X^*\times \G(\Af)/\K).$$
On suit ici [P1] chapitres 4 et 6 (dont les résultats reposent sur ceux de [AMRT]).

Dans [P1] 4.6, Pink introduit pour chaque sous-groupe parabolique maximal (admissible, mais ici ils le sont tous, car $G^{ad}$ est simple) $\P$ de $\G$
un morphisme $\omega:\SD_\C\fl\P_\C$ (noté $\omega\circ h_\infty$ par Pink).

La construction de ce morphisme repose sur celle du morphisme de Harish-Chandra\newline
 $\U^1\times\SL_2(\R)^q\fl\G^{ad}(\R)$ de [AMRT] th 2 (ii) p 177-178.
Le lemme suivant résulte facilement de la définition de ce morphisme.

\begin{lemme} Le morphisme de Harish-Chandra est l'image du morphisme\newline
$\varphi:\U^1\times\SL_2(\R)\fl\G(\R)$ défini par
$$\varphi\left(u,\left(\begin{array}{cc}a_1 & b_1 \\ c_1 & d_1 \end{array}\right),\dots,\left(\begin{array}{cc}a_q & b_q \\ c_q & d_q \end{array}\right)\right)=\left(\begin{array}{ccccccccc}
a_1 &        & 0       &   &        &   & 0    &        & ib_1 \\
      & \ddots &       &   & 0      &   &      &\begin{turn}{45}\large\ldots\end{turn}& \\
  0   &        & a_q   &   &        &   & ib_q &        & 0    \\
      &        &       & u &        & 0 &      &        &      \\
      & 0      &       &   & \ddots &   &      &   0    &      \\
      &        &       & 0 &        & u &      &        &      \\
0     &        & -ic_q &   &        &   & d_q  &        & 0    \\
      &\begin{turn}{45}\large\ldots\end{turn}& &   &   0    &   &      & \ddots &      \\
-ic_1 &        &   0   &   &        &   & 0    &        & d_1
\end{array}\right).$$

\end{lemme}

\begin{corollaire}
Pour $\P=\P_r$, $1\leq r\leq q$, le morphisme $\SD_\C\fl\P_\C$ de [P1] 4.6, que l'on notera $\omega_r$, associe à $(z,z')=1\otimes x+i\otimes y\in\SD(\C)=(\C\otimes_\R\C)^\times$, 
avec $x=(z+z')/2$ et $y=(z-z')/2i$,  la matrice
$$\left(\begin{array}{ccccc} 
Z I_r & 0 & 0 & 0 & (1\otimes 1-Z) J_r \\ 
0 & 1\otimes x I_{q-r} & 0 & i\otimes y J_{q-r} & 0 \\ 
0 & 0 & (1\otimes x  +i\otimes y)I_{p-q} & 0 & 0 \\ 
0 & i\otimes y J_{q-r} & 0 & 1\otimes x I_{q-r} & 0 \\ 
0 & 0 & 0 & 0 & I_r\end{array}\right)$$
où $Z=1\otimes Re(zz')+i\otimes Im(zz')$.
\end{corollaire}

\begin{definition}\label{def:partie_lineaire} Soient $\P$ un sous-groupe parabolique maximal de $\G$,
$\N$ son radical unipotent et $\L=\P/\N$ son quotient de Levi.
La \emph{partie hermitienne} $\L_h$ de $\L$ est le centralisateur dans $\L$ du centre $\U$ de $\N$
(pour l'action de $\L$ sur $\U$ déduite de l'action par conjugaison de $\P$ sur $\N$).
Il existe un unique sous-groupe $\L_\lin$ de $\L$, appelé \emph{partie linéaire} de $\L$,
tel que $\L_h$ et $\L_\lin$ commutent, que $\L=\L_h\L_\lin$ 
et que $\L_h(\R)\cap\L_\lin(\R)$ soit fini (cf [AMRT] p 221-222).

\end{definition}

Soient $\P$ un sous-groupe parabolique maximal de $\G$ et $\N$ son radical unipotent.
Dans [P1] 4.7, Pink introduit le plus petit sous-groupe sous-groupe distingué $\QP_{\C}$ de $\P_\C$ qui est défini sur $\Q$ et contient l'image de $\omega$
(Pink note ce sous-groupe $\P_1$).
$\QP$ contient $\N$ (d'après la preuve de [P1] 4.8),
et $\QP/\N$ est, à un facteur de type compact près et au centre près, la partie hermitienne de $\P/\L$ (cf [AMRT] III.4.1).
On peut sans rien changer agrandir le centre de $\QP/\N$ (cf la remarque (ii) de [P1] 4.11).
On utilisera donc la définition suivante de $\QP$ :

\begin{definition}\label{def_Qr} Soit $r\in\{1,\dots,q\}$.
Si $r<q$ on pose 
$$\QP_r=\left\{\left(\begin{array}{ccc}\lambda I_r & * & * \\ 0 & D & * \\ 0 & 0 & I_r\end{array}\right),D\in\GU(p-r,q-r),\lambda=c(D)\right\},$$
et, si $r=q$,
$$\QP_q=\left\{\left(\begin{array}{ccc}\lambda\overline{\lambda} I_q & * & * \\ 0 & \lambda I_{p-q} & * \\ 0 & 0 & I_q \end{array}\right),\lambda\in E^*\right\}.$$
Soit $\P$ un sous-groupe parabolique maximal de $\G$.
Il existe $g\in\G(\Q)$ et $r\in\{1,\dots,q\}$ tels que $\P=g\P_r g^{-1}$.
On pose alors $\QP=g\QP_r g^{-1}$
(cette définition ne dépend pas de $g$, car $\QP_r$ est distingué dans $\P_r$).

\end{definition}

Soit $r\in\{1,\dots,q\}$.
On note $\N_r$ le radical unipotent de $\P_r$ (qui est aussi celui de $\QP_r$)
et $\U_r$ le centre de $\N_r$.
Construisons l'espace principal homogène sous $\QP_r(\R)\U_r(\C)$ de [P1] 4.11. 
On considère l'application $\P_r(\R)$-équivariante
$$\begin{array}{rcl}\X=\P_r(\R)/(\P_r(\R)\cap \Ar\K_\infty) & \fl & \pi_0(\X)\times Hom(\SD_\C,\QP_{r,\C}) \\
x=g(\P_r(\R)\cap \Ar\K_\infty) & \fle & ([x],int(g)\circ\omega_r)\end{array}.$$
L'image de cette application est contenue dans une $\QP_r(\R)\U_r(\C)$-orbite 
($\QP_r(\R)\U_r(\C)$ agit sur $\pi_0(\X)$ par $\pi_0(\QP_r(\R)\U_r(\C))=\pi_0(\QP_r(\R))\fl\pi_0(\P_r(\R))$), 
qu'on note $\Y_r$. 
On note $h_r:\Y_r\fl Hom(\SD_c,\QP_{r,\C})$ la deuxième projection.

\begin{proposition} ([P1] 4.11)
$(\QP_r,\Y_r,h_r)$ est une donnée de Shimura mixte (cf [P1] 2.1).

\end{proposition}

$(\G_r,\X_r)=(\QP_r,\Y_r)/\N_r$ (cf [P1] 2.9 pour la définition du quotient)
est donc une donnée  de Shimura pure.
Si $r<q$, elle est isomorphe à la donnée de la section 1.2 pour $\protect{\G=\GU(p-r,q-r)}$.
Si $r=q$, elle est isomorphe à la donnée de la section 1.2 pour $\G=\GU(1)=R_{E/\Q}\Gm$.

\begin{definition}\label{def:partie_hermitienne} Si $1\leq r<q$, on note $\G'_r=\G_r$.
Pour $r=q$, on note
$$\G'_q=\left\{\left(\begin{array}{ccc}\lambda I_q & 0 & 0 \\
0 & D & 0 \\
0 & 0 & I_q\end{array}\right),D\in\GU(p-q),\lambda=c(D)\right\}.$$
On voit $\G'_q$ comme un sous-groupe de $\P_q/\N_q$; $\G_q$ est alors le centre de $\G'_q$, et $\G'_q/\G_q$ est de type compact.

\end{definition}

Si $\P$ est un sous-groupe parabolique maximal de $\G$
de radical unipotent $\N$,
il est conjugué à l'un des $\P_r$, donc on peut aussi lui associer une donnée de Shimura mixte $(\QP,\Y)$
et une donnée de Shimura pure $(\QP,\Y)/\N$.

\begin{definition}([P1] 4.11) Une \emph{composante rationnelle de bord} de $(\G,\X)$ est une donnée de Shimura mixte $(\QP,\Y)$ associée à un sous-groupe parabolique maximal de $\G$.

\end{definition}

% Soit $r\in\{1,\dots,q\}$.
% L'application $\X\fl\Y_r$ est une immersion ouverte holomorphe ([P1] 4.15). 
% On peut déterminer l'image de $\X$ dans $\Y_r$ grâce à l'application partie imaginaire ([P1] 4.14) : 
% On pose $\U_r(\R)(-1)=(2i\pi)^{-1}\U_r(\R)\subset\U_r(\C)$ ($\U_r$ étant commutatif et unipotent, on l'a identifié à son algèbre de Lie). 
% Si $x\in\Y_r$, soit $g\in\QP_r(\R)\U_r(\C)$ un représentant de $x$. 
% Il existe un unique $u_x\in\U_r(\R)(-1)$ tel que $int(u_x)^{-1}\circ int(g)\circ\omega_r$ soit défini sur $\R$; ce $u_x$ ne dépend que de $x$, et on le note $im(x)$. 
% On obtient ainsi une application
% $$im:\Y_r\fl\U_r(\R)(-1).$$

% Si $y\in\Y_r$ est représenté par $g\in\QP_r(\R)\U_r(\C)$, avec
% $$g=\left(\begin{array}{ccc}aI_r & B & C \\ 0 & D & E \\ 0 & 0 & I_r\end{array}\right),$$
% alors
% $$im(y) =a J_r + \frac {1} {2} (C+J_rC^*J_r+J_rE^*J_{p-r,q-r}E).$$

% L'image de $\X$ dans $\Y_r$ est l'image réciproque par $im$ d'un cône ouvert convexe $C(\X,\QP_r)$ de $\U(\R)(-1)$, 
% qui est une orbite de $\U(\R)(-1)$ sous l'action par conjugaison de $\P_r(\R)^\circ$ ([P1] prop 4.15). 
% Comme $im(x_0)=I_r$ (on rappelle que $x_0=\Ar\K_\infty\in\X$), $C(\X,\QP_r)$ est la classe de $\P_r(\R)^\circ$-conjugaison de $J_r$. 
% En faisant le calcul, on trouve
% $$C(\X,\QP_r)=\{AA^*J_r,A\in\GL_r(\C)\}=\Pos_r. J_r,$$
% où $\Pos_r$ est l'ensemble des matrices hermitiennes définies positives de taille $r$.

\begin{definition} On pose
$$\X^*=\X\sqcup\coprod_{(\QP,\Y)}\Y/\N,$$
où la somme est sur l'ensemble des composantes rationnelles de bord de $(\G,\X)$ et, si $(\QP,\Y)$ est une telle composante, $\N$ est le radical unipotent de $\QP$
et $(\QP,\Y)/\N=(\QP/\N,\Y/\N)$.

On munit $\X^*$ de la topologie de Satake (cf [P1] 6.2).

% Si $(\QP,\Y)$ est une composante rationnelle de bord et $\N$ est le radical unipotent de $\QP$, on note $\Ze=\Y/\N$ et $\psi:\Y\fl\Ze$ la projection.

% Soit $(\QP,\Y)$ une composante de bord. Soit $D\subset C(\X,\QP)$ un coeur («core» en anglais, cf [P1] 6.1), 
% c'est-à-dire un sous-ensemble convexe tel que, 
% pour un choix de $\U(\Z)(-1)\subset\U(\Q)(-1)=1/(2\pi i)\U(\Q)$, 
% si $D_0$ est l'enveloppe convexe de $\U(\Z)(-1)\cap C(\X,\QP)$, 
% il existe $\lambda,\mu\in\R^{+*}$ tels que $\lambda D_0\subset D\subset\mu D_0$. 
% Pour tout $V\subset\Ze$ et tout $\lambda\in\R^{+*}$, on note
% $$\Vois(V,\lambda)=\coprod_{(\QP',\Y')} \psi'(im^{-1}(\lambda D)\cap \psi^{-1}(V))\subset\X^*,$$
% où $(\QP',\Y')$ parcourt l'ensemble des composantes de bord entre $(\QP,\Y)$ et $(\G,\X)$ (c'est-à-dire telles que $\QP\subset\QP'$). 

% Soit $z\in\Ze$. Les
% $$\Vois(V,\lambda),$$
% où $V$ parcourt une base de voisinages de $z$ dans $\Ze$ et $\lambda$ parcourt $\R^{+*}$, forment une base de voisinages de $z$ dans $\X^*$.

% L'action continue de $\G(\Q)$ sur $\X$ s'étend en une action continue sur $\X^*$.

\end{definition}

\begin{remarque} Si $(\QP,\Y)$ et $(\QP',\Y')$ sont deux composantes de bord, et si on note $\Ze$ et $\Ze'$ les sous-ensembles correspondants de $\X^*$, il résulte immédiatement de la définition de la topologie dans [P1] 6.2 que $\Ze'\subset\overline{\Ze}$ si et seulement si $\QP'\subset\QP$.

\end{remarque}

% \begin{exemple} On prend $q=1$, $p=n-1\geq 2$. $\X_1$ est un point, $\U_1(\C)=\C$, $\U_1(\R)(-1)=\R$ et $C(\X,\QP_1)=\R^{+*}$; on prend $D=[1,\infty[$. L'application $im$ est donnée sur $\X$ par : si $x=g\Ar\K_\infty$ avec
% $$g=\left(\begin{array}{ccc}a & * & * \\ 0 & * & * \\ 0 & 0 & f\end{array}\right),$$
% alors $im(x)=a/f$. Comme on divise par $\Ar\K_\infty$, on peut supposer $f=1$. L'image de $g\Ar\K_\infty$ dans la réalisation bornée de $\X$ est alors
% $$\left(\frac {a-c-1} {a-c+1},\frac {-\sqrt{2} {}^tE} {a-c+1} \right)\in\C^{n-1}.$$

% Une base de voisinages de l'unique point $(1,0,\dots,0)$ de $\Ze_1$ est
% $$\left\{\left(\frac {a-c-1} {a-c+1},\frac {-\sqrt{2} {}^tE} {a-c+1}\right),a\in\R^{+*},a>\lambda,c\in\C,E\in M_{n-2,1}(\C),c+\overline{c}+E^*E=0 \right\},$$
% pour $\lambda$ parcourant $\R^{+*}$.

% \end{exemple}

\begin{definition} Soit $\K$ un sous-groupe ouvert compact de $\G(\Af)$.
La compactification de Baily-Borel de $M^{\K}(\G,\X)(\C)$ est
$$M^{\K}(\G,\X)^*(\C)=\G(\Q)\sous(\X^*\times \G(\Af)/\K).$$

Elle a naturellement une stratification par des espaces localement symétriques (isomorphes à des $\Gamma\sous\X_r$, avec $\Gamma$ un sous-groupe arithmétique de $\G_r(\Q)$).

\end{definition}

\begin{theoreme} ([P1] 6.2, [BB] th 10.4)
Il existe une unique structure d'espace analytique complexe normal compact sur $M^{\K}(\G,\X)^*(\C)$ dont la restriction à une strate isomorphe à $\Gamma\sous\X_r$ est la structure complexe induite par celle de $\X_r$.

\end{theoreme}

De plus :

\begin{proposition} ([P1] 6.2)
Soient $\K'$, $\K$ des sous-groupes compacts ouverts de $\G(\Af)$ et $g\in\G(\Af)$ tels que $\K'\subset g\K g^{-1}$. Le morphisme
$$T_g:M^{\K'}(\G,\X)(\C)\fl M^{\K}(\G,\X)(\C)$$
se prolonge par continuité en une application holomorphe finie
$$M^{\K'}(\G,\X)^*(\C)\fl M^{\K}(\G,\X)^*(\C),$$
qu'on notera $\overline{T}_g$.

\end{proposition}

\subsection{Stratification du bord}

Décrivons, en suivant [P2] 3.7, une stratification du bord de $M^{\K}(\G,\X)^*(\C)$.

Les $(\QP_r,\Y_r)$, $1\leq r\leq q$, forment un ensemble de représentants des classes de $\G(\Q)$-conjugaison de composantes de bord de $(\G,\X)$. 
Soient $r\in\{1,\dots,q\}$ et $g\in\G(\Af)$. 
Notons $\pi:\P_r\fl\P_r/\N_r=\G_r$ la projection, $\K_Q=\QP_r(\Af)\cap g\K g^{-1}$ et $\K_G=\pi(\K_Q)$. 
Le morphisme holomorphe $i_{r,g}:M^{\K_G}(\G_r,\X_r)(\C)\fl\M^\K(\G,\X)^*(\C)$ est défini par le diagramme suivant :
$$\xymatrix{\G_r(\Q)\sous\X_r\times (\G_r(\Af)/\K_G) & [(x,\pi(h)] \\
\QP_r(\Q)\sous\X_r\times (\QP_r(\Af)/\K_Q)\ar[u]^{\wr}\ar[d]_{\wr} & [(x,h)]\ar@{|->}[u]\ar@{|->}[d] \\
\P_r(\Q)\sous \X_r\times (\P_r(\Q)\QP_r(\Af)/\K_Q)\ar[d] & [(x,h)]\ar@{|->}[d] \\
\G(\Q)\sous \X^*\times (\G(\Af)/\K) & [(x,hg)]}.$$

Notons $\Hr_P=g\K g^{-1}\cap\P_r(\Q)\QP_r(\Af)$ et $\Hr_{\lin}=g\K g^{-1}\cap Cent_{\P_r(\Q)}(\X_r)\N(\Af)$. 
On fait agir $\Hr_P$ sur $M^{\K_G}(\G_r,\X_r)(\C)$ par multiplication à droite sur le deuxième facteur dans la troisième ligne du diagramme ci-dessus. 
Les sous-groupes $\K_Q$ et $\Hr_\lin$ de $\Hr_P$ agissent trivialement. 
Comme $\Hr_P/\K_Q\Hr_\lin$ est fini, $i_{r,g}$ est fini sur son image.

\begin{lemme}([P2] (3.7))
\begin{itemize}
\item[(i)] $i_{r,g}$ induit une immersion localement fermée $M^{\K_G}(\G_r,\X_r)(\C)/\Hr_P\fl M^{K}(\G,\X)^*(\C)$.
\item[(ii)] Si $\K$ est net, le groupe fini $\Hr_P/\K_Q\Hr_\lin$ agit librement sur $M^{\K_G}(\G_r,\X_r)(\C)$.
\end{itemize}

\end{lemme}

\begin{remarque}
Si $\K$ est net, $\Hr_P/\K_Q$ est un sous-groupe arithmétique net de $(\P_r/\QP_r)(\Q)$. Or, si $r=1$, tous les sous-groupes arithmétiques de $(\P_r/\QP_r)(\Q)$ sont finis; donc dans ce cas on a $\K_Q=\Hr_P$,
et le morphisme $i_{r,g}$ est une immersion localement fermée.

\end{remarque}

Les images des morphismes $i_{r,g}$ forment une partition du bord
de $M^{\K}(\G,\X)^*(\C)$,
et les images de $i_{r,g}$ et $i_{r,g'}$ sont les mêmes si et seulement si $\P_r(\Q)\QP_r(\Af)g\K=\P_r(\Q)\QP_r(\Af)g'\K$. 

On a donc obtenu une stratification du bord de la compactification de Baily-Borel par des quotients de variétés de Shimura par des groupes finis. 

On peut réécrire d'une manière un peu plus simple les définitions de $\Hr_P$ et $\Hr_\lin$.
On pose, pour $1\leq r\leq q-1$, 
$$\L_{\lin,r}=\L'_{\lin,r}=\left(\begin{array}{ccc}* & 0 & 0 \\
0 & I_{n-2r} & 0 \\
0 & 0 & *\end{array}\right)$$
et 
$$\L'_{\lin,q}=\left(\begin{array}{ccc}* & 0 & 0 \\
0 & I_{n-2q} & 0 \\
0 & 0 & * \end{array}\right)$$
$$\L_{\lin,q}=\left(\begin{array}{ccc}* & 0 & 0\\
0 & \SU(p-q) & 0 \\
0 & 0 & *\end{array}\right)$$
(les blocs diagonaux sont carrés de tailles $r,n-2r,r$ resp. $q,n-2q,q$).
Le groupe $\L'_{\lin,r}$ est, au centre près, la partie linéaire de $\L_r$ (cf la définition \ref{def:partie_lineaire}).

$\L_r$ est produit direct de $\L'_{\lin,r}$ et de $\G'_r$, 
et produit quasi-direct de $\L_{\lin,r}$ et de $\G_r$.
Si $1\leq r\leq q-1$, 
on a $\L_{\lin,r}=\L'_{\lin,r}$ et de $\G_r=\G'_r$,
et $\L_r$ est simplement produit direct de ces deux sous-groupes.

\begin{lemme} Soient $\P$ un sous-groupe parabolique maximal de $\G$ et $g\in\G(\Af)$. 
$\P$ est conjugué à un unique $\P_r$.
On note $\N$ le radical unipotent de $\P$, $\L=\P/\N$ le quotient de Levi, $\QP$ le sous-groupe distingué défini par Pink (cf plus haut), $\L_\lin,\L'_\lin$ les sous-groupes de $\L$ obtenus à partir de $\L_{\lin,r},\L'_{\lin,r}$ par conjugaison, $\Ze$ la composante de bord de $\X^*$ associée à $\P$,
$$\Hr_P=g\K g^{-1}\cap \P(\Q)\QP(\Af)$$
$$\Hr_\lin=g\K g^{-1}\cap Cent_{\P(\Q)}(\Ze)\N(\Af).$$
Alors, si $\K$ est net,
\begin{itemize}
\item[(i)] $\Hr_P=g\K g^{-1}\cap\L'_\lin(\Q)\QP(\Af)=g\K g^{-1}\cap\L_\lin(\Q)\QP(\Af)$;
\item[(ii)] $\Hr_\lin=g\K g^{-1}\cap\L'_\lin(\Q)\N(\Af)=g\K g^{-1}\cap\L_\lin(\Q)\N(\Af)$.
\end{itemize}
\end{lemme}

\begin{preuve} On se ramène par conjugaison au cas où $\P=\P_r$.
Supposons d'abord $r<q$. Alors $\P_r(\Q)=\L_{\lin,r}(\Q)\QP_r(\Q)$ et $\QP_r(\Q)\subset\QP_r(\Af)$, d'où (i) (on n'utilise pas dans ce cas la netteté de $\K$). 
Prouvons (ii). 
L'image de $\Hr_\lin$ par la projection $\P_r(\Q)\QP_r(\Af)\fl\P_r(\Q)\QP_r(\Af)/\L_{\lin,r}(\Q)\N_r(\Af)\simeq\G_r(\Af)$ 
est un sous-groupe compact (pour la topologie induite par $\G_r(\Af)$) de $\G_r(\Q)$, dont les éléments sont nets et agissent trivialement sur $\X_r$, 
donc c'est un sous-groupe du centre de $\G_r(\Q)$; comme le centre de $\G_r$ est une extension d'un tore déployé par un tore de type compact, ce sous-groupe est trivial. 

Traitons le cas $r=q$, et montrons les premières égalité de (i) et (ii).
$\G'_q$ est isomorphe à $\GU(p-q)$, 
et $\L_q\simeq\G'_q\times\L'_{\lin,q}$. 
De plus, $\G_q$ est le centre de $\G'_q$. 
L'image de $\Hr_P$ par la projection $\P_q(\Af)\fl(\P_q/\L_{\lin,q}\QP_q)(\Af)\simeq (\G'_q/\G_q)(\Af)$ est un sous-groupe net de $(\G_q/\G'_q)(\Q)$; comme $(\G_q/\G'_q)(\R)$ est compact, cette image est triviale, d'où (i). Pour (ii), on remarque que $Cent_{\P_q(\Q)}(\X_q)=\P_q(\Q)$, donc l'image de $\Hr_\lin$ dans $\G'_q(\Af)$ est un sous-groupe arithmétique net de $\G'_q(\Q)$, forcément trivial puisque $\G'_q\simeq\GU(p-q)$.
Enfin, on peut remplacer $\L'_{\lin,q}(\Q)$ par $\L_{\lin,q}(\Q)$ dans (i) et (ii)
car $\L_{\lin,q}=\L'_{\lin,q}\times\SU(p-q)$
(et $\K$ est net).

\end{preuve}

% \begin{corollaire} $\Hr_P/\Hr_\lin$ est un sous-groupe ouvert compact de $\G_r(\Af)$, qu'on notera $\K^0_G$. Il contient $\K_G$. Par définition de $\K^0_G$, $M^{\K_G}(\G_r,\X_r)(\C)/\Hr_P=M^{\K^0_G}(\G_r,\X_r)(\C)$. $i_{r,g}$ donne une immersion localement fermée holomorphe $i^0_{r,g}:M^{\K^0_G}(\G_r,\X_r)(\C)\fl M^{\K}(\G,\X)^*(\C)$.

% \end{corollaire}

Enfin, on a :

\begin{proposition}([P1] 7.6) Le morphisme $i_{r,g}:M^{\K_G}(\G_r,\X_r)(\C)\fl M^{\K}(\G,\X)^*(\C)$ se prolonge en un morphisme fini holomorphe
$$\overline{i_{r,g}}:M^{\K_G}(\G_r,\X_r)^*(\C)\fl M^{\K}(\G,\X)^*(\C),$$
dont l'image est l'adhérence de l'image de $i_{r,g}$. 
% De même, l'immersion $i^0_{r,g}$ se prolonge en un morphisme fini holomorphe
% $$\overline{i^0_{r,g}}:M^{\K^0_G}(\G_r,\X_r)^*(\C)\fl M^{\K}(\G,\X)^*(\C),$$
% dont l'image est la même que celle de $\overline{i_{r,g}}$. 

\end{proposition}

\subsection{Algébricité}

On donne ici des résultats d'algébricité (sur $\C$) pour les variétés de Shimura et les compactifications de Baily-Borel.

\begin{theoreme} ([BB] th. 10.11, cf aussi [P1] 8.2)
Soient $(\G,\X)$ une donnée de Shimura pure, et $\K$ un sous-groupe compact ouvert net de $\G(\Af)$. Alors $M^{\K}(\G,\X)^*(\C)$ est l'ensemble des points complexes d'une variété algébrique complexe projective normale, qu'on notera $M^{\K}(\G,\X)^*_\C$. 
De plus, les morphismes $\overline{T}_g$ du paragraphe précédent sont algébriques.

\end{theoreme}

\begin{corollaire} ([P1] 9.24) Avec les même hypothèses,
$M^{\K}(\G,\X)(\C)$ est l'ensemble des points sur $\C$ d'une variété algébrique complexe quasi-projective lisse, notée $M^{\K}(\G,\X)_\C$, et les morphismes $T_g$ définis plus haut sont algébriques, ainsi que l'immersion ouverte
$$M^{K}(\G,\X)(\C)\fl M^{\K}(\G,\X)^*(\C).$$

La limite projective $M(\G,\X)_\C$ des $M^{\K}(\G,\X)_\C$ est un schéma séparé quasi-compact sur $\C$, sur lequel $\G(\Af)$ agit continûment (à droite), et on a pour tout $\K$
$$M^{\K}(\G,\X)_\C=\M(\G,\X)_\C/\K.$$

\end{corollaire}

Et enfin :

\begin{proposition} La stratification de $M^{\K}(\G,\X)^*(\C)$ définie plus haut est algébrique, c'est-à-dire 
que les morphismes $i_{r,g}$ sont algébriques (ce qui implique que les morphismes $\overline{i_{r,g}}$ le sont aussi).

\end{proposition}

\section{Combinatoire des strates de la compactification de Baily-Borel}

Les résultats de ce paragraphe serviront dans la section 5.2.

Soit $\K\subset\G(\Af)$ un sous-groupe ouvert compact net. 
Pour tout $I\subset\{1,\dots,q\}$ non vide, soit
$$\P_I=\bigcap_{i\in I}\P_i,$$
et
$$\N_I=\prod_{i\in I}\N_i$$
le radical unipotent de $\P_I$.
Les $\P_I$ sont les sous-groupes paraboliques standard de $\G$.
Si $r=max(I)$, on a $\QP_r\subset\P_I\subset\P_r$, donc $\P_I/\N_r$ est produit quasi-direct de $\G_r$ et $\P_{\lin,I}$, avec $\P_{\lin,I}$ un sous-groupe parabolique de $\L_{\lin,r}$.
Si $r<q$, $\G_r=\G'_r$, et $\P_I/\N_r$ est produit direct de $\G_r$ et $\P_{\lin,I}$; on note $\P'_{\lin,I}=\P_{\lin,I}$.
Si $r=q$, $\G_q$ est le centre de $\G'_q\simeq\GU(p-q)$. 
$\P_I/\N_q$ s'écrit $\P'_{\lin,I}\times\G'_q$, avec $\P'_{\lin,I}$ un sous-groupe parabolique de $\L'_{\lin,q}$,
et on a $\P_{\lin,I}=\P'_{\lin,I}\times\SU(p-q)$.

Pour tous $g\in\G(\Af)$, $r\in\{1,\dots,q\}$ et $I\subset\{1,\dots,q\}$ non vide, on note
$$\K_{g,r}=(g\K g^{-1}\cap\QP_r(\Af))/(g\K g^{-1}\cap\N_r(\Af))$$
$$\Hr_{g,I}=g\K g^{-1}\cap\P_I(\Q)\QP_{max(I)}(\Af)=g\K g^{-1}\cap\P'_{\lin,I}(\Q)\QP_{max(I)}(\Af)$$
$$\K_{g,I}=(g\K g^{-1}\cap\P_I(\Q)\QP_{max(I)}(\Af))/(g\K g^{-1}\cap\P_{\lin,I}(\Q)\N_{max(I)}(\Af)).$$
Si $r=max(I)$, $\K_{g,r}\subset\K_{g,I}$ s'identifient à des sous-groupes compacts ouverts nets de $\G_r(\Af)$.

Pour tous $g\in\G(\Af)$ et $r\in\{1,\dots,q\}$, on a défini dans 1.4.2 un morphisme holomorphe fini sur son image
$i:M^{\K_{g,r}}(\G_r,\X_r)(\C)\fl M^{\K}(\G,\X)^*(\C)$.
% $$\xymatrix{M^{\K_{g,r}}(\G_r,\X_r)(\C)\ar@{=}[d] & \\
% \G_r(\Q)\sous (\X_r\times\G_r(\Af)/\K_{g,r})\ar@{=}[d] & \\
% \P_r(\Q)\sous (\X_r\times\P_r(\Q)\QP_r(\Af)/(g\K g^{-1}\cap\QP_r(\Af)))\ar[d] & [(x,h)]\ar@{|->}[d] \\
% \G(\Q)\sous (\X^*\times \G(\Af)/\K)\ar@{=}[d] & [(x,hg)] \\
% M^{\K}(\G,\X)^*(\C) & }$$
L'image de ce morphisme est un sous-espace localement fermé de $M^{\K}(\G,\X)^*(\C)$, qui s'identifie à $M^{\K_{g,\{r\}}}(\G_r,\X_r)(\C)$.
De plus, $M^{\K_{g,r}}(\G_r,\X_r)(\C)$ et $M^{\K_{h,r}}(\G_r,\X_r)(\C)$ ont la même image dans $M^{\K}(\G,\X)^*$ \ssi $\P_r(\Q)\QP_r(\Af)g\K=\P_r(\Q)\QP_r(\Af)h\K$.
Si $I\subset\{1,\dots,q\}$ et $r=max(I)$, on a aussi un morphisme holomorphe fini sur son image
$M^{\K_{g,I}}(\G_r,\X_r)(\C)\fl M^{\K}(\G,\X)^*(\C)$, le composé du revêtement étale fini
$T_1:M^{\K_{g,I}}(\G_r,\X_r)(\C)\fl M^{\K_{g,\{r\}}}(\G_r,\X_r)(\C)$ et de l'immersion localement fermée $M^{\K_{g,\{r\}}}(\G_r,\X_r)(\C)\fl M^{\K}(\G,\X)^*(\C)$.
Tous ces morphismes se prolongent aux compactifications de Baily-Borel, et ces prolongements 
sont des morphismes finis.

\begin{notation} Soient $g,h\in\G(\Af)$ et $I\subset J\subset\{1,\dots,q\}$ non vides tels que $max(I)=max(J)=r$. 
Si $a\in\L'_{\lin,r}(\Q)\QP_r(\Af)$ est tel que $\Hr_{g,J}\subset a\Hr_{h,I}a^{-1}$, 
alors la réduction modulo $(\P_{\lin,r}\N_r)(\Af)$ de $a$, $\overline{a}\in\G_r(\Af)$, 
vérifie $\K_{g,J}\subset\overline{a}\K_{h,I}\overline{a}^{-1}$, donc on a un morphisme
$$T_{\overline{a}}:M^{\K_{g,J}}(\G_r,\X_r)(\C)\fl M^{\K_{h,I}}(\G_r,\X_r)(\C).$$
On notera aussi $T_a$ ce morphisme.

\end{notation}

On fixe $a\in \G(\Af)$, $I\subset\{1,\dots,q\}$ non vide et $s>r=max(I)$.
On considère l'ensemble $\mathcal{E}$ des diagrammes commutatifs $D_{b,q,g,h,c}$
$$\xymatrix@C=15pt{M^{\K_{hg,I\cup\{s\}}}(\G_s,\X_s)(\C)\ar@^{^{(}->}[r]\ar[d]_{T_c} & M^{\K_{g,I}}(\G_r,\X_r)^*(\C) & M^{\K_{g,I}}(\G_r,\X_r)(\C)\ar@{_{(}->}[l]\ar[d]^{T_q} & \\
M^{\K_{a,\{s\}}}(\G_s,\X_s)(\C) & & M^{\K_{b,\{r\}}}(\G_r,\X_r)(\C)\ar@{^{(}->}[r] & M^{\K}(\G,\X)^*}$$
avec $b,g\in\G(\Af)$, $q,h\in\P_r(\Q)\QP_r(\Af)$ et $c\in\P_s(\Q)\QP_s(\Af)$ tels que $g\in bq\K$ et $hg\in ca\K$.
On munit $\mathcal{E}$ de la relation d'équivalence suivante : $D_{b,q,g,h,c}\sim D_{b',q',g',h',c'}$ \ssi 
\begin{itemize}
\item[(i)] $\P_r(\Q)\QP_r(\Af)b\K=\P_r(\Q)\QP_r(\Af)b'\K$;
\item[(ii)] $\P_I(\Q)\QP_r(\Af)g\K=\P_I(\Q)\QP_r(\Af)g'\K$;
\item[(iii)] $\P_{I\cup\{s\}}(\Q)\QP_s(\Af)hg\K=\P_{I\cup\{s\}}(\Q)\QP_s(\Af)hg\K$.

\end{itemize}
La condition (i) signifie que $M^{\K_{b,\{r\}}}(\G_r,\X_r)$ et $M^{\K_{b',\{r\}}}(\G_r,\X_r)(\C)$ ont la même image dans $M^{\K}(\G,\X)^*$, la condition (ii), qui se réécrit
$$\P_I(\Q)\QP_r(\Af)q\Hr_{b,\{r\}}=\P_I(\Q)\QP_r(\Af)q'\Hr_{b,\{s\}},$$
implique que l'isomorphisme $M^{\K_{b,\{r\}}}(\G_r,\X_r)(\C)\simeq M^{\K_{b',\{r\}}}(\G_r,\X_r)(\C)$ donné par (i) se relève (et se prolonge) en un isomorphisme $T_{\gamma}:M^{\K_{g,I}}(\G_r,\X_r)^*(\C)\iso M^{\K_{g',I}}(\G_r,\X_r)^*(\C)$, avec $\gamma\in\P_I(\Q)\QP_r(\Af)$, 
et la condition (iii) que, modulo cet isomorphisme, les images de $M^{\K_{hg,I\cup\{s\}}}(\G_s,\X_s)(\C)$ et $M^{\K_{h'g',I\cup\{s\}}}(\G_s,\X_s)(\C)$ sont les mêmes.
Remarquons enfin que la condition (iii) implique les conditions (i) et (ii), 
et que, pour tout $D_{b,q,g,h,c}\in\mathcal{E}$, on a
$$D_{b,q,g,h,c}\sim D_{ca,1,ca,1,c}.$$

\begin{proposition}\label{combi_strates} Considérons l'application $\Phi$ qui à un élément $D_{b,q,g,h,c}\in\mathcal{E}$ associe la classe de $c\in\P_s(\Q)\QP_s(\Af)$ dans
$$\P_{I\cup\{s\}}(\Q)\QP_s(\Af)\sous\P_s(\Q)\QP_s(\Af)/\Hr_{a,\{s\}}.$$
Alors $\Phi$ passe au quotient par $\sim$ et donne une bijection
$$\mathcal{E}/\sim\iso \P_{I\cup\{s\}}(\Q)\QP_s(\Af)\sous\P_s(\Q)\QP_s(\Af)/\Hr_{a,\{s\}}.$$

\end{proposition}

\begin{preuve} Soient $D=D_{b,q,g,h,c},D'=D_{b',q',g',h',c'}\in\mathcal{E}$. 
Supposons que $D\sim D'$, et montrons que $\Phi(D)=\Phi(D')$. 
Soient $k,k',l\in\K$ et $\gamma\in\P_{I\cup\{s\}}(\Q)\QP_s(\Af)$ tels que $hg=cak$, $h'g'=c'ak'$ et $\gamma hgl=h'g'$.
Alors 
$$c=hgk^{-1}a^{-1}=\gamma^{-1}h'g'l^{-1}k^{-1}a^{-1}=\gamma^{-1}c'ak'l^{-1}l^{-1}a^{-1},$$
avec $\gamma^{-1}\in\P_{I\cup\{s\}}(\Q)\QP_s(\Af)$ et $ak'lk^{-1}a^{-1}={p'}^{-1}\gamma p\in (a\K a^{-1}\cap\P_s(\Q)\QP_s(\Af))=\Hr_{a,\{s\}}$,
donc $\Phi(D)=\Phi(D')$.

Supposons maintenant que $\Phi(D)=\Phi(D')$, et montrons que $D\sim D'$.
Grâce au calcul qu'on vient de faire, on peut remplacer $D$ et $D'$ par des diagrammes équivalents.
D'après la remarque qui précède l'énoncé de la proposition, on peut supposer que $q=q'=h=h'=1$, $b=g=ca$ et $b'=g'=c'a$.
$\Phi(D)=\Phi(D')$ s'écrit 
$$\P_{I\cup\{s\}}(\Q)\QP_s(\Af)c(a\K a^{-1}\cap\P_s(\Q)\QP_s(\Af))=\P_{I\cup\{s\}}(\Q)\QP_s(\Af)c'(a\K a^{-1}\cap\P_s(\Q)\QP_s(\Af)),$$
et ceci implique évidemment que
$$\P_{I\cup\{s\}}(\Q)\QP_s(\Af)caK=\P_{I\cup\{s\}}(\Q)\QP_s(\Af)c'a\K,$$
c'est-à-dire que $D\sim D'$.

Enfin, pour tout $c\in\L'_{\lin,s}(\Q)\QP_s(\Af)$, la classe de $c$ dans $\P_{I\cup\{s\}}(\Q)\QP_s(\Af)\sous\P_s(\Q)\QP_s(\Af)/\Hr_{a,\{s\}}$ est l'image par $\Phi$ de $D_{ca,1,ca,1,c}\in\mathcal{E}$, donc $\Phi$ est surjective.

\end{preuve}

\section{Modèles canoniques}

On renvoie à [P1] chapitre 11 et à [D] pour la théorie générale des modèles canoniques. On s'intéressera ici à la donnée de Shimura $(\G,\X)$ définie dans 
la section 1.2, c'est-à-dire que $\G=\GU(p,q)$.

\subsection{Variété de Shimura}

Soit $\K$ un sous-groupe ouvert compact net de $\G(\Af)$. La théorie des modèles canoniques assure l'existence d'un modèle de $M^{\K}(\G,\X)_\C$ sur un corps de nombres $F$ appelé corps reflex de $(\G,\X)$, c'est-à-dire d'une variété algébrique normale quasi-projective $M^{\K}(\G,\X)$ sur $F$ telle que $M^{\K}(\G,\X)\otimes_F\C\simeq M^{\K}(\G,\X)_\C$ (vérifiant aussi d'autres propriétés). 
De plus, tous les morphismes $T_g$ de la section 1.1 sont définis sur $F$.

\paragraph{Corps reflex et loi de réciprocité}
$$\quad$$

\begin{fait} ([D1] 2.6, [P1] 11.1) Soit $r:\Gr_{m,\C}\fl\SD_C$ le morphisme qui à $z\in\C^\times$ associe 
$$(z,1)=1\otimes\frac {z+1} {2} +i\otimes \frac {z-1} {2i}\in\SD(\C)=(\C\otimes_\R\C)^\times.$$
 Le corps reflex de $(\G,\X)$, qu'on notera $F$,  est le corps de définition de la classe de conjugaison de
$$h_0 r:\Gr_{m,\C}\fl\G_\C.$$

\end{fait}

\begin{lemme} Le corps de reflex de $(\G,\X)$ est $E$ si $p>q$, $\Q$ si $p=q$.

\end{lemme}

\begin{preuve} On écrit comme avant $E=\Q[i\sqrt{d}]$, avec $d$ un entier positif sans facteur carré. 
On a, si $z\in\C^\times$,
$$h_0r(z)=\left(\begin{array}{ccc}1\otimes \frac{z+1} {2} I_q & 0 & i\sqrt{d}\otimes\frac {z-1} {2i\sqrt{d}}J_q \\
0 & \left(1\otimes\frac {z+1} {2}+i\sqrt{d}\otimes\frac {z-1} {2i\sqrt{d}}\right)I_{p-q} & 0\\
i\sqrt{d}\otimes\frac {z-1} {2i\sqrt{d}}J_q & 0 & 1\otimes\frac {z+1} {2} I_q\end{array}\right)\in\G(\C).$$

On cherche le plus petit sous-corps $F$ de $\C$ tel que $h_0r(\Q^\times)$ soit conjugué dans $\G(\C)$ à un sous-groupe de $\G(\Q)\subset\GL_n(E\otimes_\Q\Q)$. 

Déjà, on a $h_0r(\Q^\times)\subset\GL_n(E\otimes_\Q E)$, donc $F\subset E$.

De plus, pour tout $z\in\Q^\times$, on a
$$Tr(h_0r(z))=n\left(1\otimes \frac {z+1} {2}\right)+(p-q)\left(i\sqrt{d}\otimes\frac {z-1} {2i\sqrt{d}}\right),$$
et $(z-1)/2i\sqrt{d}\not\in\Q$ si $z\not =1$, donc, si $p>q$, $h_0r(\Q^\times)$ ne peut pas être conjugué dans $\GL_n(E\otimes_\Q\C)$ à un sous-groupe de $\GL_n(E\otimes_\Q\Q)$, donc $F\not =\Q$, et finalement $F=E$.

Supposons maintenant $p=q$. On a 
$$h_0r(z)=1\otimes\frac {z+1} {2}I_n+i\sqrt{d}\otimes\frac {z-1} {2i\sqrt{d}} J_n.$$
$J_n$ est conjugué dans $\GL_n(\C)$ à
$$J'=\left(\begin{array}{cc}0 & iJ_q \\ -iJ_q & 0\end{array}\right).$$
Comme le morphisme $\GU(p,q)(\C)\fl\GL_n(\C)$, $1\otimes X+i\sqrt{d}\otimes Y\fle X+i\sqrt{d}Y$, induit un isomorphisme $\U(p,q)(\C)\iso\GL_n(\C)$, il existe $g\in\U(p,q)(\C)$ tel que l'image de $J''=g(i\sqrt{d}\otimes (1/i\sqrt{d})J_n)g^{-1}$ dans $\GL_n(\C)$ soit $J'$. De plus, on a $c(J'')=c(i\sqrt{d}\otimes (1/i\sqrt{d})J_n)=-1$. Soit
$$J'''=\left(\begin{array}{cc}0 & 1\otimes i J_q \\ -1\otimes i J_q & 0\end{array}\right).$$
$J'''\in\GU(p,q)(\C)$, l'image de $J'''$ dans $\GL_n(\C)$ est $J'$, et $c(J''')=-1$, donc $J''=J'''$. On a donc montré que pour tout $z\in\C^\times$,
$$gh_0r(z)g^{-1}=\left(\begin{array}{cc}1\otimes \frac {1+z} {2} I_q & 1\otimes\frac {z-1} {2} J_q \\ 1\otimes\frac {1-z} {2} J_q & 1\otimes \frac {1+z} {2} I_q\end{array}\right),$$
donc $h_0r(\Q^\times)$ est conjugué dans $\GU(p,q)(\C)$ à un sous-groupe de $\GL_n(\Q\otimes_\Q\C)$, et le corps reflex est $\Q$.
\end{preuve}

Comme $M^{\K}(\G,\X)(\C)$ est l'ensemble des points complexes d'une variété définie sur $F$, le groupe $Gal(\overline{F}/F)$ agit sur $\pi_0(M^{\K}(\G,\X)(\C))$. On va calculer cette action. 

Soit $\H=\G/\G^{der}$. On note $\kappa:\G\fl\H$ la projection.

\begin{proposition}([D] 3.3.2) L'application, induite par $\kappa$,
$$\pi_0(M^{\K}(\G,\X)(\C))\fl\pi_0(\H(\Q)\sous\H(\Ad)/\kappa(\K))/\kappa(\Ar\K_\infty)$$
est bijective.

\end{proposition}

L'action de $Gal(\overline{F}/F)$ sur $\pi_0(\H(\Q)\sous\H(\Ad)/\kappa(\K))/\kappa(\Ar\K_\infty)$ se factorise par $Gal(\overline{F}^{ab}/F)$ ([D] 2.4). On identifie $Gal(\overline{F}^{ab}/F)$ à $\pi_0(F^\times\sous\Ad_F^\times)$ par l'isomorphisme du corps de classes, avec la convention suivante ([P1] 11.3) : Si $v$ est une place non archimédienne de $F$, $\varpi_v$ est une uniformisante de $F_v$ et $F'$ est une extension abélienne de $F$ non ramifiée en $v$, alors le morphisme
$$\pi_0(F^\times\sous\Ad_F^\times)\iso Gal(\overline{F}^{ab}/F)\fl Gal(F'/F)$$
envoie la composante de l'idèle $(\dots,1,\varpi_v,1,\dots)$ sur le Frobenius arithmétique en $v$.

\begin{proposition}([D] 2.4, 2.5, 2.6) L'action de $Gal(\overline{F}^{ab}/F)\simeq \pi_0(F^\times\sous\Ad_F^\times)$ sur\newline
 $\pi_0(\H(\Q)\sous\H(\Ad)/\kappa(\K))/\kappa(\Ar\K_\infty)$ provient d'un morphisme de groupes algébriques
$$\lambda:Res_{F/\Q}\Gm\fl\H,$$
appelé morphisme de réciprocité, et donné par la règle suivante : si $r':\Gr_{m,F}\fl\H_F$ est le morphisme $\kappa h_0 r$ (quitte à changer de point base, ce morphisme est défini sur $F$), $\lambda$ est le composé de
$$Res_{F/\Q}(r'):Res_{F/\Q}\Gm\fl Res_{F/\Q}\H_F$$
et de la norme de $\H$
$$Res_{F/\Q}\H_F\fl\H.$$

\end{proposition}

% Supposons $n$ impair, $n=2k+1$. On a alors $F=E$ (car $p>q$), et on peut identifier $\H$ à $Res_{E/\Q}\Gm$, avec $\kappa=c^k/det$.

% \begin{lemme} 
% Supposons $n$ impair, $n=2k+1$. 
% On a alors $F=E$ (car $p>q$), et on peut identifier $\H$ à $Res_{E/\Q}\Gm$.
% $\lambda:Res_{E/\Q}\Gm\fl Res_{E/\Q}\Gm$ est alors donné par la formule suivante : si $z\in E^\times$,
% $$\lambda(z)=z^{k-p}\overline{z}^{k-q}.$$

\begin{lemme}
On a
$$\H\simeq\{(x,\lambda),x\overline{x}=\lambda^n\}\subset Res_{E/\Q}\Gm\times\Gm,$$
et
$$\lambda(z) = (z^p\overline{z}^q,z\overline{z}).$$

\end{lemme}

\begin{corollaire} On note $\lambda':Res_{E/\Q}\Gm\fl\G$ l'application qui envoie $z$ sur
$$\left(\begin{array}{ccc}z\overline{z}I_q & 0 & 0 \\ 0 & zI_{p-q} & 0 \\ 0 & 0 & I_q\end{array}\right).$$

Comme $\lambda=\kappa\lambda'$, l'action de $Gal(\overline{F}^{ab}/F)\simeq\pi_0(F^\times\sous\Ad_F^\times)$ sur $\pi_0(M^{\K}(\G,\X)(\C))$ est induite par l'action de $Res_{E/\Q}\Gm$ sur $\G$ donnée par la multiplication à gauche par $\lambda'$.

\end{corollaire}

\paragraph{Le système projectif}

$$\quad$$

On a dit que les $T_g$ étaient définis sur $F$. En particulier, pour $\K'\subset\K$, on a 
$$T_1:M^{\K'}(\G,\X)\fl M^{\K}(\G,\X).$$
Le système projectif des $M^{\K}(\G,\X)$ est donc défini sur $F$, et on note $M(\G,\X)$ sa limite projective. C'est un $F$-schéma séparé quasi-compact, sur lequel $\G(\Af)$ agit continûment.

\begin{fait} ([P2] 3.4) Si $\K'\subset\K$ un sous-groupe ouvert distingué, alors
$$T_1:M^{\K'}(\G,\X)\fl M^{\K}(\G,\X)$$
est un recouvrement étale galoisien fini de groupe $\K/\K'$.

En conséquence,
$$M(\G,\X)\fl M^{\K}(\G,\X)$$
est un recouvrement étale galoisien profini de groupe $\K$.

\end{fait} 

\subsection{Compactification de Baily-Borel}

\begin{theoreme} ([P1] 12.3, [P2] 3.7.2) 
\begin{itemize}
\item[(1)] $M^{\K}(\G,\X)^*_\C$ a un modèle $M^{\K}(\G,\X)^*$ sur $F$, les $\overline{T}_g$ sont définis sur ces modèles et on a la propriété suivante : 

\item[(2)] Soit $(\QP,\Y)$ une composante rationnelle de bord de $(\G,\X)$ et $g\in\G(\Af)$. On note $\N$ le radical unipotent de $\QP$, $(\G',\Ze)=(\QP,\Y)/\N$, 
$\K_Q=\QP(\Af)\cap g\K g^{-1}$, $\K_G$ l'image de $\K_Q$ dans $\G'(\Af)$. 
Alors le corps reflex de $(\G',\Ze)$ est celui de $(\G,\X)$, c'est-à-dire $F$, et le morphisme (fini sur son image)
$$M^{\K_G}(\G',\Ze)(\C)\fl M^{\K}(\G,\X)^*(\C)$$
de 1.4.2 vient d'un unique morphisme des modèles canoniques.

\end{itemize}
\end{theoreme}

\begin{proposition}([P2] 3.8)
Soit $\K'\subset\K$ un sous-groupe ouvert distingué. $K$ agit à droite sur $M^{\K'}(\G,\X)$ (par les morphismes $T_g$), et le quotient est
$$M^{\K'}(\G,\X)^*/\K=M^{\K}(\G,\X)^*.$$ 

% Soient $(\QP,\Y)$ une composante de bord, $g\in\G(\Af)$, $\P$ le parabolique maximal associé; on définit $(\G',\Ze)$, $\K_P$, $\K_G$ comme avant, et on pose $\Hr_P=g\K g^{-1}\cap \P(\Q)\QP(\Af)$, $\K_N=g\K g^{-1}\cap\N(\Af)$. On note $\K'_P$,$\K'_G$, $\Hr'_P$ les analogues de ces groupes pour $\K'$. Le morphisme
% $$T_1:M^{\K'}(\G,\X)^*/g^{-1}\Hr_P g\simeq M^{\K'g^{-1}\Hr_Pg}(\G,\X)^*\fl M^{\K}(\G,\X)^*$$
% est étale au voisinage de la strate définie par $(\QP,\Y)$ et $g$. D'autre part, le morphisme entre les strates
% $$M^{\K'_G}(\G',\Ze)\fl M^{\K_G}(\G',\Ze)$$
% est un recouvrement étale galoisien de groupe $\K_G/\K'_G\simeq\K_P/\K'_P\K_N$.

\end{proposition}

\begin{corollaire} On note $M(\G,\X)^*$ la limite projective sur $\K$ des $M^{\K}(\G,\X)^*$. Le morphisme
$$M(\G,\X)^*\fl M(\G,\X)^*/\K\simeq M^{\K}(\G,\X)^*$$
est un $K$-recouvrement (au sens de [P2] 1.7). La limite projective des morphismes
$$M^{\K_G}(\G',\Ze)\fl M^{\K}(\G,\X) *$$
est une immersion (localement fermée)
$$M(\G',\Ze)\fl M(\G,\X)^*.$$

\end{corollaire}

\newpage

\chapter{Systèmes de coefficients et leurs prolongements à la compactification de Baily-Borel, d'après Pink}

\section{Variété de Shimura}

\subsection{Systèmes locaux sur $\C$ à coefficients dans $\Q_\ell$}

\subsubsection{Espaces localement symétriques}

Soient $\G$ un groupe algébrique linéaire sur $\Q$ et $ X$ un espace symétrique de $\G$. On note comme avant $\S$ le centre déployé de $\G$ et $A=\S(\R)^{\circ}$.

Soit $\Gamma$ un sous-groupe arithmétique net de $\G(\Q)$. Pour tout $\Gamma$-module à gauche discret $M$, on pose
$$\F^{\Gamma}(M)=\Gamma\sous (M\times X),$$
où, si $\gamma\in\Gamma$ et $(m,x)\in M\times X$,
$$\gamma.(m,x)=(\gamma m,\gamma x).$$
La flèche évidente $\F^{\Gamma}(M)\fl\Gamma\sous X$ fait de $\F^{\Gamma}(M)$ un système local (de groupes abéliens) sur l'espace localement symétrique $\Gamma\sous X$ (qui est une variété analytique réelle).

Il est clair que $\F^{\Gamma}$ est un foncteur exact et fidèle de la catégorie $Mod_{\Gamma}$ des $\Gamma$-modules à gauche dans celle des systèmes locaux sur $\Gamma\sous X$. On notera donc $\F^{\Gamma}$ son foncteur dérivé, défini sur $D^b(Mod_{\Gamma})$. L'isomorphisme canonique $H^0(\Gamma\sous X,\F^{\Gamma}(M))=M^{\Gamma}$ s'étend en un isomorphisme fonctoriel
$$R\Gamma(\Gamma\sous X,\quad)\F^{\Gamma}\simeq R\Gamma(\Gamma,\quad).$$

Soient $g\in\G(\Q)$ et $\Gamma,\Gamma'$ deux sous-groupes arithmétiques nets de $\G(\Q)$ tels que $\Gamma'\subset g\Gamma g^{-1}$.
On a défini dans 1.1 un morphisme analytique réel
$$\left\{\begin{array}{rcl}T_g:\Gamma'\sous X & \fl & \Gamma\sous X \\
\Gamma' x & \fle & \Gamma g^{-1}x\end{array}\right..$$
On note $\theta$ le morphisme $\gamma\fle g^{-1}\gamma g$ de $\Gamma'$ dans $\Gamma$. Le foncteur $\theta^*:Mod_{\Gamma}\fl Mod_{\Gamma'}$ de restriction des scalaires via $\theta$ est exact, et on note encore $\theta^*$ le foncteur dérivé $D^b(Mod_{\Gamma})\fl D^b(Mod_{\Gamma'})$.

\begin{fait}\begin{itemize}
\item[(i)] On a un isomorphisme canonique
$$ \F^{\Gamma'}\theta^*\simeq T_g^*\F^{\Gamma},$$
donné, pour un $\Gamma$-module $M$, par
$$\begin{array}{rcl}\Gamma'\sous (M\times X) & \fl & (\Gamma\sous (M\times X))\times_{\Gamma\sous X}\Gamma'\sous X \\
(m,x) & \fle & (m,g^{-1}x,x)\end{array}.$$

\item[(ii)] Si $M$ est un $\G(\Q)$-module et si $\Gamma=\Gamma'$ (donc $g$ normalise $\Gamma$)
on a un isomorphisme de $\Gamma$-modules $M\iso\theta^*M,m\fle g^{-1}m$,
et le composé
$$R\Gamma(\Gamma,M)=R\Gamma(\Gamma\sous X,\F^{\Gamma}(M))\stackrel{CB}{\fl}R\Gamma(\Gamma\sous X,T_g^*\F^{\Gamma}(M))\simeq R\Gamma(\Gamma\sous X,\F^{\Gamma}(M))\simeq R\Gamma(\Gamma,M)$$
est l'endomorphisme de $R\Gamma(\Gamma,M)$ provenant de l'endomorphisme $\gamma\fle g^{-1}\gamma g$ de $\Gamma$.

\end{itemize}
\end{fait}

\subsubsection{Variétés de Shimura}

On fixe un nombre premier $\ell$.
Supposons que $(\G,\X)$ est une donnée de Shimura (pure ou mixte). Notons $Rep_{\G(\Q_\ell)}$ la catégorie des représentations algébriques de $\G(\Q_\ell)$ dans un $\Q_\ell$-espace vectoriel de dimension finie.

Soit $\K$ un sous-groupe compact ouvert net de $\G(\Af)$. 
Si $V\in Rep_{\G(\Q_\ell)}$, on pose
$$\F^{\K}(V)=\G(\Q)\sous V\times\X\times (\G(\Af)/\K),$$
où $\G(\Q)$ agit sur $V\times\X\times (\G(\Af)/\K)$ de la manière suivante : si $\gamma\in\G(\Q)$, $v\in V$, $x\in\X$ et $g\in\G(\Af)$,
$$\gamma.(v,x,gK)=(\gamma.v,\gamma.x,\gamma gK).$$
$\F^{\K}(V)$ est un système local de $\Q_\ell$-espaces vectoriels sur $M^{\K}(\G,\X)(\C)$.

Comme plus haut, il est clair que $\F^{\K}$ est un foncteur exact et fidèle sur $Rep_{\G(\Q_\ell)}$. 
On notera de la même façon le foncteur dérivé.

Soient $g\in\G(\Af)$ et $K'\subset gKg^{-1}$ un sous-groupe ouvert, et soit comme avant
$$\begin{array}{rcl}T_g : M^{\K'}(\G,\X)(\C) & \fl & M^{\K}(\G,\X)(\C) \\
                          (x,hK') & \fle & (x,hgK) \end{array}.$$
Pour tout $V\in D^b(Rep_{\G(\Q_\ell)})$, on a un isomorphisme
$$\F^{\K'}(V)\simeq T_g^* \F^{\K}(V),$$ 
donné par
$$\begin{array}{rcl} \G(\Q)\sous V\times\X\times (\G(\Af)/\K') & \fl & (\G(\Q)\sous V\times \X\times (\G(\Af)/\K))\times_{M^{\K}(\G,\X)(\C)}M^{\K'}(\G,\X)(\C)\\
(v,(x,hK')) & \fle & (v,(x,hgK),(x,hK')) \end{array}.$$

\begin{remarque} On retrouve sur les composantes de $M^{\K}(\G,\X)(\C)$ les systèmes locaux définis plus haut : Fixons $g\in\G(\Af)$ et posons $\Gamma=g\K g^{-1}\cap\G(\Q)$. 
On a une immersion holomorphe
$$\begin{array}{rcl}i:\Gamma\sous\X & \fl & M^{\K}(\G,\X)(\C) \\
\Gamma x & \fle & \lbrack (x,g)\rbrack \end{array}$$
dont l'image est une composante connexe de $M^{\K}(\G,\X)(\C)$, et qui est un isomorphisme sur son image. 
Pour tout $V\in D^b(Rep_{\G(\Q_\ell)})$, $i^*\F^{\K}(V)$ est canoniquement isomorphe à $\F^{\Gamma}(M)$, où $M$ est $V$, avec la structure de $\Gamma$-module donnée par l'inclusion $\Gamma\subset\G(\Q)\subset\G(\Q_\ell)$.

\end{remarque}

\subsection{Systèmes locaux sur $\C$ à coefficients dans $\Z_\ell$}

Soient $(\G,\X)$ une donnée de Shimura, $\K$ un sous-groupe ouvert compact net de $\G(\Af)$ et $V\in Rep_{\G(\Q_\ell)}$.
On voudrait associer à $V$ un faisceau $\F^{\K}(V)(\Z_\ell)$ de $\Z_\ell$-modules sur $M^{\K}(\G,\X)(\C)$. 
On va rappeler la méthode de [Ln], p34.

On fait agir $\G(\Af)$ sur $V$ via la projection $\G(\Af)\fl\G(\Q_\ell)$, et
on choisit un $\Z_\ell$-réseau stable par $\K$ de $V$, qu'on note $V(\Z_\ell)$. 
Soit
$$\mathcal{E}=\coprod_{g\K\in\G(\Af)/K} (g.V(\Z_\ell)\times\X\times (g\K/\K)).$$
On a une application continue surjective $\mathcal{E}\fl\X\times (\G(\Af)/\K)$ 
qui envoie $(g.v,x,gK)$ sur $(x,gK)$. 
Cette application est $\G(\Q)$-équivariante, pour l'action suivante de $\G(\Q)$ sur $\mathcal{E}$ :
$$\gamma.(g.v,x,g\K)=(\gamma g.v,\gamma.x,(\gamma g)\K).$$
On pose
$$\F^{\K}(V)(\Z_\ell)=\G(\Q)\sous\mathcal{E}.$$
$\F^{\K}(V)(\Z_\ell)\fl \G(\Q)\sous\X\times (\G(\Af)/\K)=M^{\K}(\G,\X)(\C)$ est un faisceau en $\Z_\ell$-modules localement libre, et on a 
$$\F^{\K}(V)(\Z_\ell)\otimes_{\Z_\ell}\Q_\ell=\F^{\K}(V).$$
De plus, la classe d'isomorphisme de $\F^{\K}(V)(\Z_\ell)$ est indépendante du choix de $V(\Z_\ell)$.

Si $g\in\G(\Af)$ et $\K'\subset g\K g^{-1}$, on a comme avant un isomorphisme canonique
$$\F^{\K'}(V)(\Z_\ell)\simeq T_g^*\F^{\K}(V)(\Z_\ell).$$
Pour $m\in\Nat^*$, on pose
$$\F^{\K}(V)(\Z/\ell^m\Z)=\F^{\K}(V)(\Z_\ell)\otimes_{\Z_\ell}\Z/\ell^m\Z.$$ 
C'est un faisceau abélien de $\ell^m$-torsion sur $M^{\K}(\G,\X)(\C)$. 
Nous montrerons plus loin que ce faisceau provient d'un faisceau étale sur $M^{\K}(\G,\X)$.

\subsection{Rappels sur la construction de Pink}

On va introduire rapidement certains faisceaux étales, qui sont construits et étudiés dans les deux premières parties de [P2].

\begin{rappel}([P2] 1.1.2) 
Un groupe topologique $\Gamma$ est dit \emph{de type pro-}$FP_\infty$ s'il existe une famille de sous-groupes ouverts distingués $\Delta$ de $\Gamma$ telle que :
\begin{itemize}
\item les $\Delta$ sont profinis;
\item pour tout $\Delta$, le groupe discret $\Gamma/\Delta$ est de type $FP_\infty$;
\item $\Gamma=\limpd \Gamma/\Delta$.
\end{itemize}
Rappelons qu'un groupe discret $\Gamma$ est dit de type $FP_\infty$ 
si le $\Z[\Gamma]$-module trivial $\Z$ a une résolution 
$\dots\fl L_1\fl L_0\fl\Z$ par des $\Z[\Gamma]$-modules libres de type fini.

\end{rappel}

\begin{exemple}([P2] 1.1.3)
Soient $\G$ un groupe linéaire algébrique sur $\Q$, $\Gamma$ un sous-groupe arithmétique 
de $\G(\Q)$ et $\K$ un sous-groupe compact de $\G(\Af)$ normalisé par $\Gamma$. 
Alors $\Gamma$ est de type $FP_\infty$, $\K$ est profini (donc de type pro-$FP_\infty$) et $\Gamma\K$ est de type pro-$FP_\infty$.
\end{exemple}

\begin{notation}([P2] 1.1,1.2)
\begin{itemize}
 \item[$\bullet$] Si $\Gamma$ est un groupe topologique, on note $Mod_\Gamma$ la catégorie abélienne des $\Gamma$-modules à gauches discrets sur lesquels l'action de $\Gamma$ est continue (c'est-à-dire tels que le stabilisateur d'un point soit un sous-groupe ouvert de $\Gamma$). On note $Mod'_\Gamma$ la sous-catégorie pleine de $Mod_\Gamma$ engendrée 
par les limites inductives de $\Gamma$-modules de type fini sur $\Z$.
\item[$\bullet$] Tous les schémas et morphismes de schémas sont supposés quasi-séparés. 
Si $X$ est un schéma, on note $\acute{E}t_X$ la catégorie des faisceaux abéliens étales sur $X$.

\end{itemize}
\end{notation}

\begin{lemme}([P2] 1.2) Soient $X$ un schéma, $\Gamma$ un groupe de type pro-$FP_\infty$ et $\F$ un faisceau étale de $\Gamma$-modules à gauche sur $X$.  Pour tout sous-groupe $\Delta$ de $\Gamma$, on note $\F^\Delta$ le sous-faisceau de $\F$ (des $\Delta$-invariants) défini par $\F^\Delta(U)=\F(U)^\Delta$ pour tout $U\fl X$ étale. Les conditions suivantes sont équivalentes :
\begin{itemize}
\item[(1)] toutes les fibres de $\F$ sont dans $Mod_\Gamma$;
\item[(2)] $\F=\limid \F^\Delta$, où $\Delta$ parcourt l'ensemble des sous-groupes ouverts distingués de $\Gamma$;
\item[(3)] (si $X$ est quasi-compact) pour tout $U\fl X$ étale de présentation finie, $\F(U)$ est dans $Mod_\Gamma$.
\end{itemize}

\end{lemme}

\begin{notation}([P2] 1.2) On note $\acute{E}t_{X,\Gamma}$ la catégorie des faisceaux étales de $\Gamma$-modules à gauches sur $X$ satisfaisant les conditions du lemme précédent. (Si $X$ est le spectre d'un corps séparablement clos, on peut identifier $\acute{E}t_{X,\Gamma}$ et $Mod_\Gamma$.)

\end{notation}

\begin{definition} Soit $X$ un schéma sur lequel un groupe discret $\Gamma$ agit à droite fidèlement. On dit que l'action de $\Gamma$ est \emph{propre} s'il existe un recouvrement de $X$ par des ouverts affines tel que pour tout ouvert $U$ de ce recouvrement, 
$$ card(\{\gamma\in\Gamma,U.\gamma\cap U\not=\emptyset\})<\infty.$$

\end{definition}

\begin{definition}\label{def_recouvrements}([P2] 1.7) Soit $X$ un schéma muni d'une action à droite d'un groupe $\Gamma$ de type pro-$FP_\infty$. On suppose qu'il existe une famille de sous-groupes ouverts distingués $\Delta$ de $\Gamma$ telle que :
\begin{itemize}
\item pour tout $\Delta$, il existe un recouvrement ouvert affine de $X$ dont chaque ouvert est invariant par $\Delta$; en particulier, le quotient catégorique $X/\Delta$ existe;
\item $X=\limpd X/\Delta$;
\item pour tout $\Delta$, si $\Delta'$ est le noyau de l'action de $\Gamma$ sur $X/\Delta$, alors $\Gamma/\Delta'$ agit proprement sur $X/\Delta$.
\end{itemize}

Alors le quotient géométrique de $X/\Gamma$ existe. On dit que $X\fl Y$ est un $\Gamma$\emph{-recouvrement}.

\end{definition}

\begin{notation} Si $X$ est un schéma et $A$ est un groupe abélien, on note $\underline{A}_X$ ou simplement $\underline{A}$ le faisceau étale constant de fibre $A$ sur $X$.

\end{notation}

\begin{definition}([P2] 1.8) Soit $X\stackrel {\varphi}{\fl} Y$ un $\Gamma$-recouvrement. Pour un sous-groupe ouvert distingué assez petit $\Delta$ de $\Gamma$, on a un diagramme commutatif de morphismes $\Gamma$-équivariants
$$\xymatrix{X\ar[rr]\ar[rd]_\varphi & & X/\Delta\ar[dl]^{\varphi_\Delta} \\
 & Y & }.$$

Pour un tel $\Delta$, le faisceau $\varphi_{\Delta*}\underline{\Z}$ sur $Y$ est muni d'une 
action à gauche de $\Gamma$ : $\gamma\in\Gamma$ envoie une section 
$s\in\varphi_{\Delta *}\underline{\Z}(U)\simeq\Z^{\pi_0(\varphi_\Delta^{-1}(U))}$ 
sur la section $\gamma.s$ qui à $V$ dans $\pi_0(\varphi_\Delta^{-1}(U))$ 
associe $s(V.\gamma)$. Les $\varphi_{\Delta *}\underline{\Z}$ forment un système inductif 
de $\acute{E}t_{Y,\Gamma}$, dont la limite est aussi dans $\acute{E}t_{Y,\Gamma}$.

On définit alors un foncteur
$$\left\{\begin{array}{rcl}\lambda_{\Gamma,\varphi} : Mod_\Gamma & \fl & \acute{E}t_{Y,\Gamma} \\
M & \fle & \underline{M}\otimes\limid \varphi_{\Delta *}\underline{\Z}\end{array}\right.,$$

où $\Gamma$ agit sur les deux facteurs. Ce foncteur est exact et commute aux limites inductives. 

\end{definition}

\begin{lemme}\label{fibres_lambda} ([P2] 1.8) On se place dans la situation de la définition ci-dessus. Soit $M$ un objet de $Mod_\Gamma$; on va décrire les fibres de $\lambda_{\Gamma,\varphi}(M)$. Soient $x$ un point géométrique de $X$, $y$ l'image de $x$ dans $Y$, $\Gamma_1\subset\Gamma$ le stabilisateur de $x$. Le choix de $x$ détermine une bijection $\Gamma$-équivariante $\varphi^{-1}(x)\simeq\Gamma_1\sous\Gamma$, d'où un isomorphisme de $\Gamma$-modules
$$\lambda_{\Gamma,\varphi}(M)_y\simeq M\otimes Ind^\Gamma_{\Gamma_1}\Z,$$
où $\Gamma_1$ agit trivialement sur $\Z$.

\end{lemme}

\begin{proposition}([P2] 1.9.3) Soient $X\fl Y$ un $\Gamma$-recouvrement et $1\fl\Gamma'\fl\Gamma\fl\Gamma''\fl 1$ une suite exacte de groupes de type pro-$FP_\infty$ tels que $\Gamma'$ agisse trivialement sur $X$. 
Alors on a un isomorphisme de foncteurs sur $D^+(Mod'_{\Gamma})$
$$R(\Gamma',\quad)\circ\lambda_{\Gamma,\varphi}\simeq\lambda_{\Gamma'',\varphi}\circ R(\Gamma',\quad).$$

\end{proposition}

\begin{propdef}([P2] 1.10.1) Soit $\varphi:X\fl Y$ un recouvrement étale galoisien de groupe de Galois $\Gamma$. 
Alors $\lambda_{\Gamma,\varphi}(M)$ est $R\Gamma(\Gamma,\quad)$-acyclique pour tout objet $M$ de $Mod_\Gamma'$. On obtient donc un foncteur exact
$$\left\{\begin{array}{rcl}\mu_{\Gamma,\varphi}:Mod'_\Gamma & \fl & \acute{E}t_X \\
M & \fle & \Gamma(\Gamma,\lambda_{\Gamma,\varphi}(M))\end{array}\right. .$$

Pour tout $U\fl X$ étale de présentation finie, on a 
$$\mu_{\Gamma,\varphi}(M)(U)=\limid \{s:\pi_0(\varphi_\Delta^{-1}(U))\fl M,\forall\gamma\in\Gamma,\forall V\in\pi_0(\varphi_\Delta^{-1}(U)),s(V.\gamma)=\gamma^{-1}s(V)\}.$$

\end{propdef}

\begin{proposition}\label{fibres_faisceaux}([P2] 1.10.4) Soient $\varphi:X\fl Y$ comme avant, 
et $M$ un objet de $Mod'_\Gamma$. 
Soit $y_0=Spec(k)$ un point de $Y$, $\overline{k}$ une clôture séparable de $k$, $y$
le point géométrique de $Y$ correspondant,
$x$ un point géométrique de $X$ au-dessus de $y$. 
Pour tout $\sigma\in Gal(\overline{k}/k)$, on note $\psi(\sigma)$ l'unique élément de $\Gamma$ tel que $\sigma.x=x.\psi(\sigma)$. Alors 
$$\psi:Gal(\overline{k}/k)\fl\Gamma$$
est un morphisme continu, et la fibre de $\mu_{\Gamma,\varphi}(M)$ en $y$ est isomorphe à $M$ 
avec l'action de $Gal(\overline{k}/k)$ donnée par $\sigma.m=\psi(\sigma).m$. 

\end{proposition}

\begin{proposition}\label{prop:Shapiro} Soient $\til{X}$ un schéma muni d'une action à droite
d'un groupe $\Gamma$ de type pro-$FP_\infty$
et $\Gamma'$ un sous-groupe fermé de $\Gamma$.
On note $\varphi:\til{X}\fl\til{X}/\Gamma=X$,
$\varphi':\til{X}\fl\til{X}/\Gamma'=X'$ et $f:X'\fl X$
les morphismes évidents.
On a un morphisme canonique de foncteurs (provenant de morphismes de changement de base, cf [P1] 1.11.3)
$$Res_{\Gamma'}^{\Gamma}f^*\lambda_{\Gamma,\varphi}\fl\lambda_{\Gamma',\varphi'}Res_{\Gamma'}^{\Gamma}.$$
On suppose que le stabilisateur de tout point géométrique de $\til{X}$ est inclus dans $\Gamma'$.
Alors:
\begin{itemize}
\item[(i)] Le morphisme canonique
$$\begin{array}{rcl}f^*R\Gamma(\Gamma,\quad)\lambda_{\Gamma,\varphi} & \fl & f^*R\Gamma(\Gamma',\quad)Res_{\Gamma'}^{\Gamma}\lambda_{\Gamma,\varphi}=R\Gamma(\Gamma',\quad)f^*Res_{\Gamma'}^{\Gamma}\lambda_{\Gamma,\varphi} \\
& \fl & R\Gamma(\Gamma',\quad)\lambda_{\Gamma',\varphi'}Res_{\Gamma'}^{\Gamma}\end{array}$$
(cf [P1] 1.11.4)
est un isomorphisme de foncteurs sur $D^+(Mod'_{\Gamma})$.
\item[(ii)] Si de plus $\Gamma'$ est d'indice fini dans $\Gamma$, alors le morphisme canonique (déduit du morphisme de (i) par adjonction)
$$R\Gamma(\Gamma,\quad)\lambda_{\Gamma,\varphi}Ind_{\Gamma'}^{\Gamma}\fl f_*R\Gamma(\Gamma',\quad)\lambda_{\Gamma',\varphi'}$$
est un isomorphisme de foncteurs sur $D^+(Mod'_{\Gamma'})$.

\end{itemize}
\end{proposition}

\begin{preuve} (i) est [P1] 1.11.5.

(ii) se prouve de la même façon.
Soit $\til{x}$ un point géométrique de $\til{X}$.
Notons $x$ le point géométrique de $X$ image de $\til{x}$
et $\Gamma_1$ le stabilisateur de $\til{x}$ dans $\Gamma$.
On a $\Gamma_1\subset\Gamma'$ par hypothèse.
On fixe un système de représentants $(\delta_i)_{i\in I}$ de $\Gamma'\sous\Gamma$.
Soit $M\in Mod'_\Gamma$.
D'après le lemme \ref{fibres_lambda}, le choix de $\til{x}$ et de $(\delta_i)$ détermine des isomorphismes
$$\left(R\Gamma\left(\Gamma,\lambda_{\Gamma,\varphi}Ind_{\Gamma'}^{\Gamma}M\right)\right)_x\simeq Ind_{\Gamma'}^{\Gamma}M\otimes Ind_{\Gamma_1}^{\Gamma}\Z$$
$$\left(Rf_*R\Gamma\left(\Gamma',\lambda_{\Gamma',\varphi'}M\right)\right)_x\simeq
\bigoplus_{i\in I}\left(R\Gamma\left(\Gamma',\lambda_{\Gamma',\varphi'}(M)\right)\right)_{\varphi'(\til{x}\delta_i)}\simeq\bigoplus_{i\in I}M\otimes Ind_{\Gamma_1}^{\Gamma'}\Z,$$
et le morphisme 
$$u:\left(R\Gamma\left(\Gamma,\lambda_{\Gamma,\varphi}Ind_{\Gamma'}^{\Gamma}M\right)\right)_x\fl \left(Rf_*R\Gamma\left(\Gamma',\lambda_{\Gamma',\varphi'}M\right)\right)_x$$
correspond au composé du morphisme fonctoriel
$R\Gamma(\Gamma,\quad)\fl R\Gamma(\Gamma',\quad)Res^{\Gamma}_{\Gamma'}$ de [P1] 1.11.1
et du morphisme 
$$\begin{array}{rcl} v: Ind_{\Gamma'}^{\Gamma}M\otimes Ind_{\Gamma_1}^{\Gamma}\Z & \fl & \bigoplus_{i\in I}M\otimes Ind_{\Gamma_1}^{\Gamma'}\Z \\
f\otimes g & \fle & \sum_{i\in I}f(\delta_i)\otimes g_{|\Gamma'}\end{array}$$
où $f\in Ind_{\Gamma'}^{\Gamma}M$ est une fonction $f:\Gamma\fl M$ et $g\in Ind_{\Gamma_1}^{\Gamma}\Z$ est une fonction $\Gamma\fl \Z$.

Comme tous les foncteurs considérés commutent aux limites inductives,
on peut supposer $M$ de type fini sur $\Z$.
Alors on a des isomorphismes canoniques (cf [P1] 1.10.2)
$$Ind_{\Gamma'}^{\Gamma}M\otimes Ind_{\Gamma_1}^{\Gamma}\Z\simeq Ind_{\Gamma_1}^{\Gamma}Res_{\Gamma_1}^{\Gamma}Ind_{\Gamma'}^{\Gamma}M$$
$$\bigoplus_{i\in I}M\otimes Ind_{\Gamma_1}^{\Gamma'}\Z\simeq\bigoplus_{i\in I}Ind_{\Gamma_1}^{\Gamma'}Res_{\Gamma_1}^{\Gamma'}M$$
et le morphisme $v$ ci-dessus devient
$$\left(f:\Gamma\fl Res_{\Gamma_1}^{\Gamma}Ind_{\Gamma'}^{\Gamma}M\right)\fle\sum_{i\in I}(f_{\delta_i}:\gamma'\fle f(\gamma')(\delta_i)).$$

Comme dans la preuve de [P1] 1.11.5, on peut remplacer $M$ par un complexe $M^\bullet=Ind^{\Gamma'}_{\{1\}}A^\bullet$.
D'après [P1] 1.6.1, les composantes de $Ind_{\Gamma_1}^{\Gamma}Res^{\Gamma}_{\Gamma_1}Ind_{\Gamma'}^{\Gamma}M^\bullet$ sont $R\Gamma(\Gamma,\quad)$-acycliques,
et celles de $Ind_{\Gamma'}^{\Gamma_1}Res_{\Gamma_1}^{\Gamma'}M^\bullet$ sont $R\Gamma(\Gamma',\quad)$-acycliques,
donc le morphisme $u$ devient dans ce cas
$$\left(Ind_{\Gamma_1}^{\Gamma}Res_{\Gamma_1}^{\Gamma}Ind_{\{1\}}^{\Gamma}A^\bullet\right)^{\Gamma}\fl \left(Ind_{\Gamma_1}^{\Gamma}Res_{\Gamma_1}^{\Gamma}Ind_{\{1\}}^{\Gamma}A^\bullet\right)^{\Gamma'}\fl\left(\bigoplus_{i\in I}Ind_{\Gamma_1}^{\Gamma'}Res_{\Gamma_1}^{\Gamma'}Ind_{\{1\}}^{\Gamma'}A^\bullet\right)^{\Gamma'}.$$
Or on a
$$\left(Ind_{\Gamma_1}^{\Gamma}Res_{\Gamma_1}^{\Gamma}Ind_{\{1\}}^{\Gamma}A^\bullet\right)^{\Gamma}=\left(Res_{\Gamma_1}^{\Gamma}Ind_{\{1\}}^{\Gamma}A^\bullet\right)^{\Gamma_1}$$
$$\left(\bigoplus_{i\in I}Ind_{\Gamma_1}^{\Gamma'}Res_{\Gamma_1}^{\Gamma'}Ind_{\{1\}}^{\Gamma'}A^\bullet\right)^{\Gamma'}=\left(\bigoplus_{i\in I}Res_{\Gamma_1}^{\Gamma'}Ind_{\{1\}}^{\Gamma'}A^\bullet\right)^{\Gamma_1}$$
et $u$ est induit par le morphisme
$$\begin{array}{rcl}Res_{\Gamma_1}^{\Gamma}Ind_{\{1\}}^{\Gamma}A^\bullet & \fl & \bigoplus_{i\in I}Res_{\Gamma_1}^{\Gamma'}Ind_{\{1\}}^{\Gamma'}A^\bullet \\
(f:\Gamma\fl A^\bullet) & \fle & \sum_{i\in I}(f_{\delta_i}:\gamma'\fle f(\gamma'\delta_i))\end{array}$$
Ce dernier morphisme est clairement un isomorphisme.

\end{preuve}

\subsection{Systèmes locaux étales sur les modèles canoniques}

\subsubsection{Systèmes locaux provenant de représentations du groupe}

Revenons à la situation de 2.1.1. 
On fixe donc un nombre premier $\ell$, une donnée de Shimura $(\G,\X)$, un sous-groupe compact ouvert $\K\subset\G(\Af)$ et $V\in Rep_{\G(\Q_\ell)}$. 
Soit $m\in\Nat^*$.
On veut construire un faisceau étale sur $M^{\K}(\G,\X)$ qui sur $\C$ redonne le faisceau $\F^{\K}(V)(\Z/\ell^m\Z)$ de 2.1.2. 
Pour cela, on applique la construction de Pink au $\K$-recouvrement 
$$\varphi:M(\G,\X)\fl M^{\K}(\G,\X).$$
On fait agir $\K$ sur $V(\Q_\ell)$ via
$$\K\subset\G(\Af)\fl\G(\Q_\ell).$$
On choisit comme avant un $\Z_\ell$-réseau $\K$-invariant $V(\Z_\ell)$ de $V$. 
Alors $V(\Z/\ell^m\Z)=V(\Z_\ell)\otimes_{\Z_\ell}\Z/\ell^m\Z= V(\Z_\ell)/\ell^m V(\Z_\ell)$ est dans $Mod'_\K$ et on peut lui appliquer le foncteur $\mu_{\K,\varphi}$.

\begin{proposition} On a un isomorphisme canonique
$$\F^{\K}(V)(\Z/\ell^m\Z)\iso\mu_{\K,\varphi}V(\Z/\ell^m\Z)(\C),$$
où $\mu_{\K,\varphi}V(\Z/\ell^m\Z)(\C)$ est le faisceau en groupes abéliens sur $M^{\K}(\G,\X)(\C)$ provenant du faisceau étale $\mu_{\K,\varphi}V(\Z/\ell^m\Z)$ sur $M^{\K}(\G,\X)$.

\end{proposition}

\begin{preuve} ([Ln] p 38) 

Comme $V(\Z/\ell^m\Z)$ est fini, il existe un sous-groupe ouvert distingué $\K_0$ de $\K$ qui agit trivialement sur $V(\Z/\ell^m\Z)$. Si on note $\varphi_0$ la projection $T_1:M^{\K_0}(\G,\X)\fl M^{\K}(\G,\X)$, on a
$$\lambda_{\K,\varphi}V(\Z/\ell^m\Z)=\underline{V(\Z/\ell^m\Z)}\otimes\varphi_{0*}\underline{\Z}$$
et
$$\mu_{\K,\varphi}V(\Z/\ell^m\Z)=(\underline{V(\Z/\ell^m\Z)}\otimes\varphi_{0*}\underline{\Z})^{\K/\K_0},$$
donc 
$$\mu_{\K,\varphi}V(\Z/\ell^m\Z)(\C)=(\K/\K_0)\sous V(\Z/\ell^m\Z)\times M^{\K_0}(\G,\X)(\C),$$
où $\K/\K_0$ agit sur $V(\Z/\ell^m\Z)\times M^{\K_0}(\G,\X)(\C)$ par
$$k.(v,x,g\K_0)=(k^{-1}.v,x,gk\K_0).$$
On pose
$$\mathcal{E}=\coprod_{g\K\in\G(\Af)/\K} (g.V(\Z/\ell^m\Z)\times\X\times g\K/\K).$$
$\G(\Q)$ agit sur $\mathcal{E}$ comme dans 2.1.2, 
et on a $\F^{\K}(V)(\Z/\ell^m\Z)=\G(\Q)\sous\mathcal{E}$. 
On considère l'application de $\mathcal{E}$ dans $\mu_{\K,\varphi}V(\Z/\ell^m\Z)(\C)$, 
qui à $(g.v,x,g\K)$ associe $(v,x,g\K_0)$. Cette application est bien définie, car, si on remplace $g$ par $gk$ avec $k\in\K$, l'image de
$$(g.v,x,gk\K)=(gkk^{-1}.v,x,gk\K)$$
qu'on obtient est encore
$$(k^{-1}.v,x,gk\K_0)=(v,x,g\K_0).$$
De plus, si $\gamma\in\G(\Q)$, l'image de
$$\gamma.(g.v,x,g\K)=(\gamma g.v,\gamma.x,(\gamma g)\K)$$
est
$$(v,\gamma.x,\gamma g\K_0)=(v,x,g\K_0),$$
donc $\mathcal{E}\fl\mu_{\K,\varphi}V(\Z/\ell^m\Z)(\C)$ donne
$$\Phi:\G(\Q)\sous\mathcal{E}=\F^{\K}(V)(\Z/\ell^m\Z)\fl\mu_{\K,\varphi}V(\Z/\ell^m\Z)(\C),$$
qui est l'isomorphisme cherché :
\begin{itemize}
\item[$\bullet$] $\Phi$ est un homéomorphisme local;
\item[$\bullet$] $\Phi$ est surjective;
\item[$\bullet$] $\Phi$ est injective : Soient $(v,x,g),(v',x',g')\in V(\Z/\ell^m\Z)\times\X\times\G(\Af)$ 
tels que $\Phi(g.v,x,gK)=\Phi(g'.v',x',g'\K)$. 
Alors il existe $\gamma\in\G(\Q)$ et $k\in\K$ tels que $x'=\gamma.x$, $g'=\gamma g k$ et $v'=k^{-1}.v$, donc
$$(g'.v',x',g'\K)=(\gamma g.v,\gamma.x,\gamma g\K)=(g.v,x,g\K).$$
\end{itemize}

\end{preuve}

\begin{notation} Le faisceau $\ell$-adique 
$$\limpa {m\in\Nat^*} \mu_{\K,\varphi}(V(\Z/\ell^m\Z))\otimes\Q_\ell $$
sur $M^{\K}(\G,\X)$ sera noté $\F^{\K}V$ dans la suite. 

\end{notation}

\subsubsection{Généralisation}

Dans la suite, on utilisera des faisceaux $\ell$-adiques sur certains quotients de variétés de Shimura par des groupes finis, qui proviennent de représentations d'un groupe plus grand que celui de la donnée. 
On se place dans la situation suivante : $\B$ est un groupe linéaire connexe sur $\Q$, $\N_0$ est son radical unipotent, $\L=\B/\N_0$, $\G$ et $\G_\lin$ sont deux sous-groupes fermés réductifs connexes de $\L$ commutant entre eux, d'intersection finie et vérifiant $\L=\G\G_\lin$. 
On note $\QP_0$ l'image inverse $\G\N_0$ de $\G$ dans $\B$.
On considère une donnée de Shimura $(\G,\X)$. 
On suppose que l'action de $\G(\Q)$ sur $\X$ se prolonge en une action de $\L(\Q)$ 
telle que $\G_\lin(\Q)$ agisse trivialement
et on prolonge l'action à droite de $\G(\Af)$ sur $M(\G,\X)$ à $\L(\Q)\QP_0(\Af)$
en faisant agir trivialement $\G_\lin(\Q)\N_0(\Af)$.

Enfin, soit $\K$ un sous-groupe ouvert compact net de $\B(\Af)$; on note $\Hr=\K\cap\B(\Q)\QP_0(\Af)$, $\Hr_\lin=\K\cap Cent_{\B(\Q)}(\X)\N_0(\Af)$, $\K_Q=\K\cap\QP_0(\Af)\subset\Hr$, $\K_N=\K\cap\N_0(\Af)\subset\Hr_\lin$ et $\K_G=\K_Q/\K_N$. Le morphisme quotient
$$\varphi:M(\G,\X)\fl M(\G,\X)/\Hr=M$$
est un revêtement étale profini de groupe $\Hr/\Hr_\lin$, et $\Hr/\Hr_\lin\K_Q$ est fini, donc $M$ est le quotient de la variété de Shimura $M^{\K_G}(\G,\X)$ par un groupe fini. On sait qu'à tout $\Hr$-module $A\in Mod'_{\Hr}$ on peut associer un complexe de faisceaux étales sur $M$,
$$R\Gamma(\Hr,\lambda_{\Hr,\varphi}(A))\simeq\mu_{\Hr/\Hr_\lin,\varphi}(R\Gamma(\Hr_\lin,A)).$$ 
On voudrait, à partir d'une représentation algébrique de $\B(\Q_\ell)$ dans un $\Q_\ell$-espace vectoriel de dimension finie (ou d'un complexe borné de telles représentations), sur laquelle $\B(\Af)$ et ses sous-groupes agissent comme d'habitude via le morphisme $\B(\Af)\fl\B(\Q_\ell)$,
construire un complexe borné de faisceaux $\ell$-adiques sur $M$.

\begin{lemme}\label{lemme:cohomologie_N}([P2] 5.2.2) Soient $\U$ un groupe algébrique unipotent connexe sur $\Q$
et $\K_U$ un sous-groupe ouvert compact de $\U$.
Alors, pour tout $V\in D^b(Rep_{\U(\Q_\ell)})$, le morphisme de restriction
$$R\Gamma_{cont}(\K_U,V)\fl R\Gamma(\K_U\cap\U(\Q),V)$$
est un isomorphisme.

\end{lemme}

Le lemme ci-dessous, qui est une version algébrique du théorème de van Est, est prouvé dans [GHM] 24.

\begin{lemme}\label{lemme:van_Est}Soient $\U$ un groupe algébrique unipotent connexe sur $\Q$ et $\Gamma_U$ un sous-groupe arithmétique de $\U(\Q)$. 
On note $Mod_\U$ la catégorie des limites inductives de représentations algébriques rationnelles de dimension finie de $\U$ (ou, ce qui revient au même, de $Lie(\U)$),
et $A_U$ la $\Q$-algèbre des polynômes à coefficients rationnels sur $\U(\Q)$, qui contient la sous-algèbre $\Q$ des constantes.
Pour tout $M$, on définit une suite exacte
$$0\fl M\stackrel{u_0}{\fl}I^0(M)\stackrel{u_1}{\fl} I^1(M)\stackrel{u_2}{\fl}I^2(M)\stackrel{u_3}{\fl}\dots$$
par :
\begin{itemize}
\item[$\bullet$] $u_0:M\fl I^0(M)$ est l'injection évidente $M\simeq M\otimes\Q\fl M\otimes A_U=I^0(M)$;
\item[$\bullet$] pour tout $i\in\Nat$, $u_{i+1}:I^i(M)\fl I^{i+1}(M)$ est le morphisme évident $I^i(M)\fl Coker(u_i)\fl Coker(u_i)\otimes A_U=I^{i+1}(M)$.

\end{itemize}

Alors :
\begin{itemize}
\item[(i)] Pour tout $M\in Mod_\U$, $M\otimes A_U$ est un objet injectif de $Mod_\U$ et un objet acyclique pour le foncteur $(\quad)^{\Gamma_U}$ (invariants par $\Gamma_U$).
\item[(ii)] Pour tout $M\in Mod_\U$, le morphisme suivant est un isomorphisme :
$$R\Gamma(Lie(\U),M)=R\Gamma(\U(\Q),M)\simeq I^\bullet(M)^{\U(\Q)}\fl I^{\bullet}(M)^{\Gamma_U}\simeq R\Gamma(\Gamma_U,M).$$

\end{itemize}

\end{lemme}

\begin{corollaire} Soit $V\in D^b(Rep_{\B(\Q_\ell)})$. Alors on a un isomorphisme canonique
$$R\Gamma(\Hr_\lin,V)\simeq R\Gamma(\Hr_\lin/\K_N,R\Gamma(Lie(\N_0),V)).$$

\end{corollaire}

Si $V\in D^b(Rep_{\B(\Q_\ell)})$, $R\Gamma(Lie(\N_0),V)\in D^b(Rep_{\L(\Q_\ell)})$, donc le corollaire permet de se ramener au cas $\B=\L$. Traitons ce cas.

Comme $\K$ est net, $\Gamma_\lin=\Hr_\lin$ est net, et $\Gamma_\lin\cap\G(\Q)=\{1\}$ (car le stabilisateur dans $\G(\R)$ d'un élément de $\X$ est compact modulo le centre de $\G(\R)$). $\Gamma_\lin$ s'injecte donc dans $(\L/\G)(\Q)$, et son image est un sous-groupe arithmétique. 
On a un autre sous-groupe arithmétique de $(\L/\G)(\Q)$, $\Gamma'_\lin=\Hr/\K_G$, et $\Gamma_\lin$ s'identifie à un sous-groupe distingué de $\Gamma'_\lin$ (forcément d'indice fini).
D'après [BS] 11.1, la compactification de Borel-Serre partielle d'un espace symétrique de $\L/\G$ permet de construire 
une résolution bornée $L_\bullet\fl\Z$ du $\Z[\Gamma'_\lin]$-module trivial $\Z$ par des $\Z[\Gamma'_\lin]$-modules libres de type fini 
(les $L_i$ sont donc aussi des $\Z[\Gamma_\lin]$-modules libres de type fini).
On fait agir $\Hr$ sur les $L_i$ via la projection $\Hr\fl\Gamma'_\lin$. Si $A$ est un $\Hr$-module fini, alors $Hom_{\Z[\Gamma_\lin]}(L_\bullet,A)$, où $\Hr/\Gamma_\lin$ agit par $(h\Gamma_\lin.f)(x)=hf(h^{-1}x)$, est un complexe borné de $\Hr/\Gamma_\lin$-modules finis qui représente l'objet $R\Gamma(\Gamma_\lin,A)$ de $D^b(Mod_{\Hr/\Gamma_\lin})$. 

\begin{definition}\label{def_faisceaux} Soit $V\in Rep_{\L(\Q_\ell)}$.  
On choisit un réseau $\K$-invariant $\Lambda\subset V$. Alors le système projectif des $\mu_{\Hr/\Gamma_\lin,\varphi}(Hom_{\Z[\Gamma_\lin]}(L_\bullet,\Lambda/\ell^m\Lambda)),m\in\Nat^*$, avec les morphismes de transition évidents, est un complexe borné de faisceaux $\ell$-adiques constructibles sur $M$. On note $\F^{\Hr/\Gamma_\lin}R\Gamma(\Gamma_\lin,V)$ l'objet correspondant de $D^b_c(M,\Q_\ell)$, dont la classe d'isomorphisme ne dépend pas du choix de $\Lambda$. Cette construction s'étend trivialement au cas où $V$ est un complexe borné de représentations de $\L(\Q_\ell)$, et donne un foncteur triangulé
$$\F^{\Hr/\Gamma_\lin}R\Gamma(\Gamma_\lin,\quad):D^b(Rep_{\L(\Q_\ell)})\fl D^b_c(M,\Q_\ell).$$

Si on revient au cas général du début de ce paragraphe ($\L$ est le quotient de Levi d'un groupe connexe $\B$), on obtient un foncteur triangulé
$$\begin{array}{rcl}\F^{\Hr/\Hr_\lin}R\Gamma(\Hr_\lin,\quad):D^b(Rep_{\B(\Q_\ell)}) & \fl & D^b_c(M,\Q_\ell) \\
V & \fle & \F^{\Hr/\Hr_\lin}R\Gamma(\Hr_\lin/\K_N,R\Gamma(Lie(\N_0),V))\end{array}.$$

\end{definition}

\begin{remarque}\label{rq:def_faisceaux} De même, si $A$ est un $\Q_\ell$-espace vectoriel avec une action de $\Hr/\Hr_\lin$ qui admet un $\Z_\ell$-réseau $\Lambda$ invariant par $\Hr/\Hr_\lin$, alors on lui associe le faisceau $\ell$-adique lisse sur $M$
$$\F^{\Hr/\Hr_\lin}A=(\mu_{\Hr/\Hr_\lin,\varphi}(\Lambda/\ell^m\Lambda))_{m\in\Z},$$
dont la classe d'isomorphisme ne dépend pas du choix de $\Lambda$. 

Par exemple, dans le cas $\B=\L$, si $V\in Rep_{\L(\Q_\ell)}$, on a pour tout $i\in\Z$ un faisceau $\F^{\Hr/\Gamma_\lin}H^i(\Gamma_\lin,V)$, qui s'identifie à $H^i(\F^{\Hr/\Gamma_\lin}R\Gamma(\Gamma_\lin,V))$ (grâce à l'exactitude de $\mu_{\Hr/\Gamma_\lin,\varphi}$).

La notation $\F^{\Hr/\Hr_\lin}$ qu'on vient d'introduire est une généralisation de la notation $\F^{\K}$ du paragraphe précédent.

\end{remarque}

\section{Le théorème de Pink (prolongement des faisceaux sur la compactification de Baily-Borel)}

Soient $\G=\GU(p,q)$, $(\G,\X)$ la donnée de Shimura pure définie dans la section 1.2. 
On va énoncer un résultat très légèrement plus général que celui de Pink (mais qui se prouve exactement de la même façon), en se plaçant dans la situation de la définition \ref{def_faisceaux}.
 
Soient donc $\B$ un groupe algébrique connexe sur $\Q$, $\N_0$ son radical unipotent, $\L=\B/\N_0$ son quotient de Levi. On suppose que $\G$ s'identifie à un sous-groupe fermé de $\L$, et qu'on a un sous-groupe réductif connexe $\G_\lin\subset\L$ commutant avec $\G$, d'intersection finie avec $\G$, et tel que $\L=\G\G_\lin$. Enfin, on note $\QP_0$ l'image réciproque de $\G$ dans $\B$ (donc $\G=\QP_0/\N_0$). Si $\G'$ est un sous-groupe de $\L$, on notera souvent $\G'\N_0$ son image réciproque dans $\B$. 
On suppose qu'on peut prolonger l'action à droite de $\G(\Af)$ sur $M(\G,\X)$ à 
$\B(\Q)\QP_0(\Af)$ en faisant agir trivialement $\G_\lin(\Q)\N_0(\Af)$. Soit $\K$ un sous-groupe compact ouvert net de $\B(\Af)$. On note $\Hr=\K\cap\B(\Q)\QP_0(\Af)$, $\Hr_\lin=\K\cap Cent_{\B(\Q)}(\X)\N_0(\Af)$, et $\varphi$ le recouvrement étale $\til{M}=M(\G,\X)\fl M=M(\G,\X)/\Hr$, qui est galoisien profini de groupe $\Hr/\Hr_\lin$. L'action de $\B(\Q)\QP_0(\Af)$ s'étend évidemment à $M(\G,\X)^*$, et on note $j$ l'immersion ouverte $M\fl M^*=M(\G,\X)^*/\Hr$. Dans [P2], Pink a calculé, pour un $\Hr$-module de torsion $M$ et dans le cas $\B=\G$, la restriction de $Rj_*R\Gamma(\Hr,\lambda_{\Hr,\varphi}M)$ aux strates de $M^*$. 
Nous allons énoncer ce théorème dans le cas où $\B$ n'est plus forcément égal à $\G$. 

Soient $\P$ un sous-groupe parabolique maximal de $\G$ et $(\QP,\Y)$ la composante de bord associée. On note comme avant $\N$ le radical unipotent de $\P$, $\L_P=\P/\N$ le quotient de Levi, et $(\G_1,\X_1)=(\QP,\Y)/\N$. Soit $g\in\B(\Q)\QP_0(\Af)$. On pose  
$$\K_{Q}=(\QP\N_0)(\Af)\cap g\K g^{-1}$$
$$\K_{N}=(\N\N_0)(\Af)\cap g\K g^{-1}$$
$$\K_{G}=\K_{Q}/\K_{N}\subset\G_1(\Af)$$
$$\Hr_{P}=g\K g^{-1}\cap Stab_{\B(\Q)}(\X_1)(\QP\N_0)(\Af)=g\K g^{-1}\cap 
(\G_\lin\P)(\Q)(\QP\N_0)(\Af)$$
$$\Hr_{P,\lin}=g\K g^{-1}\cap Cent_{\B(\Q)}(\X_1)(\N\N_0)(\Af)=g\K g^{-1}\cap 
(\G_\lin\L_{P,\lin})(\Q)(\N\N_0)(\Af),$$
où $\L_{P,\lin}$ est la partie linéaire de $\L_P$.

On a un diagramme commutatif :

$$\xymatrix{\til{M}= M(\G,\X)\ar@{^{(}->}[r]^{\til{j}}\ar@<3.5ex>[d]_\varphi & \til{M}^*= M(\G,\X)^*\ar@<3.5ex>[d]_{\overline{\varphi}} & M(\G_1,\X_1)=\til{M}_1 \ar@{_{(}->}[l]_{\til{i}}\ar@<-3ex>[d]^{\varphi_1}  \\
            M= M(\G,\X)/\Hr\ar@{^{(}->}[r]^j & M^*= M(\G,\X)^*/\Hr & M(\G_1,\X_1)/\Hr_P=M_1 \ar@{_{(}->}[l]_{i}}$$
où $\varphi_1$ est un recouvrement étale galoisien profini de groupe $\Hr_P/\Hr_{P,\lin}$.

Soit $TorMod_{\Hr}$ la sous-catégorie pleine de $Mod'_{\Hr}$ dont les objets sont les $\Hr$-modules de torsion. 
On note encore $\lambda_{\Hr,\varphi}$ la restriction de $\lambda_{\Hr,\varphi}:Mod'_{\Hr}\fl \acute{E}t_{M,\Hr}$ à $TorMod_{\Hr}$.

\begin{theoreme}\label{th_Pink}(cf [P2] 4.2.1) On note 
$\theta:\Hr_{P}\fl\Hr$, $h\fle g^{-1}hg$, qui donne $\theta^*:TorMod_{\Hr}\fl TorMod_{\Hr_P}$. 

On a un isomorphisme canonique de foncteurs $D^+(TorMod_{\Hr})\fl D^+(\acute{E}t_{M_1})$
$$i^* Rj_*R\Gamma(\Hr,\quad)\lambda_{\Hr,\varphi}\simeq R\Gamma(\Hr_P,\quad)\lambda_{\Hr_P,\varphi_1}\theta^*,$$
ou, autrement dit,
$$i^*Rj_*\mu_{\Hr/\Hr_\lin,\varphi}R\Gamma(\Hr_\lin,\quad)\simeq\mu_{\Hr_P/\Hr_{P,\lin},\varphi_1}R\Gamma(\Hr_{P,\lin},\quad)\theta^*.$$

$$\xymatrix@C+30pt{D^+(TorMod_{\Hr})\ar[d]^{\mu_{\Hr/\Hr_\lin,\varphi}R\Gamma(\Hr_\lin,\quad)}\ar[r]^{\theta^*} & D^+(TorMod_{\Hr_{P}})\ar[r]^{R\Gamma(\Hr_{P,\lin},\quad)} & D^+(TorMod_{\Hr_P/\Hr_{P,\lin}})\ar[d]^{\mu_{\Hr_P/\Hr_{P,\lin},\varphi_1}} \\
            D^+(\acute{E}t_{M})\ar[r]_{Rj_*} & D^+(\acute{E}t_{M^*})\ar[r]_{i^*} & D^+(\acute{E}t_{M_1})}$$

\end{theoreme}

Voyons ce que ce théorème donne pour les systèmes de coefficients définis plus haut.

\begin{corollaire}\label{restriction_bord}(cf [P2] 5.3.1) Soit $V\in D^b(Rep_{\B(\Q_\ell)})$, sur lequel on fait comme d'habitude agir $\Hr$ via le morphisme $\Hr\subset\B(\Af)\fl\B(\Q_\ell)$. 
On a construit plus haut (voir la définition \ref{def_faisceaux}) un complexe de faisceaux 
$\ell$-adiques $\F^{\Hr/\Hr_\lin}R\Gamma(\Hr_\lin,V)\in D^b_c(M,\Q_\ell)$. De même, 
en utilisant cette fois l'inclusion $\G_1\subset\L_P$, on peut associer à tout complexe $V'\in D^b(Rep_{(\P\N_0)(\Q_\ell)})$ un complexe de faisceaux $\ell$-adiques $\F^{\Hr_P/\Hr_{P,\lin}}R\Gamma(\Hr_{P,\lin},V')\in D^b_c(M_1,\Q_\ell)$. 

On a des isomorphismes canoniques
$$\begin{array}{rcl}i^*Rj_*\F^{\Hr/\Hr_\lin}R\Gamma(\Hr_\lin,V) & \simeq & \F^{\Hr_P/\Hr_{P,\lin}}R\Gamma(\Hr_{P,\lin},V) \\
& \simeq & \F^{\Hr_P/\Hr_{P,\lin}}R\Gamma(\Hr_{P,\lin}/\K_N,R\Gamma(\K_N,V)) \\
& \simeq & \F^{\Hr_P/\Hr_{P,\lin}}R\Gamma(\Hr_{P,\lin}/\K_N,R\Gamma(Lie(\N\N_0),V)).\end{array}$$

\end{corollaire}

Le corollaire résulte du théorème ci-dessus, et
des lemmes \ref{lemme:cohomologie_N} et \ref{lemme:van_Est}.

\newpage

\chapter{Modèles entiers}

Le but de cette partie est de trouver un ensemble fini de nombres premiers $\Sigma$
tel que les variétés et les systèmes de coefficients sur le corps reflex $F$ des parties 1 et 2 s'étendent
(en gardant les mêmes propriétés) sur l'anneau $\O_F[1/\Sigma]$.
Pour la variété de Shimura elle-même, $\Sigma$ est explicite et ne dépend que du discriminant de $E$ et du niveau (voir la section 3.1).
Malheureusement, comme on ne dispose pas de compactifications toroïdales des modèles entiers,
on doit dès que l'on veut travailler avec la compactification de Baily-Borel inverser 
un certain ensemble d'autres nombres premiers sur lequel on n'a aucune information
(à part le fait qu'il est fini).

\section{Variété de Shimura}

$(\G,\X)$ est toujours la donnée de Shimura définie dans la section 1.2, 
et on utilise les notations du chapitre 1.

Les variétés de Shimura qu'on considère ici sont des variétés PEL,
c'est-à-dire 
qu'elles s'interprètent comme des espaces de modules de variétés abéliennes avec des polarisations, des endomorphismes imposés et des structures de niveau. Cette description modulaire permet de construire des modèles des variétés de Shimura sur l'anneau des entiers de son corps reflex, où on a inversé certains nombres premiers. 
La référence qu'on utilisera ici est [K2].

\subsection{Schémas abéliens et structures de niveau}

Commençons par mettre une structure entière sur le groupe $\G$. Il en existe une évidente : si $A$ est une $\Z$-algèbre, on pose
$$\G(A)=\{g\in\GL_n(A\otimes_\Z\O_E),g^*J_{p,q}g=\lambda J_{p,q},\lambda\in A^\times\}.$$

On note $D$ le discriminant de $E=\Q[\sqrt{-d}]$ (donc $D=-4d$ si $d\not=3\mbox{ mod }4$ et $D=-d$ si $d=3\mbox{ mod }4$). On rappelle qu'on a noté $F$ le corps reflex de $(\G,\X)$ ($F=E$ si $p>q$ et $F=\Q$ si $p=q$). On pose $\alpha=\sqrt{-d}$ si $d\not =3\mbox{ mod }4$ et $\alpha=(1+\sqrt{-d})/2$ si $d=3\mbox{ mod }4$; on a donc
$$\O_E=\Z\oplus\Z\alpha=\Z[\alpha].$$

Soit $V$ le $E$-espace vectoriel $E^{p+q}$, muni de la forme hermitienne $J$ de matrice $J_{p,q}$ dans la base canonique. On a une forme alternée $\psi:V\times V\fl \Q$ associée à $J$ par la formule
$$J(v,w)=\psi(v,\alpha w)-\overline{\alpha}\psi(v,w).$$
Un calcul simple donne
$$\psi(v,w)=\left\{\begin{array}{ll}\frac{1}{\sqrt{d}}Im J(v,w) & \mbox{ si }d\not=3\mbox{ mod }4 \\
\frac{2}{\sqrt{d}}Im J(v,w) & \mbox{ si }d=3\mbox{ mod }4\end{array}\right.$$
En particulier, le $\O_E$-réseau $\Lambda=\O_E^{p+q}$ de $V$ est autodual pour $\psi$.

\begin{notation} 
\begin{itemize}
\item[(1)] Si $A$ et $B$ sont deux schémas abéliens sur un schéma $S$, 
on note $Hom_S(A,B)$ (ou $Hom(A,B)$) le groupe abélien des morphismes de $A$ dans $B$. 
\item[(2)] Si $A\fl S$ est un schéma abélien, on note $\widehat{A}\fl S$ le schéma abélien dual. 
Si $f:A\fl B$ est un morphisme de schémas abéliens sur $S$, on note $\widehat{f}:\widehat{B}\fl\widehat{A}$ le morphisme dual.
\item[(3)] Soient $A\fl S$ un schéma abélien et $\lambda$ une polarisation de $A$.
On définit une involution $\iota_\lambda$ de $End(A)$, dite involution de Rosati associée à $\lambda$ : si $f\in End(A)$, $\iota_\lambda(f)=\widehat{\lambda}\widehat{f}\lambda$.

$$\xymatrix{A & A\ar[d]^{\lambda} \\
\widehat{A}\ar[u]^{\widehat{\lambda}} & \widehat{A}\ar[l]^{\widehat{f}}}$$

\item[(4)] Soit $R$ un anneau commutatif. La catégorie des schémas abéliens sur $S$ à $R$-polarisation près est la catégorie additive (et même $R$-linéaire) dont les objets sont les schémas abéliens sur $S$ et où le groupe des morphismes de $A$ dans $B$ est $Hom_S(A,B)\otimes_\Z R$. Une $R$-isogénie de $A$ dans $B$ est un isomorphisme dans cette catégorie. Comme plus haut, une $R$-isogénie de $A$ dans $\widehat{A}$ donne une involution de  Rosati sur $End(A)\otimes_\Z R$.

\item[(5)] Soit $\pi:A\fl S$ un schéma abélien, $\pi':\widehat{A}\fl S$ le schéma abélien dual 
 et $\ell$ un nombre premier inversible sur $S$. 
On a une forme bilinéaire non dégénérée canonique
$$R^1\pi_*\Q_\ell\times R^1\pi'_*\Q_\ell\fl\Q_\ell(1)_S.$$
Si $\lambda:A\fl\widehat{A}$ est une $\Q$-isogénie,
elle induit un isomorphisme $R^1\pi_*\Q_\ell\iso R^1\pi'_*\Q_\ell$,
d'où une forme alternée non dégénérée
$$R^1\pi_*\Q_\ell\times R^1\pi_*\Q_\ell\fl\Q_\ell(1)_S,$$
qu'on appellera accouplement de Weil associé à $\lambda$.

\end{itemize}
\end{notation}

\begin{notation} Pour tout nombre premier $\ell$, on note $\Lambda_\ell=\Lambda\otimes_\Z\Z_\ell$  et
$$\K_{0,\ell}=\{g\in\G(\Q_\ell)\mbox{ tq }g(\Lambda_\ell)=\Lambda_\ell\}=\G(\Z_\ell).$$
$\K_{0,\ell}$ est un sous-groupe ouvert compact de $\G(\Q_\ell)$.

Soit $\Sigma$ un ensemble de nombres premiers. On note
$$\Z[1/\Sigma]=\Z[1/\ell,\ell\in \Sigma]$$
$$\Z_\Sigma=\prod_{\ell\in \Sigma}\Z_\ell$$
$$\Z^\Sigma=\prod_{\ell\not\in \Sigma}\Z_\ell$$
$${\Af}_\Sigma=\Z_\Sigma\otimes_\Z\Q$$
$$\Af^\Sigma=\Z^\Sigma\otimes_\Z\Q$$
$$\K^\Sigma_0=\prod_{\ell\not\in \Sigma}\K_{0,\ell}.$$
Si $R$ est un anneau commutatif, on note
$$R[1/\Sigma]=R\otimes_\Z\Z[1/\Sigma].$$

\end{notation}

La définition des structures de niveau que nous donnons ici est directement inspirée de celle de [K2] 5.

\begin{definition}\label{def:struct_niveau} Soit $A\fl S$ un schéma abélien, avec $S$ connexe, et $\Sigma$ un ensemble de nombres premiers inversibles sur $S$. On se donne une $\Z[1/\Sigma]$-isogénie $\lambda:A\fl\widehat{A}$ et un morphisme d'anneaux $i:\O_E[1/\Sigma]\fl End(A)\otimes_\Z\Z[1/\Sigma]$ qui envoie la conjugaison complexe sur l'involution de Rosati associée à $\lambda$.
Enfin, soit $\K$ un sous-groupe ouvert compact de $\G(\Af)$ de la forme $\K=\K_\Sigma\K_0^\Sigma$, avec $\K_\Sigma\subset\G({\Af}_\Sigma)$.

Soit $s$ un point géométrique de $S$. Une structure de niveau $(\K,\Sigma)$ sur $(A,\lambda,i)$ 
est une $\K$-orbite $\overline{\eta}$ stable par $\pi_1(S,s)$ d'isomorphismes $E$-linéaires 
$\eta:H^1(A_s,{\Af}_\Sigma)\fl V\otimes_\Q{\Af}_\Sigma$ 
(avec l'action de $E$ sur $H^1(A_s,{\Af}_\Sigma)$ provenant de $i$) 
qui envoient l'accouplement de Weil associé à $\lambda$ sur un multiple 
(par un scalaire de ${\Af}_\Sigma^\times$) de la forme alternée $\psi$. 
Cette notion est essentiellement indépendante du choix de $s$.

\end{definition}

\begin{remarque} \begin{itemize}
\item[(a)] On parlera de structure de niveau $\K$ si la définition de $\Sigma$ est claire.
\item[(b)] Si $\K\subset\K_0=\prod_\ell\K_\ell$ est un sous-groupe compact ouvert, il existe un $\Sigma$ fini qui vérifie la condition de la définition. 

\end{itemize}
\end{remarque}

\subsection{Le problème de modules}

\begin{definition}\label{def:prob_mod} Soient $\Sigma$ un ensemble de nombres premiers qui contient tous les facteurs premiers de $D$, et $\K\subset\G(\Af)$ un sous-groupe ouvert compact vérifiant la condition de la définition \ref{def:struct_niveau}. 
On note $\Mod^{\K}_\Sigma(\G,\X)$ la catégorie fibrée en groupoïdes sur la catégorie $Sch/\O_F[1/\Sigma]$ des $\O_F[1/\Sigma]$-schémas localement noethériens suivante : pour tout $\O_F[1/\Sigma]$-schéma localement noethérien $S$ 
\begin{itemize}
\item[(i)] les objets de $\Mod^{\K}_\Sigma(\G,\X)(S)$ sont les quadruplets $(A,\lambda,i,\overline{\eta})$, où $A\fl S$ est un schéma abélien, $\lambda$ est une polarisation de $A$ qui est une $\Z[1/\Sigma]$-isogénie, $i:\O_E[1/\Sigma]\fl End(A)[1/\Sigma]$ est un morphisme d'anneaux qui envoie la conjugaison complexe sur l'involution de Rosati associée à $\lambda$, et $\overline{\eta}$ est la donnée d'une structure de niveau $(\K,\Sigma)$ au-dessus de chaque composante connexe de $S$, vérifiant la condition suivante (dite condition du déterminant) :
$$det(X_1+\alpha X_2,Lie(A)^*)=(X_1+\alpha X_2)^p(X_1+\overline{\alpha}X_2)^q\in\Gamma(S,\O_S)[X_1,X_2]$$
\item[(ii)] les morphismes de $(A,\lambda,i,\overline{\eta})$ dans $(A',\lambda',i',\overline{\eta}')$ sont les $\Z[1/\Sigma]$-isogénies $A\fl A'$ qui commutent à l'action de $\O_E[1/\Sigma]$, envoient $\lambda$ sur un multiple (par un scalaire de $\Z[1/\Sigma]^\times$) de $\lambda'$ et $\overline{\eta}$ sur $\overline{\eta}'$.

\end{itemize}

\end{definition}

\begin{remarque}\begin{itemize}
\item[(a)] Les schémas abéliens qui interviennent dans le problème de modules précédent sont forcément de dimension $p+q$ (c'est une conséquence de la condition du déterminant ou de l'existence de la structure de niveau).
 
\item[(b)] $\Mod^{\K}_\Sigma(\G,\X)$ est un champ fppf.

\item[(c)] Si $S$ est un schéma sur $\O_F[1/\Sigma]$, le polynôme $(X_1+\alpha X_2)^p(X_1+\overline{\alpha}X_2)^q$
peut bien être vu comme un polynôme à coefficients dans $\Gamma(S,\O_S)$,
car, si $p=q$ (le seul cas où $F=\Q$), ce polynôme est à coefficients dans $\Z$.\newline
Demander que $(A,i)$ vérifie la condition du déterminant revient à demander que 
le $\O_S$-module localement libre $Lie(A)^*\otimes\Z[1/\Sigma]$, muni de l'action 
de $\O_E$ provenant de $i$, soit de type $(p,q)$ au sens de [Be] I.2.

\item[(d)] Soit $p$ un nombre premier qui ne divise pas $D$, $\K\subset\G(\Af)$ un 
sous-groupe compact ouvert de la forme $\K=\K_p\K^p$, avec $\K^p\subset\G(\Af^p)$, 
et $\Sigma$ l'ensemble des nombres premiers différents de $p$. 
Alors $\Mod^{\K}_\Sigma(\G,\X)$ est le problème de modules de [K2] 5, excepté que
Kottwitz impose aux schémas abéliens d'être projectifs (mais on peut tout de même utiliser les résultats de [K2] sur les points du problème de modules à valeurs dans un corps)
et que nous avons remplacé $H_1$ par $H^1$ et $Lie(A)$ par $Lie(A)^*$.

Indiquons la correspondance entre les notations utilisées ici et celles de [K2] 5 : 
l'algèbre simple $B$ de [K2] est le corps $E$, le $\Z_{(p)}$-ordre $\O_B$ est 
$\O_E\otimes_\Z\Z_{(p)}$, l'involution * sur $B$ est la conjugaison complexe, 
le $B$-module $V$ est aussi noté par $V$ ici, la forme alternée $(.,.)$ sur $V$ est $\psi$, 
le réseau autodual $\Lambda_0$ de $V_{\Q_p}$ est $\Lambda_p$. 
Les symboles $\G$ et $\K^p$ désignent les mêmes objets que dans [K2], et le $h:\C\fl\G(\R)$ de [K2] est le prolongement évident du $h_0:\C^\times\fl\G(\R)$ de la partie I.
Le corps reflex (noté $E$ dans [K2]) est noté $F$ ici. Enfin, la base $(\alpha_1,\dots,\alpha_t)$ du $\Z_{(p)}$-module $\O_B$ est ici $(1,\alpha)$.

\end{itemize}

\end{remarque}

\begin{lemme} Pour tout $N\in\Nat^*$, notons
$$\K(N)=\{g\in\G(\Af),(g-1)(\Lambda\otimes_\Z\widehat{\Z})\subset N(\Lambda\otimes_\Z\widehat{\Z})\}.$$
Soit $(\K,\Sigma)$ comme plus haut. S'il existe $N\geq 3$ tel que $\K\subset\K(N)$, alors les objets de $\Mod^{\K}_\Sigma(\G,\X)$ n'ont pas d'automorphismes non triviaux.

\end{lemme}

\begin{proposition}\label{independance_de_T} Soient $\Sigma_1\subset \Sigma_2$ des ensembles de nombres premiers (contenant les facteurs premiers de $D$), 
et $\K\subset\G(\Af)$ un sous-groupe ouvert compact qui s'écrit $\K=\K_{\Sigma_1}\K_0^{\Sigma_1}$, avec $\K_{\Sigma_1}\subset\G({\Af}_{\Sigma_1})$. 
On a un $1$-morphisme évident de $\Mod^{\K}_{\Sigma_2}(\G,\X)$ dans $\Mod^{\K}_{\Sigma_1}(\G,\X)_{|Sch/\O_F[1/\Sigma_2]}$ : 
soit $S$ un schéma localement noethérien sur $\O_F[1/\Sigma_2]$; 
si $(A,\lambda,i,\overline{\eta})$ est un objet de $\Mod^{\K}_{\Sigma_2}(\G,\X)(S)$, 
on lui associe l'objet $(A,\lambda,i,\overline{\eta}')$ de $\Mod^{\K}_{\Sigma_1}(\G,\X)(S)$, 
où $\overline{\eta}'$ est obtenue en restreignant les $\eta\in\overline{\eta}$ à 
$H^1(A_s,{\Af}_{\Sigma_1})$.

Ce morphisme est un isomorphisme.

\end{proposition}

\begin{corollaire} Soit $(\K,\Sigma)$ comme dans la définition \ref{def:prob_mod} et $p\not\in \Sigma$.
($\K$ est bien sûr de la forme $\K_p\K^p$, avec $\K^p\subset\G(\Af^p)$.)
Alors $\Mod^{\K}_\Sigma(\G,\X)_{|Sch/\O_F\otimes\Z_{(p)}}$ est canoniquement isomorphe au problème de modules de [K2] 5.

\end{corollaire}

\begin{theoreme} Avec les notations de la définition \ref{def:prob_mod}, $\Mod^{\K}_\Sigma(\G,\X)$ est un champ de Deligne-Mumford connexe, lisse et de dimension relative $pq$ sur $\O_F[1/\Sigma]$. 
S'il existe $N\geq 3$ tel que $\K\subset\K(N)$, $\Mod^{\K}_\Sigma(\G,\X)$ est représentable par un espace algébrique.

\end{theoreme}

\begin{theoreme} Soient $(\K,\Sigma)$ comme dans la définition \ref{def:prob_mod}. Alors $\Mod^{\K}_\Sigma(\G,\X)_F$ est le modèle canonique de $M^{\K}(\G,\X)(\C)$, autrement dit, $\Mod^{\K}_\Sigma(\G,\X)_F$ et $M^{\K}(\G,\X)$ sont canoniquement isomorphes.

\end{theoreme}

\begin{preuve} On applique [K2] 8, en remarquant que $ker^1(\Q,\G)=\{1\}$ ([Sh] 5.8).

\end{preuve}

\subsection{Changement de niveau}

Soient $\Sigma$ un ensemble de nombres premiers contenant les facteurs premiers de $D$, $g\in\G({\Af}_\Sigma)$, $\K_\Sigma,\K'_\Sigma$ des sous-groupes compacts ouverts de $\G({\Af}_\Sigma)$ tels que $\K'_\Sigma\subset g\K_\Sigma g^{-1}$. Notons $\K=\K_\Sigma\K_0^\Sigma$ et $\K'=\K'_\Sigma\K_0^\Sigma$; on a $\K'\subset g\K g^{-1}$. 
On définit un $1$-morphisme
$$T_g:\Mod^{\K'}_\Sigma(\G,\X)\fl\Mod^{\K}_\Sigma(\G,\X)$$
de la manière suivante : Soit $S$ un $\O_F[1/\Sigma]$-schéma localement noethérien, qu'on peut supposer connexe. Si $(A,\lambda,i,\overline{\eta})$ est un objet de $\Mod^{\K'}_\Sigma(\G,\X)(S)$, la $\K$-orbite de $\eta g$ ne dépend pas du choix de $\eta\in\overline{\eta}$; on la note $\overline{\eta}'$, et on pose $T_g(A,\lambda,i,\overline{\eta})=(A,\lambda,i,\overline{\eta}')$.

\begin{proposition}\label{extension_operateurs_Hecke} On garde les notations ci-dessus. $T_g$ est un morphisme fini étale, et sa restriction à la fibre générique est le morphisme $T_g$ entre les modèles canoniques de la section 1.6.

Si $\K'_\Sigma\subset\K_\Sigma$ est un sous-groupe ouvert distingué, on obtient une action de $\K_\Sigma/\K'_\Sigma=\K/\K'$ sur $\Mod^{\K'}_\Sigma(\G,\X)$, et le quotient de $\Mod^{\K'}_\Sigma(\G,\X)$ par cette action est $\Mod^{\K}_\Sigma(\G,\X)$. En fait, $T_1:\Mod^{\K'}_\Sigma(\G,\X)\fl\Mod^{\K}_\Sigma(\G,\X)$ est un revêtement fini étale galoisien de groupe $\K/\K'$. 

\end{proposition}

\section{Compactification de Baily-Borel }

Le but de ce paragraphe est de montrer que, pour $(\K,\Sigma)$ comme plus haut, quitte à ajouter à $\Sigma$ un nombre fini de nombres premiers, le modèle entier $\Mod^{\K}_\Sigma(\G,\X)$ de $M^{\K}(\G,\X)$ s'étend en un ``bon'' modèle de la compactification de Baily-Borel $M^{\K}(\G,\X)^*$. 
On veut en particulier que le modèle entier de la compactification ait une stratification dont les strates sont les modèles entiers des strates de $M^{\K}(\G,\X)^*$.
Pour que cela soit possible, on a besoin du lemme suivant.

\begin{lemme}\label{strat_mod_p} Soit $(\K,\Sigma)$ comme plus haut, c'est-à-dire tel que
$\K=\K_\Sigma\K_0^\Sigma$, avec $\K_0^\Sigma=\G(\Z^\Sigma)$ et $\K_\Sigma\subset\G({\Af}_{\Sigma})$,
et tel que $\Sigma$ contienne les diviseurs premiers du discriminant $D$ de $E$. 
Soit $\P_r$, $r\in\{1,\dots,q\}$ un sous-groupe parabolique maximal standard de $\G$.
On rappelle qu'on a défini dans 1.4.1 et 1.4.2 des sous-groupes distingués
$\G_r\subset\G'_r$ et $\L'_{\lin,r}\subset\L_{\lin,r}$ du quotient de Levi $\L_r=\P_r/\N_r$ de $\P_r$
tels que $\L_r$ soit produit direct de $\G'_r$ et $\L'_{\lin,r}$ et produit quasi-direct de $\G_r$ et $\L_{\lin,r}$,
que $\G_r$ soit le centre de $\G'_r$,
que $\L'_{\lin,r}=\L_{\lin,r}$ et $\G'_r=\G_r$ si $r\not=q$
et que $\L_{\lin,q}$ soit isomorphe à $\L'_{\lin,q}\times\SU(p-q)$.

\begin{itemize} 
\item[(i)] Si $r<q$, on peut trouver un système de représentants du double quotient \newline
$\P_r(\Q)\QP_r(\Af)\sous\G(\Af)/\K$ dans $(\L_{\lin,r}\N_r)(\Af)\G({\Af}_\Sigma)$.
Si $r=q$, on peut trouver un système de représentants du double quotient ci-dessus
dans $\P_r(\Af)\G({\Af}_\Sigma)$.
\item[(ii)] Soit $\RP$ un sous-groupe parabolique de $\P$ tel que $\QP_r\subset\RP$.
On note $\N_R$ le radical unipotent de $\RP$ 
et $\RP_\lin$ le sous-groupe parabolique de $\L_{\lin,r}$ correspondant à $\RP$ 
($\RP/\N_r$ est produit quasi-direct de $\RP_\lin$ et $\G_r$).
Si $g\in (\L_{\lin,r}\N_r)(\Af)\G({\Af}_\Sigma)$, et si $\K_G$ est l'un des sous-groupes ouverts compact de $\G_r(\Af)$ ci-dessous
$$(g\K g^{-1}\cap\QP_r(\Af))/(g\K g^{-1}\cap\N_r(\Af))$$
$$(g\K g^{-1}\cap\RP(\Q)\QP_r(\Af))/(g\K g^{-1}\cap\RP_\lin(\Q)\N_r(\Af))$$
$$(g\K g^{-1}\cap\N_R(\Q)\QP_r(\Af))/(g\K g^{-1}\cap\N_R(\Q)\N_r(\Af))$$
alors $\K_G$ s'écrit
$\K_G=\K_{G,\Sigma}\G_r(\Z^\Sigma)$,
avec $\K_{G,\Sigma}\subset\G_r({\Af}_\Sigma)$.

Si $r=q$, cette propriété reste vraie pour $g\in\P_r(\Af)\G({\Af}_\Sigma)$.
\end{itemize}
\end{lemme}

\begin{preuve} 
\begin{itemize}
\item[(i)] Pour tout $p\not\in \Sigma$, on a la décomposition d'Iwasawa
$$\G(\Q_p)=\P_r(\Q_p)\K_p.$$
On en déduit que 
$$\G(\Af^\Sigma)=\P_r(\Af^\Sigma)\K_0^\Sigma,$$
donc qu'il existe un système de représentants dans 
$\P_r(\Af)\G({\Af}_\Sigma)$ du double quotient $\P_r(\Q)\QP_r(\Af)\sous\G(\Af)/\K$.
Si $r<q$, $\P_r(\Af)=\P_{\lin,r}(\Af)\QP_r(\Af)$, donc on peut prendre le système de représentants 
dans $(\P_{\lin,r}\N_r)(\Af)\G({\Af}_\Sigma)$.

\item[(ii)] Soit $g\in (\L_{\lin,r}\N_r)(\Af^\Sigma)\G({\Af}_\Sigma)$.

Il suffit de montrer que
$$\K_1=(g\K_0^\Sigma g^{-1}\cap\QP_r(\Af^\Sigma))/(g\K_0^\Sigma g^{-1}\cap\N_r(\Af^\Sigma))$$
et
$$\K_2=(g\K_0^\Sigma g^{-1}\cap\RP(\Q)\QP_r(\Af^\Sigma))/(g\K_0^\Sigma g^{-1}\cap\RP_\lin(\Q)\N_r(\Af^\Sigma))$$
sont tous les deux égaux à $\G_r(\Z^\Sigma)$. 
Comme $\K_1\subset\K_2$, on va montrer que $\G_r(\Z^\Sigma)\subset\K_1$ et $\K_2\subset\G_r(\Z^\Sigma)$.

Supposons d'abord que $r<q$.
Alors il suffit de traiter le cas où $g\in (\L_{\lin,r}\N_r)(\Af^\Sigma)$.
Notons $\pi$ la projection de $\P_r(\Af^\Sigma)$ sur $\G_r(\Af^\Sigma)$; alors
$$\K_1=\pi(g\K_0^\Sigma g^{-1}\cap\QP_r(\Af^\Sigma))$$
$$\K_2=\pi(g\K_0^\Sigma g^{-1}\cap\RP(\Q)\QP_r(\Af^\Sigma)).$$
Soit $h\in\G_r(\Z^\Sigma)$. 
Il existe $h'\in\QP_r(\Af^\Sigma)\cap\K_0^\Sigma=\QP_r(\Z^\Sigma)$ tel que $\pi(h')=h$. 
Alors $gh' g^{-1}\in g\K_0^\Sigma g^{-1}\cap\QP_r(\Af^\Sigma)$, et $\pi(g h' g^{-1})=h$ (car $\pi(g)=1$). Donc $\G_r(\Z^\Sigma)\subset\K_1$.\newline
Soit $k\in\K_0^\Sigma$ tel que $g kg^{-1}\in g\K_0^\Sigma g^{-1}\cap\RP(\Q)\QP_r(\Af^\Sigma)$. 
Alors $k\in\P_r(\Af^\Sigma)\cap\K_0^\Sigma$ (car $g\in\P_r(\Af^\Sigma)$), donc $\pi(k)=\pi(gkg^{-1})\in\G_r(\Z^\Sigma)$. Donc $\K_2\subset\G_r(\Z^\Sigma)$.

Traitons enfin le cas $r=q$.
On peut supposer que $g\in\P_q(\Af^\Sigma)$.
$\K_1$ est le projeté sur $\L_q(\Af^\Sigma)$ de $g(\K_0^\Sigma\cap\QP_q(\Af^\Sigma))g^{-1}=g\QP_q(\Z\Sigma)g^{-1}$.
Comme $\G_q$ est central dans $\L_q$, $\K_1$ contient $\G_q(\Z^\Sigma)$.
D'autre part, on sait que $\L_q=\L'_{\lin,q}\times\G'_q$.
Notons $\pi$ la projection de $\P_q$ sur $\G'_q$.
On a 
$$\K_2=(g\K_0^\Sigma g^{-1}\cap\L'_{\lin,q}(\Q)\QP_q(\Af^\Sigma))/(g\K_0^{\Sigma}g^{-1}\cap\L'_{\lin,q}(\Q)\N_q(\Af^\Sigma))=\pi(g\K_0^\Sigma g^{-1}\cap\L'_{\lin,q}(\Q)\QP_q(\Af^\Sigma)).$$
Soit $k\in K_0^\Sigma$ tel que $gk g^{-1}\in g\K_0^\Sigma g^{-1}\cap\L'_{\lin,q}(\Q)\QP_q(\Af^\Sigma)$.
Alors $k\in\P_q(\Af^\Sigma)$, 
et $\pi(gkg^{-1})=\pi(g)\pi(k)\pi(g)^{-1}=\pi(k)$,
car $\pi(k)\in\G_q(\Z^\Sigma)$ et $\G_q$ est le centre de $\G'_q$.
Donc $\K_2\subset\G_q(\Z^\Sigma)$.

\end{itemize}
\end{preuve}

On voudrait de plus que tous les modèles entiers qui apparaissent soient des espaces algébriques. C'est possible grâce au lemme suivant.

\begin{lemme}\label{K_petit} Soit $(\K,\Sigma)$ comme dans la définition \ref{def:struct_niveau}. Rappelons qu'on a fixé un système $\P_1,\dots,\P_q$ de sous-groupes paraboliques maximaux standard de $\G$.
Alors il existe un sous-groupe ouvert compact $\K'\subset\K$, un ensemble de nombres premiers $\Sigma'\supset \Sigma$ et, pour tout $i\in\{1,\dots,q\}$, un système de représentants $(g_{ij})_{i\in J_i}$ du double quotient\newline
 $\P_i(\Q)\QP_i(\Af)\sous\G(\Af)/K'$  tels que
\begin{itemize}
\item[(i)] $\Sigma'-\Sigma$ est fini;
\item[(ii)] $(\K',\Sigma')$ vérifie la condition de la définition \ref{def:struct_niveau};
\item[(iii)] les $g_{ij}$ vérifient la condition du (ii) du lemme \ref{strat_mod_p}; 
\item[(iv)] pour tout $i\in\{1,\dots,q\}$, pour tout $j\in J_i$, le sous-groupe compact ouvert
$$\K'_{ij}= (g_{ij}\K' g_{ij}^{-1}\cap\P_i(\Q)\QP_i(\Af))/(g_{ij}\K' g_{ij}^{-1}\cap\L_{\lin,i}(\Q)\N_i(\Af))$$
de $\G_i(\Af)$ vérifie la condition suivante, qu'on notera (*) : pour tout sous-groupe compact ouvert $\K_G\subset\K'_{ij}$
tel que $(\K_G,\Sigma')$ vérifie la condition de la définition \ref{def:struct_niveau},
le champ $\Mod^{\K_G}_{\Sigma'}(\G_i,\X_i)$ est représentable par un espace algébrique.

\end{itemize}

En particulier, si $\RP\subset\P_i$ est un sous-groupe parabolique de $\P_i$ tel que $\QP_i\subset\RP$ 
et si $\K_G$ est un des sous-groupes ouverts compacts de $\G_i(\Af)$ qui sont associés à $\RP$ comme dans le lemme \ref{strat_mod_p}, le champ $\Mod^{\K_G}_{\Sigma'}(\G_i,\X_i)$ est représentable par un espace algébrique.

\end{lemme}

\begin{preuve} Pour tout $i\in\{1,\dots,q\}$ soit $(g_{ij})_{j\in J_i}$ un système de représentants du double quotient $\P_i(\Q)\QP_i(\Af)\sous\G(\Af)/\K$ vérifiant la condition du (i) lemme \ref{strat_mod_p}. 
Si $\K'$ est un sous-groupe compact ouvert de $\K$, on note 
$$\K'_{ij}=(g_{ij}\K' g_{ij}^{-1}\cap\P_i(\Q)\QP_i(\Af))/(g_{ij}\K' g_{ij}^{-1}\cap\L_{\lin,i}(\Q)\N_i(\Af)).$$
Si $\K'$ est assez petit, il existe pour tous $i,j$ un entier $N\geq 3$ tel que
$$\K'_{ij}\subset\{g\in\G_i(\Af)\mbox{ tq }(g-1)(\O_E^{p+q-2i}\otimes_\Z\widehat{\Z})\subset N(\O_E^{p+q-2i}\otimes_\Z\widehat{\Z})\}.$$
Fixons un $\K'\subset\K$ assez petit au sens ci-dessus et distingué dans $\K$. 
Quitte à diminuer encore $\K'$, on peut trouver un ensemble fini $\Sigma_0$ de nombres premiers tel que $\K'=\K'_{\Sigma_0}\K_0^{\Sigma_0}$, avec $\K'_{\Sigma_0}\subset\G({\Af}_{\Sigma_0})$. On pose $\Sigma'=\Sigma\cup \Sigma_0$.
Pour tous $i,j$, $(\K'_{ij},\Sigma')$ vérifie la condition de la définition \ref{def:struct_niveau}, et $\K'_{ij}$ vérifie la condition (*).
Pour tout $i$, on peut choisir un système de représentants de $\P_i(\Q)\QP_i(\Af)\sous\G(\Af)/\K'$ dans $\displaystyle{\bigcup_{j\in J_i}g_{ij}\K}$. Comme $\K'$ est distingué dans $\K$, ceci finit la démonstration du lemme.

\end{preuve}

\begin{remarque}\label{rq:K_petit} Plus généralement, la méthode de la preuve ci-dessus permet de montrer qu'on peut rendre les sous-groupes compacts ouverts des $\G_i(\Af)$ associés aux strates de bords aussi petits que l'on veut (en les nombres premiers de $\Sigma$) en diminuant $\K$ (ce n'est pas entièrement évident, car le nombre de strates de bord augmente lorsque $\K$ diminue).

\end{remarque}

Soit $(\K,\Sigma)$ comme plus haut; on suppose que $\Mod^{\K}_\Sigma(\G,\X)$ est un espace algébrique. On pose
$$\omega=\bigwedge^{pq}\Omega_{\Mod^{\K}_\Sigma(\G,\X)/\O_F[1/\Sigma]}.$$
C'est un faisceau inversible sur $\Mod^{\K}_\Sigma(\G,\X)$, et 
il existe $N\in\Nat^*$ tel que $\omega^{\otimes N}$ soit très ample sur la fibre générique $M^{\K}(\G,\X)$ ([BB] 10.11).
On a un morphisme (injectif sur les fibres génériques)
$$\Mod^{\K}_\Sigma(\G,\X)\fl P={\bf Proj}\bigoplus_{m\in\Nat}\omega^{\otimes mN}.$$

\begin{definition} On note 
$$\Mod^{\K}_\Sigma(\G,\X)^*$$
l'adhérence de l'image de $\Mod^{\K}_\Sigma(\G,\X)$ dans $P$, et
$$j:\Mod^{\K}_\Sigma(\G,\X)\fl\Mod^{\K}_\Sigma(\G,\X)^*$$
le morphisme évident.
\end{definition}

La fibre générique de $\Mod^{\K}_\Sigma(\G,\X)^*$ est canoniquement isomorphe à $M^{\K}(\G,\X)^*$, et la restriction de $j$ aux fibres génériques est l'immersion ouverte $M^{\K}(\G,\X)\fl M^{\K}(\G,\X)^*$ (mais $j$ lui-même n'est pas forcément une immersion).

Soit $(\K,\Sigma)$ comme plus haut. On suppose que $\K$ vérifie les mêmes conditions que le $\K'$ du lemme \ref{K_petit}.
Pour tout $i\in\{1,\dots,q\}$, soit $(g_{ij})_{j\in J_i}$ un système de représentants du double quotient $\P_i(\Q)\QP_i(\Af)\sous\G(\Af)/\K$ comme dans le lemme \ref{K_petit}. 
On a construit dans 1.4.2 des sous-groupes compacts ouverts $\K_{ij}\subset\G_i(\Af)$ et une stratification 
$$M^{\K}(\G,\X)^*=M^{\K}(\G,\X)\sqcup\coprod_{i,j}M_{ij}$$
telle que $M_{ij}$ s'identifie à la variété de Shimura $M^{\K_{ij}}(\G_i,\X_i)$. 
De plus, on a choisi les $g_{ij}$ pour que $(\K_{ij},\Sigma)$ vérifie la condition de la définition \ref{def:struct_niveau}, donc $M^{\K_{ij}}(\G_i,\X_i)$ a un modèle $\Mod^{\K_{ij}}_\Sigma(\G_i,\X_i)$ sur $\O_F[1/\Sigma]$, qui est un espace algébrique.

On va définir, par récurrence descendante sur $i$, des sous-espaces algébriques localement fermés $\Mod_{ij}$ de $\Mod^{\K}_\Sigma(\G,\X)^*$ dont les fibres génériques sont les $M_{ij}$.
Remarquons que, pour tout $i$, les $M_{ij}$ sont fermés dans
$$M^{\K}(\G,\X)^*-\left(M^{\K}(\G,\X)\cup\bigcup_{i'>i,j'}M_{i'j'}\right).$$
Ceci suggère une manière de faire. Pour tout $j\in J_q$, on prend pour $\Mod_{qj}$ l'adhérence de $M_{qj}$ dans $\Mod^{\K}_\Sigma(\G,\X)^*-j(\Mod^{\K}_\Sigma(\G,\X))$, munie de la structure réduite.
Soit $i\in\{1,\dots,q-1\}$, et supposons les $\Mod_{i'j'}$ définis pour tout $i'>i$. Alors, pour tout $j\in J_i$, on prend pour $\Mod_{ij}$ l'adhérence (munie de la structure réduite) de $M_{ij}$ dans
$$\Mod^{\K}_\Sigma(\G,\X)^*-\left(j(\Mod^{\K}_\Sigma(\G,\X))\cup\bigcup_{i'>i,j'}\Mod_{i'j'}\right).$$

\begin{proposition}\label{modele_entier_BB} Si on remplace $\Sigma$ par son union avec un ensemble fini de nombres premiers 
(ce qui revient à inverser un nombre fini de nombres premiers sur la base $\O_F[1/\Sigma]$), 
alors $\Mod^{\K}_\Sigma(\G,\X)^*$ est un schéma normal projectif sur $\O_F[1/\Sigma]$, 
$j:\Mod^{\K}_\Sigma(\G,\X)\fl\Mod^{\K}_\Sigma(\G,\X)^*$ est une immersion ouverte, $(\Mod_{ij})$ est une stratification du bord $\Mod^{\K}_\Sigma(\G,\X)^*-\Mod^{\K}_\Sigma(\G,\X)$ 
et, pour tous $i,j$, l'isomorphisme $M_{ij}\simeq M^{\K_{ij}}(\G_i,\X_i)$ se prolonge en un isomorphisme $\Mod_{ij}\simeq\Mod^{\K_{ij}}_\Sigma(\G_i,\X_i)$.

\end{proposition}

\section{Systèmes de coefficients}

Dans ce paragraphe, on va étendre les complexes de la définition \ref{def_faisceaux} aux modèles entiers. 
On a donc un groupe connexe $\L$ qui est produit quasi-direct de deux sous-groupes
$\G$ et $\G_\lin$ et une donnée de Shimura $(\G,\X)$, et l'action de $\G(\Q)$
sur $\X$ s'étend en une action de $\L(\Q)$ telle que $\G_\lin(\Q)$ agisse trivialement.
Soit $\K_L\subset\L(\Af)$ un sous-groupe ouvert compact net.
On note $\Hr=\K_L\cap\L(\Q)\G(\Af)$, $\Gamma_\lin=\K_L\cap\G_\lin(\Q)$.
Dans les cas qui nous intéressent ici, $\Hr/\Gamma_\lin$ s'identifie à un sous-groupe ouvert compact de $\G(\Af)$,
et $M(\G,\X)/\Hr=M^{\Hr/\Gamma_\lin}(\G,\X)$.

Soit $\Sigma$ un ensemble de nombres premiers.
On suppose que $\Hr/\Gamma_\lin$ s'écrit $\G_r(\Z^\Sigma)\K_{\Sigma}$, avec $\K_{\Sigma}\subset\G({\Af}_\Sigma)$.
Alors $(\Hr/\Gamma_\lin,\Sigma)$ vérifie la condition de la définition \ref{def:struct_niveau}, donc le modèle entier $\Mod_\Sigma^{\Hr/\Gamma_\lin}(\G,\X)$ existe.

\begin{proposition}\label{lissite_faisceaux} On fixe un nombre premier $\ell\in \Sigma$.
Alors les complexes de faisceaux $\ell$-adiques lisses $\F^{\Hr/\Gamma_\lin}R\Gamma(\Gamma_\lin,V)$
de la définition \ref{def_faisceaux} 
(resp. les faisceaux $\ell$-adiques lisses $\F^{\Hr/\Gamma_\lin}A$ de la remarque \ref{rq:def_faisceaux}) 
s'étendent en des complexes de faisceaux $\ell$-adiques lisses (resp. des faisceaux $\ell$-adiques lisses) sur $\Mod^{\Hr/\Gamma_\lin}_\Sigma(\G,\X)$, qu'on désignera par les mêmes notations.

\end{proposition}

\begin{preuve} Notons $pr$ la projection $\L(\Af)\fl\L(\Q_\ell)$.
Les faisceaux $\ell$-adiques lisses sur $M^{\Hr/\Gamma_\lin}(\G,\X)$ qu'on cherche à prolonger sont définis par la méthode de Pink, 
en utilisant le système projectif de revêtements finis étales $M^{\Hr'/\Gamma'_\lin}(\G,\X)\fl M^{\Hr/\Gamma_\lin}(\G,\X)$, 
où $\K'_L\subset\K_L$, $\Hr'=\K'_L\cap\L(\Q)\G(\Af)$ et $\Gamma'_\lin=\K'_L\cap\G_\lin(\Q)$.
Comme tous les $\Hr$-modules auxquels on applique la méthode de Pink proviennent de $pr(\Hr)$-modules, 
on n'utilise en fait que les revêtements $M^{\Hr'/\Gamma'_\lin}(\G,\X)\fl M^{\Hr/\Gamma_\lin}(\G,\X)$ tels que $\K'_L\supset\K_L\cap ker(pr)$. 
Or, si $\K'_L\supset\K_L\cap ker(pr)$, 
$\Hr'/\Gamma'_\lin$ s'écrit $\G(\Z^\Sigma)\K'_{\Sigma}$, avec $\K'_{\Sigma}\subset\G({\Af}_\Sigma)$,
 donc le revêtement fini étale $M^{\Hr'/\Gamma'_\lin}(\G,\X)\fl M^{\Hr/\Gamma_\lin}(\G,\X)$ s'étend en un revêtement fini étale
 $\Mod_\Sigma^{\Hr'/\Gamma'_\lin}(\G,\X)\fl\Mod_\Sigma^{\Hr/\Gamma_\lin}(\G,\X)$.
On peut donc, en appliquant la méthode de Pink à ces revêtements, construire des faisceaux lisses $\ell$-adiques sur $\Mod^{\Hr/\Gamma_\lin}_\Sigma(\G,\X)$ qui prolongent ceux déjà construits sur $M^{\Hr/\Gamma_\lin}(\G,\X)$.

\end{preuve}

\section{Poids}

On s'intéresse maintenant aux poids des faisceaux obtenus.

On identifie le groupe des caractères de $\Gm$ avec $\Z$ en faisant correspondre à $m\in\Z$ le caractère $\chi_m:x\fle x^m$. 
Soient $\H$ un groupe algébrique, $\lambda:\Gm\fl\H$ un cocaractère et $\rho:\H(\Q_\ell)\fl\GL(V)$ une représentation algébrique de $\H(\Q_\ell)$ dans un $\Q_\ell$-espace vectoriel de dimension finie. Alors la représentation $\rho\circ \lambda$ de $\Gm$ s'écrit comme une somme de caractères :
$$\rho\circ w=\bigoplus_{m\in\Z}\chi_m^{\otimes a_m},$$
et les poids de $\rho$ (ou $V$) relativement à $\lambda$ sont les $m$ tels que $a_m\not=0$. 
On dit que la représentation $\rho$ est pure si elle n'a qu'un seul poids. Si $\lambda$ est un cocaractère central de $\H$ (c'est le cas qui nous intéresse), alors toute représentation de $\H(\Q_\ell)$ est somme de représentations pures.

Rappelons qu'on a un cocaractère central $h_0\circ w:\Gm\fl\G,\lambda\fle\lambda I_n$ (le cocaractère $w:\Gm\fl \SD$ est l'inclusion $\R^\times\fl\C^\times$ sur les points réels). Dans la suite, les poids des représentations de $\G$ seront toujours relativement à ce cocaractère.
Le résultat suivant a été prouvé par Pink ([P2] 5.6.2, voir aussi [LR2] 6) :

\begin{proposition}\label{poids_faisceaux} On suppose ici que $\L=\G$. Soient $(\K,\Sigma)$ et $\ell$ comme ci-dessus. 
Alors, quitte à remplacer $\Sigma$ par sa réunion avec un ensemble fini de nombres premiers, on a le résultat suivant : 
pour tout $p\not\in \Sigma$, 
pour toute représentation algébrique de $\G(\Q_\ell)$ dans un $\Q_\ell$-espace vectoriel $V$ de dimension finie, si $V$ est pure de poids $w$, 
le faisceau lisse $\F^{\K}(V)$ sur la réduction modulo $p$ de $\Mod^{\K}_\Sigma(\G,\X)$ est pur de poids $-w$ (au sens de [BBD 5]).

\end{proposition}

\begin{corollaire}\label{poids_faisceaux_generalises} On se place dans la situation 
de la section 3.3, et on prend $(\K_L,\Sigma)$ et $\ell$ comme ci-dessus.
On fixe $p\not\in \Sigma$, et on travaille sur la réduction modulo $p$ des variétés de Shimura.
Comme $\L=\G\G_\lin$ et que $\G$ et $\G_\lin$ commutent, on peut voir $h_0\circ w$ comme un cocaractère central de $\L$, et les poids des représentations de $\L$ seront pris relativement à ce cocaractère.
Soit $V$ une représentation de $\L(\Q_\ell)$ dans un $\Q_\ell$-espace vectoriel de dimension finie qui est pure de poids $w$. 
On note $C=\F^{\Hr/\Gamma_\lin}R\Gamma(\Gamma_\lin,V)$, vu comme un complexe de faisceaux $\ell$-adiques sur $M^{\Hr/\Gamma_\lin}(\G,\X)$.
Alors on a un isomorphisme (non canonique)
$$C\simeq\bigoplus_{i\in\Z} H^i(C)[-i],$$
et, pour tout $i\in\Z$, le faisceau $H^i(C)=\F^{\Hr/\Gamma_\lin}H^i(\Gamma_\lin,V)$ est lisse pur de poids $-w$ sur $M^{\Hr/\Gamma_\lin}(\G,\X)$.

\end{corollaire}

\begin{preuve} Notons $\K=\Hr\cap\G(\Af)$.
$\Hr/\K\Gamma_\lin$ est un groupe fini, donc on un revêtement fini étale
$M^{\K}(\G,\X)\fl M^{\Hr/\Gamma_\lin}(\G,\X)$,
et l'image inverse de $C$ à $M^{\K}(\G,\X)$ est le complexe $\ell$-adique associé 
au complexe $R\Gamma(\Gamma_\lin,V)$ de représentations du groupe réductif $\G(\Q_\ell)$.
L'isomorphisme 
$$C\simeq\bigoplus_{i\in\Z}H^i(C)[-i]$$
est une conséquence du premier des lemmes ci-dessous, appliqué au morphisme fini étale $M^{\K}(\G,\X)\fl M^{\Hr/\Gamma_\lin}(\G,\X)$.
L'isomorphisme canonique
$$H^i(C)=\F^{\Hr/\Gamma_\lin}H^i(\Gamma_\lin,V)$$
vient de l'exactitude du foncteur $\F^{\Hr/\Gamma_\lin}$, et la lissité de $H^i(C)$ de la proposition \ref{lissite_faisceaux}.
Enfin, le fait que les $H^i(C)$ soient purs de poids $-w$ découle de la proposition ci-dessus et du deuxième des lemmes ci-dessous, appliqué à $M^{\K}(\G,\X)\fl M^{\Hr/\Gamma_\lin}(\G,\X)$.

\end{preuve}

\begin{lemme} Soient $\varphi:X\fl Y$ un morphisme fini étale entre schémas de type fini sur un corps
et $K\in D^b_c(Y,\Q_\ell)$. S'il existe un isomorphisme
$$\varphi^*K\simeq\bigoplus_{i\in\Z}H^i(\varphi^*K)[-i]=\bigoplus_{i\in\Z}\varphi^*H^i(K)[-i],$$
alors on a un isomorphisme
$$K\simeq\bigoplus_{i\in\Z}H^i(K)[-i].$$

\end{lemme}

\begin{preuve} Supposons que
$$\varphi^*K\simeq\bigoplus_{i\in\Z}\varphi^*H^i(K)[-i].$$
Pour tout $i\in\Z$, on a un triangle distingué
$$\tau_{\leq i-1}K\fl \tau_{\leq i}K\fl H^i(K)[-i]\stackrel{+1}{\fl}.$$
Si on montre que la flèche de degré $1$ de ce diagramme est nulle pour tout $i\in\Z$, cela impliquera que l'on a des isomorphismes 
$\tau_{\leq i}K\simeq \tau_{\leq i-1}K\oplus H^i(K)[-i]$, 
et la conclusion du lemme en découlera par récurrence sur le cardinal de $\{i\in\Z,H^i(K)\not=0\}$.

Soit donc $i\in\Z$. L'image par $\varphi^*$ de la flèche de degré $1$ du triangle distingué
$$\tau_{\leq i-1}K\fl\tau_{\leq i}K\fl H^i(K)[-i]\stackrel{+1}{\fl}$$
est la flèche de degré $1$ du triangle distingué
$$\tau_{\leq i-1}\varphi^*K\fl \tau_{\leq i}\varphi^*K\fl H^i(\varphi^*K)[-i]\stackrel{+1}{\fl}.$$
Cette dernière flèche est nulle, car l'hypothèse implique que $\tau_{\leq i}\varphi^*K\simeq \tau_{\leq i-1}\varphi^*\oplus H^i(\varphi^*K)[-i]$. D'après le sous-lemme ci-dessous, le morphisme
$$\Ext^1(H^i(K)[-i],\tau_{\leq i-1}K)\fl \Ext^1(H^i(\varphi^*K)[-i],\tau_{\leq i-1}\varphi^*K)$$
induit par $\varphi^*$ est injectif, d'où le résultat cherché.

\end{preuve}

\begin{souslemme} Soit $\varphi:X\fl Y$ un morphisme fini étale. Alors, pour tous $K,L\in D^b_c(Y,\Q_\ell)$ et pour tout $k\in\Z$, le morphisme induit par $\varphi^*$
$$\Ext^k(K,L)\fl \Ext^k(\varphi^*K,\varphi^*L)$$
est injectif.

\end{souslemme}

\begin{preuve} Soient $K,L\in D^b_c(Y,\Q_\ell)$ et $k\in\Z$. Le morphisme $\Ext^k(K,L)\fl \Ext^k(\varphi^*K,\varphi^*L)$ est le composé du morphisme $\Ext^k(K,L)\fl \Ext^k(K,\varphi_*\varphi^*L)$ provenant du morphisme d'adjonction $L\fl\varphi_*\varphi^*L$ et de l'isomorphisme d'adjonction $\Ext^k(K,\varphi_*\varphi^*L)\simeq \Ext^k(\varphi^*K,\varphi^*L)$. Or on a le morphisme trace $Tr:\varphi_*\varphi^*L\fl L$, qui est tel que le composé
$$L\stackrel{adj}{\fl}\varphi_*\varphi^*L\stackrel{Tr}{\fl}L$$
soit un isomorphisme (puisqu'on a pris des complexes à coefficients dans $\Q_\ell$). Le morphisme d'adjonction $L\fl\varphi_*\varphi^*L$ identifie donc $L$ à un facteur direct de $\varphi_*\varphi^*L$, ce qui implique que le morphisme $\Ext^k(K,L)\fl \Ext^k(K,\varphi_*\varphi^*L)$ est injectif.

\end{preuve}

\begin{lemme}\label{fp_lemme4} Soit $\varphi:X\fl Y$ un morphisme fini étale entre deux schémas lisses de type fini sur un corps fini et $\F$ un faisceau $\ell$-adique lisse sur $Y$. Si $\varphi^*\F$ est pur de poids $w$, alors $\F$ est pur de poids $w$.

\end{lemme}

\begin{preuve} Supposons $\varphi^*\F$ pur de poids $w$. D'après [BBD] 4.1.3 et 5.1.14, $\varphi_*\varphi^*\F$ est pur de poids $w$.
Or on a comme dans la preuve du sous-lemme ci-dessus le morphisme trace $Tr:\varphi_*\varphi^*\F\fl\F$, qui est tel que
$$\F\stackrel{adj}{\fl}\varphi_*\varphi^*\F\stackrel{Tr}{\fl}\F$$
soit un isomorphisme. $\F$ s'identifie donc à un facteur direct de $\varphi_*\varphi^*\F$, donc, d'après [BBD] 5.3.1, il est pur de poids $w$.

\end{preuve}

\section{Théorème de Pink}

On se place, comme dans les deux sections ci-dessus, dans la situation de la définition \ref{def_faisceaux}. 
Soit $\Sigma$ un ensemble de nombres premiers.
On suppose que 
$(\Hr/\Gamma_\lin,\Sigma)$ est comme dans la proposition \ref{modele_entier_BB}, dont on utilisera les notations (c'est-à-dire que $\Mod^{\Hr/\Gamma_\lin}_\Sigma(\G,\X)$ a une ``bonne'' compactification de Baily-Borel).
On note $j$ l'immersion ouverte $\Mod^{\Hr/\Gamma_\lin}_\Sigma(\G,\X)\fl\Mod^{\Hr/\Gamma_\lin}_\Sigma(\G,\X)^*$.

On fixe $\ell\in \Sigma$ et on note $pr$ la projection $\L(\Af)\fl\L(\Q_\ell)$.

La proposition suivante et son corollaire ont été montrés par Wildeshaus dans [W] : 

\begin{proposition}\label{prolongement_sur_modeles_entiers} Quitte à remplacer $\Sigma$ par sa réunion avec un ensemble fini de nombres premiers, on a le résultat suivant :
la formation de $Rj_*$ commute à tout changement de base $S\fl Spec(\O_F[1/\Sigma])$, et 
pour tout $\K_L$-module $M$ qui provient d'un $pr(\K_L)$-module annulé par une puissance de $\ell$, pour tout $q\in\Nat$, pour toute strate $\Mod_{ij}$ de $\Mod^{\Hr/\Gamma_\lin}_\Sigma(\G,\X)^*$, 
la restriction de $R^qj_*\F^{\Hr/\Gamma_\lin}R\Gamma(\Gamma_\lin,M)$ à $\Mod_{ij}$ est un faisceau localement constant.

\end{proposition}

\begin{preuve} La première condition et la propriété cherchée peuvent être obtenues pour un nombre fini de $\K_L$-modules en ajoutant à $\Sigma$ un nombre fini de nombres premiers, 
grâce au théorème de changement de base générique de Deligne (SGA 4 1/2 Th finitude), cf [W] 3.8.

Notons $\Cat$ la catégorie des $\K_L$-modules $M$ qui proviennent d'un $pr(\K_L)$-module annulé par une puissance de $\ell$. 
Si la propriété cherchée est vrai pour les objets simples de $\Cat$, alors elle est vraie pour tous les objets de $\Cat$.
Or $\Cat$ n'a qu'un nombre fini d'objets simples ([W] 3.7), d'où la proposition.

\end{preuve}

Comme un faisceau localement constant sur le modèle canonique d'une variété de Shimura a au plus une extension localement constante au modèle entier, on en déduit le :

\begin{corollaire} Dans la situation de la proposition ci-dessus, le théorème \ref{th_Pink} et son corollaire restent vrais sur les modèles entiers, 
à condition de se restreindre aux $\Hr$-modules qui proviennent de $\K_L/(\K_L\cap\G(\Af^\ell))$-modules annulés par une puissance de $\ell$.

\end{corollaire}

Soit $\K\subset\G(\Af)$ un sous-groupe compact ouvert net.
Pour construire les complexes pondérés sur (la réduction modulo $p$ de) la compactification de Baily-Borel de $M^{\K}(\G,\X)$, on utilisera un certain nombre (fini) de variétés de Shimura :
\begin{itemize}
\item[(i)] $M^{\K}(\G,\X)$;
\item[(ii)] les strates de bord de $M^{\K}(\G,\X)^*$, et les revêtements finis étales de ces strates qui correspondent aux sous-groupes ouverts compacts du lemme \ref{strat_mod_p} (ii); 
\item[(iii)] pour toute variété de Shimura de (ii), les strates de bord de 
sa compactification de Baily-Borel, et les revêtements de ces strates analogues aux revêtements de (ii);
\item[(iv)] et ainsi de suite.

\end{itemize} 

\begin{definition}\label{bon_modele_entier} On dira que $\K$ est assez petit s'il existe un ensemble fini de nombres premiers $\Sigma$ et $\ell\in \Sigma$ tels que $(\K,\Sigma)$ vérifie la condition de la définition \ref{def:struct_niveau}, et que tous les modèles entiers des variétés de Shimura énumérées ci-dessus 
(ces modèles existent sur $\O_F[1/\Sigma]$ d'après le lemme \ref{strat_mod_p}) vérifient la conclusion du lemme \ref{K_petit} 
et que les modèles entiers de leurs compactifications de Baily-Borel vérifient les conclusions des propositions \ref{modele_entier_BB}, \ref{poids_faisceaux} et \ref{prolongement_sur_modeles_entiers}.
\end{definition}

Il est toujours possible de rendre $\K$ assez petit dans ce sens (utiliser la remarque \ref{rq:K_petit}).

Si on a un couple $(\K,\Sigma)$ vérifiant toutes les conditions ci-dessus et si $p$ est un nombre premier qui n'appartient pas à $\Sigma$, on peut réduire toutes les variétés modulo $p$ :
on obtient des modèles modulo $p$ des variétés et des complexes de faisceaux $\ell$-adiques à faisceaux de cohomologie lisses purs de poids connus sur les variétés de Shimura, 
dont les prolongements aux compactifications de Baily-Borel vérifient le corollaire \ref{restriction_bord}. 
De plus, d'après la proposition \ref{independance_de_T}, ces objets ne dépendent pas du choix de $\Sigma$.
On les notera souvent de la même façon que les objets sur le corps reflex $F$ (en précisant à l'avance si l'on est en caractéristique $0$ ou $p$).

Dans la suite, lorsque nous parlerons de la réduction modulo $p$ de $M^{\K}(\G,\X)$, nous ferons toujours référence à la situation ci-dessus.

\begin{remarque} Si l'on disposait de compactifications toroïdales des $\Mod_\Sigma^{\K}(\G,\X)$ vérifiant les propriétés de [P2] 3.9 et 3.10,
alors les résultats de [P2] s'étendraient automatiquement aux modèles entiers
(sans qu'il soit besoin d'inverser des nombres premiers supplémentaires).
Malheureusement, on ne sait construire ces compactifications que si le groupe est $\GU(1,1)$ (auquel cas les variétés de Shimura associées sont des courbes modulaires)
ou $\GU(2,1)$ (voir [L] et le chapitre I de [Be]).

\end{remarque}

\newpage

\chapter{Prolongement intermédiaire des faisceaux pervers purs}

Cette partie est indépendante des autres.

On fixe un corps fini $\Fi_q$ et un nombre premier $\ell$ inversible dans $\Fi_q$.
Tous les schémas sont séparés de type fini sur $\Fi_q$.
Si $X$ est un schéma, on note $D^b_m(X,\Q_\ell)$ la catégorie des complexes $\ell$-adiques mixtes sur $X$ (au sens de [BBD] 5.1.5; en particulier, les complexes sont à poids entiers),
munie de la t-structure donnée par la perversité autoduale.

Soient $X$ un schéma et $j:U\fl X$ un ouvert non vide de $X$.
Le premier objectif de cette partie est d'écrire $j_{!*}K$, où $K$ est un faisceau pervers pur sur $U$,
comme un certain tronqué par le poids de $Rj_*K$ (théorème \ref{th:calcul_IC_general}).
Le deuxième objectif est
 de calculer la trace d'une puissance $\Phi$ de l'endomorphisme de Frobenius
sur la cohomologie de $j_{!*}K$ en fonction de
la trace de $\Phi$ sur la cohomologie de $Rj_*K$ et de complexes 
de même type supportés par $X-U$
(théorème \ref{th:simplification_formule_traces}).

\section{Troncature par le poids dans $D^b_m(X,\Q_\ell)$}

\begin{proposition}\label{t_structures:prop1} Soit $X$ un schéma sur $\Fi_q$.
Pour tout $a\in\Z\cup\{\pm\infty\}$, on note $\DP^{\leq a}(X)$, ou $\DP^{\leq a}$ s'il n'y a pas d'ambiguïté sur $X$, (resp. $\DP^{\geq a}(X)$ ou $\DP^{\geq a}$)
la sous-catégorie pleine de $D^b_m(X,\Q_\ell)$ dont les objets sont les complexes mixtes $K$ 
tels que pour tout $i\in\Z$, $\Hp^iK$ soit de poids $\leq a$ (resp. $\geq a$).
Alors:
\begin{itemize}
\item[(i)] $\DP^{\leq a}$ et $\DP^{\geq a}$ sont des sous-catégories stables par décalage et extensions (cf [BBD] 1.2.6) de $D^b_m(X,\Q_\ell)$ (en particulier, ce sont des sous-catégories triangulées).\newline
La dualité de Poincaré échange $\DP^{\leq a}$ et $\DP^{\geq -a}$.\newline
Si $a<a'$, on a $\DP^{\leq a}\cap \DP^{\geq a'}=0$.
\item[(ii)] On a $\DP^{\leq a}(1)=\DP^{\leq a-2}$ et $\DP^{\geq a}(1)=\DP^{\geq a-2}$
(où $(1)$ est le twist à la Tate).
\item[(iii)] Pour tous $K\in\DP^{\leq a}$ et $L\in\DP^{\geq a+1}$, on a
$R\Hom(K,L)=0$.
\item[(iv)] Pour tout $a\in\Z\cup\{\pm\infty\}$, $(\DP^{\leq a},\DP^{\geq a+1})$ est une t-structure sur $D^b_m(X,\Q_\ell)$.
\end{itemize}
\end{proposition}

Cette proposition sera démontrée dans la section 4.2.

D'après [BBD] 1.3.3,
l'inclusion $\DP^{\leq a}\subset D^b_m(X,\Q_\ell)$ (resp. $\DP^{\geq a}\subset D^b_m(X,\Q_\ell)$) 
admet un adjoint à droite (resp. à gauche), qu'on notera $w_{\leq a}$ (resp. $w_{\geq a}$),
et pour tout $K\in D^b_m(X,\Q_\ell)$, il existe un unique morphisme $w_{\geq a+1}K\fl (w_{\leq a})K[1]$ qui fait du triangle suivant un triangle distingué
$$w_{\leq a}K\fl K\fl w_{\geq a+1}K\fl (w_{\leq a}K)[1].$$
Ce triangle est, à isomorphisme unique près, l'unique triangle distingué $A\fl K\fl B\stackrel{+1}{\fl}$ avec $A\in\DP^{\leq a}$ et $B\in\DP^{\geq a+1}$ (toujours d'après [BBD] 1.3.3).

Comme la dualité de Poincaré échange $\DP^{\leq a}$ et $\DP^{\geq -a}$, 
elle échange aussi $w_{\leq a}$ et $w_{\geq -a}$.

\begin{remarques}
\begin{itemize}
\item[(1)] $\DP^{\leq a}$ n'est pas la catégorie des complexes mixtes de poids $\leq a$ :
un complexe mixte $K$ est de poids $\leq a$ \ssi $\Hp^iK$ est de poids $\leq a+i$ pour tout $i\in\Z$ ([BBD] 5.4.1),
et cette condition est différente de la condition qui caractérise les objets de $\DP^{\leq a}$.
\item[(2)] Le décalage de $1$ dans la définition de la t-structure est nécessaire : 
$(\DP^{\leq a},\DP^{\geq a})$ n'est pas une t-structure.
\item[(3)] $(\DP^{\leq a},\DP^{\geq a+1})$ est une t-structure quelque peu étrange : 
son coeur est nul, et elle ne donne donc pas lieu à une théorie cohomologique intéressante.
De plus, $\DP^{\leq a}$ et $\DP^{\geq a+1}$ sont des sous-catégories triangulées de $D^b_m(X,\Q_\ell)$
(en particulier, elles sont stables par le foncteur $[1]$),
ce qui est inhabituel.
\item[(4)] Si $K$ est un faisceau pervers mixte, $w_{\leq a}$ est simplement le plus grand sous-faisceau pervers de $K$ de poids $\leq a$,
et $w_{\geq a}K$ est le plus grand quotient pervers de $K$ de poids $\geq a$ ([BBD] 5.3.5).\newline
Beilinson a montré dans [B] que le foncteur ``réalisation'' ([BBD] 3.1.9) de la catégorie dérivée
bornée de la catégorie des faisceaux pervers mixtes sur $X$ dans $D^b_m(X,\Q_\ell)$
est une équivalence de catégories.
Si $K\in D^ b_m(X,\Q_\ell)$ est représenté par un complexe $C^\bullet$ de faisceaux 
pervers mixtes, alors $w_{\leq a}K$ (resp. $w_{\geq a}K$) est représenté par le complexe
$(w_{\leq a}C^n)$ (resp. $(w_{\geq a}C^n)$).

\end{itemize}
\end{remarques}

La proposition suivante donne quelques propriétés des foncteurs $w_{\leq a}$ et $w_{\geq a}$.

\begin{proposition}\label{t_structures:prop2}
\begin{itemize}
\item[(i)] Pour tout $K\in D^b_m(X,\Q_\ell)$, on a $w_{\leq a}(K(1))=(w_{\leq a+2}K)(1)$ et $w_{> a}(K(1))=(w_{> a+2}K)(1)$.
\item[(ii)] Soit $K\in D^b_m(X,\Q_\ell)$. 
L'image par $\Hp^i$ de la flèche de cobord $w_{\geq a+1}K\fl (w_{\leq a}K)[1]$ est nulle pour tout $i\in\Z$,
donc la suite exacte longue de cohomologie perverse du triangle distingué
$$w_{\leq a}K\fl K\fl w_{\geq a+1}K\stackrel{+1}{\fl}$$
donne des suites exactes courtes de faisceaux pervers
$$0\fl\Hp^iw_{\leq a}K\fl\Hp^iK\fl \Hp^iw_{\geq a+1}K\fl 0.$$
\item[(iii)] $w_{\leq a}$ et $w_{\geq a}$ commutent au foncteur de décalage $[1]$,
et ils envoient les triangles distingués
de $D^b_m(X,\Q_\ell)$ sur des triangles distingués.
\item[(iv)] $w_{\leq a}$ et $w_{\geq a}$ envoient la catégorie abélienne des faisceaux pervers mixtes dans elle-même, 
et leurs restrictions à cette catégorie sont des foncteurs exacts.
Pour tout $K\in D^b_m(X,\Q_\ell)$, pour tout $i\in\Z$, on a
$$w_{\leq a}(\Hp^iK)=\Hp^iw_{\leq a}K$$
$$w_{\geq a}(\Hp^iK)=\Hp^iw_{\geq a}K.$$
\item[(v)] Soit $f:X\fl Y$ un morphisme. 
Si la dimension des fibres de $f$ est inférieure ou égale à $d$, alors
$$Rf_!(\DP^{\leq a}(X))\subset \DP^{\leq a+d}(Y)$$
$$Rf_*(\DP^{\geq a}(X))\subset \DP^{\geq a-d}(Y)$$
$$f^*(\DP^{\leq a}(Y))\subset \DP^{\leq a+d}(X)$$
$$f^!(\DP^{\geq a}(Y))\subset\DP^{\geq a-d}(X).$$

\end{itemize}
\end{proposition}

La démonstration des propositions \ref{t_structures:prop1} et \ref{t_structures:prop2} sera donnée dans la section 4.2.

\begin{theoreme}\label{th:calcul_IC_general} Soient $a\in\Z$, $X$ un schéma séparé 
de type fini sur $\Fi_q$, $j:U\fl X$ un ouvert non vide de $X$,
et $K$ un faisceau pervers pur de poids $a$ sur $U$.
Alors les flèches canoniques
$$w_{\geq a}j_!K\fl j_{!*}K\fl w_{\leq a}Rj_*K$$
sont des isomorphismes.

\end{theoreme}

Si $K$ est le faisceau constant $\Q_\ell$, cette formule est à rapprocher des formules 4.5.7 et 4.5.9
de l'article [S] de Morihiko Saito.

\begin{preuve} Il suffit de montrer que la deuxième flèche est un isomorphisme 
(le cas de la première flèche en résulte par dualité).\newline
Cela découle des trois points suivants ($i$ est l'inclusion de $X-U$ dans $X$) :
\begin{itemize}
\item[(1)] Un complexe $K\in D^b_m(U,\Q_\ell)$ a au plus un prolongement $L\in D^b_m(X,\Q_\ell)$ tel que
$i^*L\in\DP^{\leq a}$ et $i^!L\in\DP^{\geq a+1}$.\newline
En effet, soient $L,L'\in D^b_m(X,\Q_\ell)$.
On a un triangle distingué
$$R\Hom(i^*L,i^!L')\fl R\Hom(L,L')\fl R\Hom(j^*L,j^*L')\stackrel{+1}{\fl}.$$
Si $i^*L\in\DP^{\leq a}$ et $i^!L\in\DP^{\geq a+1}$, alors $R\Hom(i^*L,i^!L')=0$ d'après le (iii) de la proposition \ref{t_structures:prop1},
donc on a un isomorphisme
$$R\Hom(L,L')\iso R\Hom(j^*L,j^*L').$$
\item[(2)] Soit $K\in D^b_m(U,\Q_\ell)$. On note $L=w_{\leq a}Rj_*K$. 
Alors $i^*L\in\DP^{\leq a}$ par le (iv) de la proposition \ref{t_structures:prop2}.
De plus, on a un triangle distingué
$$L\fl Rj_*K\fl w_{\geq a+1}Rj_*K\stackrel{+1}{\fl},$$
d'où un isomorphisme
$$i^!w_{\geq a+1}Rj_*K[-1]\iso i^!L.$$
D'après le (iv) de la proposition \ref{t_structures:prop2}, $i^!L\in\DP^{\geq a+1}$.
\item[(3)] Soit $K$ un faisceau pervers pur de poids $a$ sur $U$.
Alors $j_{!*}K$ est pervers pur de poids $a$ sur $X$ par [BBD] 5.4.3, 
donc, par [BBD] 5.1.14 et 5.4.1, pour tout $k\in\Z$, 
$\Hp^ki^*j_{!*}K$ est de poids $\leq a+k$
et $\Hp^ki^!j_{!*}K$ est de poids $\geq a+k$.
D'après le lemme \ref{fp_lemme1} de la section 5 à la fin de cette partie, 
$\Hp^ki^*j_{!*}K=0$ si $k\geq 0$ et $\Hp^ki^!j_{!*}K=0$ si $k\leq 0$,
donc, pour tout $k\in\Z$,
$\Hp^ki^*j_{!*}K$ est de poids $\leq a-1$ et $\Hp^ik^!j_{!*}K$ est de poids $\geq a+1$.
Autrement dit, $i^*j_{!*}K\in\DP^{\leq a-1}\subset\DP^{\leq a}$ et $i^!j_{!*}K\in\DP^{\geq a+1}$.

\end{itemize}
\end{preuve}

\section{Preuve des propositions de la section 4.1}

Dans cette section, nous allons prouver les propositions \ref{t_structures:prop1} et \ref{t_structures:prop2} de la section précédente.

\begin{preuvepn}{\ref{t_structures:prop1}}
\begin{itemize}
\item[(i)] Il est clair que $\DP^{\leq a}$ et $\DP^{\geq a}$ sont stables par décalage
et que la dualité de Poincaré échange $\DP^{\leq a}$ et $\DP^{\geq -a}$.\newline
Pour montrer la stabilité par extensions de $\DP^{\leq a}$ (resp. $\DP^{\geq a})$,
il suffit de prouver que la catégorie des faisceaux pervers de poids $\leq a$ (resp. $\geq a$)
est une sous-catégorie épaisse de la catégorie des faisceaux pervers,
c'est-à-dire stable par noyaux, conoyaux et extensions.
Par dualité, il suffit de traiter le premier cas.
La stabilité par noyaux et conoyaux résulte de [BBD] 5.3.1,
et la stabilité par extensions se prouve facilement à partir de [BBD] 5.1.9.\newline
Supposons que $a<a'$.
En appliquant la fin de [BBD] 5.1.8 aux objets de cohomologie perverse, on voit que $\DP^{\leq a}\cap\DP^{\geq a'}=0$.
\item[(ii)] Évident.
\item[(iii)] Voir le premier des lemmes ci-dessous. 
\item[(iv)] La condition (ii) de [BBD] 1.3.1 est évidente. 
La condition (i) résulte du point (iii) ci-dessus, et la condition
(iii) est prouvée dans le deuxième des lemmes ci-dessous.
\end{itemize}
\end{preuvepn}

\begin{lemme}\label{lemme_annulation} Soient $X$ un schéma séparé de type fini sur $\Fi_q$, $K,L\in D^b_m(X,\Q_\ell)$ et $a\in\Z$. On suppose que pour tout $i\in\Z$, $\Hp^i(K)$ est de poids $\leq a$ et $\Hp^i(L)$ de poids $\geq a+1$. Alors $R\Hom(K,L)=0.$

\end{lemme}

\begin{preuve} Pour tout $i\in\Z$, on a un triangle distingué
$${}^{p}\tau_{\leq i-1}L\fl {}^{p}\tau_{\leq i}L\fl\Hp^i L[-i]\stackrel{+1}{\fl},$$
d'où un triangle distingué
$$R\Hom(K,{}^{p}\tau_{\leq i-1}L)\fl R\Hom(K,{}^{p}\tau_{\leq i}L)\fl R\Hom(K,\Hp^i L)[-i]\stackrel{+1}{\fl}.$$
Si on montre le résultat pour $L$ pervers, on pourra, grâce à ces triangles, en déduire le résultat pour $L$ quelconque en faisant une récurrence sur le cardinal de $\{i\in\Z\mbox{ tq }\Hp^i L\not=0\}$. 
On peut donc supposer $L$ pervers. On se ramène de même au cas où $K$ est pervers.

Notons $F$ le morphisme de Frobenius géométrique.
D'après [BBD] 5.1.2.5, on a pour tout $i\in\Z$ une suite exacte
$$0\fl \Ext^{i-1}(K_{\overline{\Fi}_q},L_{\overline{\Fi}_q})_F\fl \Ext^i(K,L)\fl \Ext^i(K_{\overline{\Fi}_q},L_{\overline{\Fi}_q})^F\fl 0.$$

$K$ est de poids $\leq a$ et $L$ de poids $\geq a+1$, donc $\Ext^i(K_{\overline{\Fi}_q},L_{\overline{\Fi}_q})$ est de poids $>i$ pour tout $i\in\Z$ ([BBD] 5.1.15 (i)). On en déduit que si $i\geq 0$, $\Ext^i(K_{\overline{\Fi}_q},L_{\overline{\Fi}_q})$ est de poids $>0$, donc que
$$\Ext^i(K_{\overline{\Fi}_q},L_{\overline{\Fi}_q})_F=\Ext^i(K_{\overline{\Fi}_q},L_{\overline{\Fi}_q})^F=0$$
(comme dans la preuve de [BBD] 5.1.15). 

D'autre part, $K$ et $L$ sont pervers (donc $K_{\overline{\Fi}_q}$ et $L_{\overline{\Fi}_q}$ aussi), d'où $\Ext^i(K_{\overline{\Fi}_q},L_{\overline{\Fi}_q})=0$ pour $i<0$. 

Finalement, on a obtenu : pour tout $i\in\Z$,
$$\Ext^i(K_{\overline{\Fi}_q},L_{\overline{\Fi}_q})_F=\Ext^i(K_{\overline{\Fi}_q},L_{\overline{\Fi}_q})^F=0.$$
Les suites exactes ci-dessus impliquent que, pour tout $i\in\Z$,
$$\Ext^i(K,L)=0,$$
ce qui est le résultat cherché.

\end{preuve}

\begin{remarque} Le lemme ci-dessus reste valable,
avec la même preuve, 
pour une perversité arbitraire (vérifiant les conditions de [BBD] 2.2.1).

\end{remarque}

\begin{lemme}\label{lemme:filtration_poids} Soit $X$ un schéma séparé de type fini sur $\Fi_q$.
Alors pour tout $a\in\Z$ et tout $K\in D^b_m(X,\Q_\ell)$, il existe un triangle distingué dans $D^b_m(X,\Q_\ell)$
$$K_1\fl K\fl K_2\stackrel{+1}{\fl}$$
tel que pour tout $i\in\Z$, $\Hp^i K_1$ soit de poids $\leq a$ et $\Hp^i K_2$ de poids $\geq a+1$.

\end{lemme}

\begin{preuve} On fixe $a\in\Z$ et on raisonne par récurrence sur $card(\{i\in\Z\mbox{ tq }\Hp^iK\not=0\})$.

Supposons qu'il existe $i\in\Z$ tel que $K\simeq\Hp^iK[-i]$.
Alors, d'après [BBD] 5.3.5, il existe un sous-faisceau pervers $L\subset\Hp^iK$ de poids $\leq a$ tel que $\Hp^iK/L$ soit de poids $\geq a+1$.
Il suffit de poser $K_1=L[-i]$ et $K_2=(\Hp^iK/L)[-i]$.

Soit $K\in D^b_m(X,\Q_\ell)$ tel que $n=card\{i\in\Z\mbox{ tq }\Hp^iK\not=0\}\geq 2$. 
Supposons le lemme prouvé pour tous les $K'\in D^b_m(X,\Q_\ell)$ tels que $card\{i\in\Z\mbox{ tq }\Hp^iK'\not=0\}<n$, et montrons-le pour $K$.
Il suffit de montrer le résultat suivant : si on a un triangle distingué de $D^b_m(X,\Q_\ell)$
$$K'\fl K\fl K''\stackrel{+1}{\fl}$$
et que la conclusion du lemme vaut pour $K'$ et $K''$, alors elle vaut aussi pour $K$.

On se donne donc des triangles distingués dans $D^b_m(X,\Q_\ell)$
$$\xymatrix{ & & & \\
K'_2\ar[u]^{+1} & & K''_2\ar[u]^{+1} & \\
K'\ar[r]\ar[u] & K\ar[r] & K''\ar[r]^{+1}\ar[u] & \\
K_1'\ar[u] &  & K''_1\ar[u] & }$$
et on suppose que pour tout $i\in\Z$, $\Hp^iK'_1$ et $\Hp^iK''_1$ sont de poids $\leq a$ et $\Hp^iK'_2$ et $\Hp^iK''_2$ de poids $\geq a+1$;
autrement dit, $K'_1$ et $K''_1$ sont dans $\DP^{\leq a}$ et $K'_2$ et $K''_2$ sont dans $\DP^{\geq a+1}$.
D'après le lemme \ref{lemme_annulation}, $R\Hom(K''_1,K'_2[1])=0$, donc il existe un unique morphisme $K''_1\fl K'_1[1]$ qui fait commuter le carré
$$\xymatrix{K''_1\ar[r]\ar[d] & K'_1[1]\ar[d] \\
K''\ar[r] & K'[1]}$$
D'après [BBD] 1.1.11, on peut compléter le diagramme commutatif en traits pleins
pour obtenir un diagramme dont les lignes et les colonnes sont des triangles distingués
et dont tous les carrés sont commutatifs, sauf le carré marqué $-$, qui est anticommutatif :
$$\xymatrix{K''_1[1]\ar[r] & K'_1[2]\ar@{-->}[r] & K_1[2]\ar@{-->}[r]\ar@{}[rd]|{-} & K''_1[2] \\
K''_2\ar@{-->}[r]\ar[u] & K'_2[1]\ar@{-->}[r]\ar[u] & K_2[1]\ar@{-->}[r]\ar@{-->}[u] & K''_2[1]\ar[u] \\
K''\ar[r]\ar[u] & K'[1]\ar[r]\ar[u] & K[1]\ar[r]\ar@{-->}[u] & K''[1]\ar[u] \\
K''_1\ar[u]\ar[r] & K'_1[1]\ar@{-->}[r]\ar[u] & K_1[1]\ar@{-->}[r]\ar@{-->}[u] & K''_1[1]\ar[u]}$$
Comme $\DP^{\leq a}$ et $\DP^{\geq a+1}$ sont stables par extensions (proposition \ref{t_structures:prop1} (i)),
$K_1\in\DP^{\leq a}$ et $K_2\in\DP^{\geq a+1}$.

\end{preuve}

\begin{preuvepn}{\ref{t_structures:prop2}}
\begin{itemize}
\item[(i)] Soit $K\in D^b_m(X,\Q_\ell)$.
On a un triangle distingué
$$(w_{\leq a+2}K)(1)\fl K(1)\fl (w_{>a+2}K)(1)\stackrel{+1}{\fl}$$
avec $(w_{\leq a+2}K)(1)\in\DP^{\leq a}$ et $(w_{>a+2}K)(1)\in\DP^{>a}$, d'où des isomorphismes canoniques
$w_{\leq a}(K(1))=(w_{\leq a+2}K)(1)$ et $w_{>a}(K(1))=(w_{>a+2}K)(1)$.
\item[(ii)] Soient $K\in D^b_m(X,\Q_\ell)$ et $i\in\Z$.
Comme $\Hp^iw_{\geq a+1}K$ est de poids $\geq a+1$ et que $\Hp^{i+1}w_{\leq a}K$ est de poids $\leq a$,
la flèche $\Hp^iw_{\geq a+1}K\fl \Hp^{i+1}w_{\leq a}K$ est nulle par [BBD] 5.3.1.
\item[(iii) et (iv)] Il suffit de prouver les assertions pour $w_{\leq a}$, celles pour $w_{\geq a}$ en résultant par dualité.\newline
$w_{\leq a}$ commute au décalage parce que $\DP^{\leq a}$ est invariante par décalage.\newline
Soit $K$ un faisceau pervers mixte. 
Si $L$ est le plus grand sous-faisceau pervers de $K$ de poids $\leq a$, alors $K/L$ est de poids $\geq a+1$ ([BBD] 5.3.5).
Donc $w_{\leq a}K=L$, et $w_{\leq a}K$ est pervers.\newline
L'exactitude de la restriction de $w_{\leq a}$ à la catégorie des faisceaux pervers mixtes provient de la fin de [BBD] 5.3.5 
(le fait que les morphismes sont strictement compatibles à la filtration par le poids).\newline
Soient $K\in D^b_m(X,\Q_\ell)$ et $i\in\Z$.
D'après (ii), on a une suite exacte de faisceaux pervers
$$0\fl\Hp^iw_{\leq a}K\fl\Hp^i K\fl\Hp^iw_{\geq a+1}K\fl 0$$
avec $\Hp^iw_{\leq a}K$ de poids $\leq a$ et $\Hp^iw_{\geq a+1}K$ de poids $\geq a+1$,
donc $\Hp^iw_{\leq a}K=w_{\leq a}(\Hp^iK)$.\newline
Le fait que $w_{\leq a}$ envoie les triangles distingués sur des triangles distingués
résulte de la preuve du lemme \ref{lemme:filtration_poids}.
\item[(v)] Les inclusions résultent de [BBD] 4.2.4 et 5.1.14.\newline
Traitons par exemple le cas de $Rf_!$.
Soit $K\in \DP^{\leq a}(X)$.
Pour tous $i,j\in\Z$, $\Hp^jK$ est de poids $\leq a$ par définition de $\DP^{\leq a}(X)$,
donc $Rf_!\Hp^jK$ est de poids $\leq a$ par [BBD] 5.1.14,
et $\Hp^i(Rf_!\Hp^jK)$ est de poids $\leq a+i$ par [BBD] 5.4.1.
Or, par [BBD] 4.2.4, $\Hp^i(Rf_!\Hp^jK)=0$ si $i>d$,
donc $\Hp^i(Rf_!\Hp^jK)$ est de poids $\leq a+d$ pour tous $i,j\in\Z$.
En utilisant la suite spectrale
$$E_2^{ij}=\Hp^i(Rf_!\Hp^jK)\Longrightarrow \Hp^{i+j}Rf_!K,$$
on voit que $\Hp^kRf_!K$ est de poids $\leq a+d$ pour tout $k\in\Z$.

\end{itemize}
\end{preuvepn}

\section{t-structures recollées}

Nous utiliserons la notion suivante de stratification :

\begin{definition}\label{definition_stratification} Soit $X$ un schéma séparé de type fini sur $\Fi_q$. 
Une \emph{stratification} de $X$ est une partition finie $(S_i)_{0\leq i\leq n}$ de $X$ par des sous-schémas localement fermés (les strates) 
telle que
pour tout $i\in\{0,\dots,n\}$, $S_{i}$ est ouvert dans
$\displaystyle{X-\bigcup_{0\leq j<i}S_{j}}$.

\end{definition}

Soit $X$ un schéma muni d'une stratification $(S_k)_{0\leq k\leq n}$.
Pour tout $k\in\{0,\dots,n\}$, on note $i_k$ l'inclusion de $S_k$ dans $X$.
$U=S_0$ est un ouvert de $X$; on note $j=i_0:U\fl X$ l'inclusion.

La proposition suivante est un cas particulier de [BBD] 1.4.10.

\begin{proposition} Soit $\as=(a_0,\dots,a_n)\in (\Z\cup\{\pm\infty\})^{n+1}$.
On note $\DP^{\leq\as}(X)$ ou $\DP^{\leq\as}$ (resp. $\DP^{>\as}(X)$ ou $\DP^{>\as}$) la sous-catégorie pleine de $D^b_m(X,\Q_\ell)$ 
dont les objets sont les complexes mixtes $K$ tels que
pour tout $k\in\{0,\dots,n\}$,
$i_k^*K\in\DP^{\leq a_k}(S_k)$ (resp. $i_k^!K\in\DP^{>a_k}(S_k)$).

Alors $(\DP^{\leq\as},\DP^{>\as})$ est une t-structure sur $D^b_m(X,\Q_\ell)$.

\end{proposition}

On note aussi
$\DP^{\geq (a_0,\dots,a_n)}=\DP^{>(a_0-1,\dots,a_n-1)}$.
Il est évident d'après la définition de $\DP^{\leq\as}$ et $\DP^{>\as}$ que ce sont des sous-catégories stables par décalage et extensions de $D^b_m(X,\Q_\ell)$,
que la dualité de Poincaré échange $\DP^{\leq\as}$ et $\DP^{\geq -\as}$,
et que 
$$\DP^{\leq (a_0,\dots,a_n)}(1)=\DP^{\leq (a_0-2,\dots,a_n-2)}$$
$$\DP^{>(a_0,\dots,a_n)}(1)=\DP^{>(a_0-2,\dots,a_n-2)}.$$

D'après [BBD] 1.3.3, l'inclusion $\DP^{\leq\as}\subset D^b_m(X,\Q_\ell)$ (resp. $\DP^{>\as}\subset D^b_m(X,\Q_\ell)$) 
admet un adjoint à droite (resp. à gauche), qu'on note
$w_{\leq\as}$ (resp. $w_{>\as}$).
Pour tout $K\in D^b_m(X,\Q_\ell)$, il existe un unique morphisme $w_{>\as}K\fl (w_{\leq\as}K)[1]$ tel que le triangle
$$w_{\leq\as}K\fl K\fl w_{>\as}K\fl (w_{\leq\as}K)[1]$$
soit distingué.\newline
De plus, à isomorphisme unique près, il existe un unique triangle distingué $K'\fl K\fl K''\stackrel{+1}{\fl}$
avec $K'\in\DP^{\leq\as}$ et $K''\in\DP^{>\as}$.

Enfin, la dualité de Poincaré échange $w_{\leq\as}$ et $w_{\geq -\as}$ (car elle échange $\DP^{\leq\as}$ et $\DP^{\geq -\as}$).

\begin{lemme} Si $a_0=\dots=a_n=a$, 
alors $(\DP^{\leq\as},\DP^{>\as})$ est la t-structure $(\DP^{\leq a}(X),\DP^{>a}(X))$ de la section 4.1.

\end{lemme}

\begin{preuve} Soit $K\in\DP^{\leq a}$.
D'après le (iv) de la proposition \ref{t_structures:prop2},
pour tout $k\in\{0,\dots,n\}$, $i_k^*K\in\DP^{\leq a}(S_k)$.
Donc $K\in\DP^{\leq\as}$.
Par dualité, on a $\DP^{>a}\subset\DP^{>\as}$.

Soit $K\in\DP^{\leq\as}$.
On a un triangle distingué
$$w_{\leq a}K\fl K\fl w_{>a}K\stackrel{+1}{\fl}.$$
D'après ce qui précède, 
$w_{\leq a}K\in\DP^{\leq\as}$ et $w_{>a}\in\DP^{>\as}$,
donc $w_{\leq a}K=w_{\leq\as}K=K$, et $K\in\DP^{\leq a}$.
Par dualité, on a $\DP^{>\as}\subset\DP^{>a}$.

\end{preuve}

\begin{proposition}\label{prop1:t_structures_recollees} Pour tous $a\in\Z\cup\{\pm\infty\}$ et $k\in\{0,\dots,n\}$, on note
$$w_{\leq a}^k=w_{\leq (+\infty,\dots,+\infty,a,+\infty,\dots,+\infty)}$$
$$w_{\geq a}^k=w_{\geq (-\infty,\dots,-\infty,a,-\infty,\dots,-\infty)}$$
où, dans les deux formules, le $a$ est en $k$-ième position (et on commence à compter à $0$).

\begin{itemize}
\item[(i)] On a
$$w_{\leq\as}=w^n_{\leq a_n}\circ\dots\circ w^0_{\leq a_0}$$
$$w_{\geq\as}=w^n_{\geq a_n}\circ\dots\circ w^0_{\geq a_0}.$$
% Comme $T_k^{\leq+\infty}=T_k^{\geq-\infty}=id$ pour tout $k$, on peut retirer de la première formule
% les $T_k^{\leq a_k}$ tels que $a_k=+\infty$, 
% et de la deuxième formule les $T_k^{\geq a_k}$ tels que $a_k=-\infty$.
\item[(ii)] Soient $a\in\Z\cup\{\pm\infty\}$, $k\in\{0,\dots,n\}$ et $K\in D^b_m(X,\Q_\ell)$.
Alors on a des triangles distingués (uniques à isomorphisme unique près)
$$w_{\leq a}^kK\fl K\fl Ri_{k*}w_{>a}i_k^*K\stackrel{+1}{\fl}$$
$$i_{k!}w_{\leq a}i_k^!K\fl K\fl w^k_{>a}K\stackrel{+1}{\fl}.$$
\end{itemize}
\end{proposition}

\begin{preuve}
\begin{itemize}
\item[(i)] Il suffit d'appliquer plusieurs fois [BBD] 1.4.13.1.
\item[(ii)] Il suffit de traiter le cas de $w_{\leq a}^k$ (celui de $w_{>a}^k$
en résulte par dualité).
On définit $\as\in (\Z\cup\{\pm\infty\})^{n+1}$ par :
$a_r=+\infty$ et si $r\not=k$, et $a_k=a$.

On note
$L=Ri_{k*}w_{>a}i_k^*K\in\DP^{>\as}$,
et on complète le morphisme évident $K\stackrel{adj}{\fl}Ri_{k*}i_k^*K \fl L$ en un triangle distingué
$L'\fl K\fl L\stackrel{+1}{\fl}$.

Il suffit de montrer que $L\in \DP^{>\as}$ et que $L'\in\DP^{\leq\as}$.
Si $r\not=k$, $i_r^!L=0$, car $i_r^!Ri_{k*}=0$,
donc $i_r^!L\in\DP^{>a_r}(S_r)$ (et il est évident que $i_r^*L'\in\DP^{\leq a_r}(S_r)$).
De plus, $i_k^*L=i_k^!L=w_{>a}i_k^*K\in\DP^{>a_k}(S_k)$ et 
on a un triangle distingué
$i_k^*L'\fl i_k^*K\fl i_k^*L=w_{>a}i_k^*K\stackrel{+1}{\fl}$,
donc $i_k^*L'\simeq w_{\leq a}i_k^*K\in\DP^{\leq a_k}(S_k)$.

\end{itemize}
\end{preuve}

\begin{theoreme}\label{th:simplification_formule_traces} 
Pour tout $\as\in (\Z\cup\{\pm\infty\})^{n+1}$,
pour tout $K\in\DP^{\leq a_0}(U)$ 
on a
$$[w_{\leq\as}Rj_*K]=\sum_{1\leq n_1<\dots <n_r\leq n}(-1)^{r}[Ri_{n_r*}w_{>a_{n_r}}i_{n_r}^*\dots Ri_{n_1*}w_{>a_{n_1}}i_{n_1}^*Rj_*K]$$
dans le groupe de Grothendieck de $D^b_m(X,\Q_\ell)$.

\end{theoreme}

% Dans la suite, on utilisera les notations suivantes :

% \begin{notation} Soient $k\in\{1,\dots,n\}$ et $a\in\Z\cup\{\pm\infty\}$.
% On pose
% $$\sigma^k_{>a}=\left\{\begin{array}{ll}Ri_{k*}w_{>a}i_k^* & \mbox{ si }a>-\infty \\
% id & \mbox{ si }a=-\infty\end{array}\right.$$
% En particulier, $\sigma^k_{>+\infty}=0$, et, si $a>-\infty$, 
% $$\sigma^k_{>a}=w_{>(+\infty,\dots,+\infty,a,+\infty,\dots,+\infty)}$$
% (le $a$ est en $k$-ième position si on commence à compter à $0$),
% donc on a pour tout $K\in D^b_m(X,\Q_\ell)$ un triangle distingué canonique
% $$w_{\leq a}^kK\fl K\fl \sigma^k_{>a}K\stackrel{+1}{\fl}.$$
% De plus, grâce au (iii) de la proposition \ref{t_structures:prop2}, 
% $\sigma^k_{>a}$ est un foncteur triangulé.

% Si $\as\in(\Z\cup\{\pm\infty\})^{n+1}$, pour tout $I\subset\{1,\dots,n\}$, on définit $\as^I$ par
% $$a_k^I=\left\{\begin{array}{ll}a_k & \mbox{ si }k\in I \\
% -\infty & \mbox{ si }k\not\in I\end{array}\right.$$

% \end{notation}

\begin{preuve} 

D'après la proposition ci-dessus, on a, dans l'anneau des endomorphismes du groupe 
de Grothendieck de $D^b_m(X,\Q_\ell)$ :
$$w_{\leq\as}=w_{\leq a_n}^n\circ\dots\circ w_{\leq a_0}^0=(1-Ri_{n*}w_{>a_n}i_n^*)\circ\dots\circ (1-Ri_{1*}w_{>a_1}i_1^*)\circ(1-Rj_*w_{>a_0}j^*).$$
Le théorème résulte de cette égalité et du fait que $Rj_*w_{>a_0}j^*Rj_*K=0$ 
(car $K\in\DP^{\leq a_0}(U)$).

\end{preuve}

\section{Propriétés supplémentaires des t-structures recollées}

\begin{proposition}\label{prop2:t_structures_recollees}
\begin{itemize}
\item[(i)] Si $K\in D^{\leq\as}$ et $L\in\DP^{>\as}$, alors $R\Hom(K,L)=0$.
\item[(ii)] $w_{\leq\as}$ et $w_{>\as}$ commutent au foncteur de décalage $[1]$,
et ils envoient les triangles distingués sur des triangles distingués.\newline
Pour tout $K\in D^b_m(X,\Q_\ell)$, on a
$$w_{\leq (a_0,\dots,a_n)}(K(1))=(w_{\leq (a_0+2,\dots,a_n+2)}K)(1)$$
$$w_{>(a_0,\dots,a_n)}(K(1))=(w_{>(a_0+2,\dots,a_n+2)}K)(1).$$
\item[(iii)] Soit $\as'=(a'_0,\dots,a'_n)\in (\Z\cup\{\pm\infty\})^{n+1}$ tel que $a_k\leq a'_k$
pour tout $k\in\{1,\dots,n\}$.
Alors, pour tous $K\in\DP^{\leq\as}$ et $L\in\DP^{>\as'}$, le morphisme canonique
$$R\Hom(K,L)\fl R\Hom(j^*K,j^*L)$$
est un isomorphisme.
\item[(iv)] Soient $f:Y\fl X$ un morphisme et $(S'_k)_{0\leq k\leq n}$ une stratification de $Y$
telle que, pour tout $k\in\{0,\dots,n\}$, $f(S'_k)\subset S_k$.
On suppose que la dimension des fibres de $f$ est inférieure ou égale à $d$.
Alors
$$f^*(\DP^{\leq (a_0,\dots,a_n)}(X))\subset\DP^{\leq (a_0+d,\dots,a_n+d)}(Y)$$
$$f^!(\DP^{>(a_0,\dots,a_n)}(X))\subset\DP^{>(a_0-d,\dots,a_n-d)}(Y)$$
$$Rf_*(\DP^{>(a_0,\dots,a_n)}(Y))\subset\DP^{>(a_0-d,\dots,a_n-d)}(X)$$
$$Rf_!(\DP^{\leq (a_0,\dots,a_n)}(Y))\subset\DP^{\leq (a_0+d,\dots,a_n+d)}(X).$$
\end{itemize}
\end{proposition}

\begin{preuve}
\begin{itemize}
\item[(i)] Montrons le résultat par récurrence sur $n$.
Si $n=0$, c'est le lemme \ref{lemme_annulation}.
Soit $n\geq 1$, et supposons le résultat vrai pour $n'<n$.
On note $\as'=(a_0,\dots,a_{n-1})$, $V=\displaystyle{\bigcup_{k=0}^{n-1}S_k}$, $Y=S_n=X-V$, 
$j_1$ l'immersion ouverte de $V$ dans $X$ et 
$i$ l'immersion fermée de $S_n$ dans $X$.
Soient $K\in\DP^{\leq\as}$ et $L\in\DP^{>\as}$.
On a un triangle distingué
$$R\Hom(i^*K,i^!L)\fl R\Hom(K,L)\fl R\Hom(j_1^*K,j_1^*L)\stackrel{+1}{\fl}.$$
Or $j_1^*K\in\DP^{\leq\as'}(U)$, $j_1^*L\in\DP^{>\as'}(U)$, $i^*K\in\DP^{\leq a_n}(Y)$ et $i_k^!L\in\DP^{>a_n}(Y)$,
donc, d'après l'hypothèse de récurrence,
$$R\Hom(j_1^*K,j_1^*L)=R\Hom(i^*K,i^!L)=0,$$
d'où 
$$R\Hom(K,L)=0.$$
\item[(ii)] $w_{\leq\as}$ et $w_{>\as}$ commutent au foncteur de décalage car $\DP^{\leq\as}$ et $\DP^{>\as}$ sont stables par décalage.
Soit $K\fl K'\fl K''\stackrel{+1}{\fl}$ un triangle distingué.
D'après [BBD] 1.1.11, on peut construire un diagramme commutatif dont les lignes et les colonnes sont distinguées
$$\xymatrix{w_{\leq\as}K\ar[r]\ar[d] & w_{\leq\as}K'\ar[r]\ar[d] & L\ar[r]^{+1}\ar[d] & \\
K\ar[r]\ar[d] & K'\ar[r]\ar[d] & K''\ar[r]^{+1}\ar[d] & \\
w_{>\as}K\ar[r]\ar[d]^{+1} & w_{>\as}K'\ar[r]\ar[d]^{+1} & L'\ar[r]^{+1}\ar[d]^{+1} & \\
 & & & }$$
Comme $\DP^{\leq\as}$ et $\DP^{>\as}$ sont des sous-catégories stables par extensions de $D^b_m(X,\Q_\ell)$,
$L\in\DP^{\leq\as}$ et $L'\in\DP^{>\as}$, donc 
$L=w_{\leq\as}K''$ et $L'=w_{>\as}K''$.\newline
La dernière assertion se prouve exactement comme la propriété analogue dans le (i) de la proposition \ref{t_structures:prop2}.
\item[(iii)] Notons $i$ l'immersion fermée $X-U=S_1\cup\dots\cup S_n\subset X$, $\underline{b}=(a_1,\dots,a_n)$ et $\underline{b}'=(a'_1,\dots,a'_n)$.
Soient $K\in\DP^{\leq\as}$ et $L\in\DP^{>\as'}$.
On a un triangle distingué canonique
$$R\Hom(i^*K,i^!L)\fl R\Hom(K,L)\fl R\Hom(j^*K,j^*L)\stackrel{+1}{\fl}.$$
Or $i^*K\in\DP^{\leq\underline{b}}(X-U)$ et $i^!L\in\DP^{>\underline{b}'}(X-U)$,
donc, d'après le point (i), $R\Hom(i^*K,i^!L)=0$.
\item[(iv)] Il suffit de traiter les cas de $f^*$ et $Rf_!$, car
ceux de $f^!$ et $Rf_*$ en résultent par dualité.\newline
Pour tout $k\in\{0,\dots,n\}$, on note $i'_k$ l'inclusion $S'_k\subset Y$ et 
$f_k:S'_k\fl S_k$ la restriction de $f$.
Soit $K\in\DP^{\leq\as}(X)$.
Pour tout $k\in\{0,\dots,n\}$, ${i'_k}^*f^*K=f_k^*i_k^*K$; 
$i_k^*K\in\DP^{\leq a_k}(S_k)$ par la définition de $\DP^{\leq\as}(X)$,
donc, d'après le (iv) de la proposition \ref{t_structures:prop2}, $f_k^*i_k^*K\in\DP^{\leq a+d}(S'_k)$.
On a donc bien $f^*K\in\DP^{\leq (a_0+d,\dots,a_n+d)}(Y)$.

Soit $K\in\DP^{\leq\as}(Y)$.
Fixons $k\in\{0,\dots,n\}$. 
Comme $S'_k=f^{-1}(S_k)$,
le diagramme suivant est cartésien aux nilpotents près
$$\xymatrix{S'_k\ar[r]^{i'_k}\ar[d]_{f_k} & Y\ar[d]_f \\
S_k\ar[r]^{i_k} & X}$$
donc, d'après le théorème de changement de base propre, ${i_k}^*Rf_!K\simeq Rf_{k!}{i'_k}^*K$.
Or ${i'_k}^*K\in\DP^{\leq a_k}(S'_k)$, donc, d'après le (iv) de la proposition \ref{t_structures:prop2},
$i_k^*Rf_!K\simeq Rf_{k!}{i'_k}^*K\in\DP^{\leq a_k+d}(S_k)$.
On a donc bien $Rf_!K\in\DP^{\leq (a_0+d,\dots,a_n+d)}(X)$.
\end{itemize}
\end{preuve}

% \begin{lemme} Soient $a\in\Z\cup\{\pm\infty\}$, $k\in\{0,\dots,n\}$ et $K\in D^b_m(X,\Q_\ell)$.
% Alors on a des isomorphismes canoniques
% $$i_k^*w^k_{\leq a}K=w_{\leq a}i_k^*K$$
% $$i_k^!w^k_{>a}K=w_{>a}i_k^!K$$
% $$i_k^!w^k_{\leq a}K=w_{>a}i_k^*K[-1]\mbox{ si }i_k^!K=0$$
% $$i_k^*w^k_{>a}K=w_{\leq a}i_k^!K[1]\mbox{ si }i_k^*K=0.$$
% Si $0\leq r<k$, on a des isomorphismes canoniques
% $$i_r^*w^k_{\leq a}K=i_r^*w^k_{>a}K=i_r^*K$$
% $$i_r^!w^k_{\leq a}K=i_r^!w^k_{>a}K=i_r^!K.$$
% Si $k<r\leq n$, on a des isomorphismes canoniques
% $$i_r^!w^k_{\leq a}K=i_r^!K$$
% $$i_r^*w^k_{>a}K=i_r^*K.$$

% \end{lemme}

% \begin{preuve}
% En appliquant $i_k^*$ et $i_k^!$ aux triangles distingués du (ii) de la proposition \ref{prop1:t_structures_recollees}, 
% on trouve des triangles distingués
% $$i_k^*w^k_{\leq a}K\fl K\fl w_{>a}i_k^*K\stackrel{+1}{\fl}$$
% $$i_k^!w_{\leq a}i_k^!K\fl K\fl i_k^!w^k_{>a}K\stackrel{+1}{\fl}$$
% $$i_k^!w^k_{\leq a}K\fl i_k^!K\fl w_{>a}i_k^*K\stackrel{+1}{\fl}$$
% $$w_{\leq a}i_k^!K\fl K\fl i_k^*w^k_{>a}K\stackrel{+1}{\fl},$$
% d'où les quatre premiers isomorphismes. \newline
% Soit $r\in\{0,\dots,k-1\}$.
% Alors $i_r^*Ri_{k*}=i_r^!Ri_{k*}=i_r^*i_{k!}=i_r^!i_{k!}=0$, donc on obtient les isomorphismes cherchés
% en appliquant $i_r^*$ et $i_r^!$ aux deux triangles distingués du (ii) de la proposition \ref{prop1:t_structures_recollees}.\newline
% Si $r\in\{k+1,\dots,n\}$, on procède de même, en remarquant que $i_r^!Ri_{k*}=i_r^*i_{k!}=0$.

% \end{preuve}

La proposition suivante est une reformulation de [BBD] 1.4.14 dans le cas particulier considéré.

\begin{proposition}\label{prop3:t_structures_recollees}
Soit $\as=(a_0,\dots,a_n)\in (\Z\cup\{\pm\infty\})^{n+1}$.
On note $\protect{\as'\!=\!(a_0,a_1+1,\dots,a_n+1)}$.
Alors, pour tout $K\in\DP^{\leq a_0}(U)\cap\DP^{\geq a_0}(U)$, 
$w_{\geq\as'}j_!K=w_{\leq\as}Rj_*K$ est l'unique prolongement de $K$ dans $\DP^{\leq\as}\cap\DP^{\geq\as'}$.

En particulier, si $K$ est pervers pur de poids $a$ sur $U$, on a
$$w_{\geq (a,a+1,\dots,a+1)}j_!K=j_{!*}K=w_{\leq a}Rj_*K,$$
et par dualité on obtient aussi
$$w_{\geq a}j_!K=j_{!*}K=w_{\leq (a,a-1,\dots,a-1)}Rj_*K.$$
(On retrouve le résultat du théorème \ref{th:calcul_IC_general}.)

\end{proposition}

% \begin{preuve}
% Soit $K\in\DP^{\leq a_0}(U)\cap\DP^{\geq a_0}(U)$.
% D'après le (iii) de la proposition \ref{prop2:t_structures_recollees}, 
% $K$ a au plus un prolongement dans $\DP^{\leq\as}\cap\DP^{\geq\as'}$.
% Notons $L=w_{\leq\as}Rj_*K$ et $L'=w_{\geq\as'}j_!K$.
% Comme $K\in\DP^{\leq a_0}(U)\cap\DP^{\geq a_0}(U)$, d'après le (ii) de la proposition \ref{prop1:t_structures_recollees}, $j^*L=j^*L'=K$.
% Il suffit donc de montrer que $L\in\DP^{\geq\as'}$ et $L'\in\DP^{\leq\as}$.
% Soit $k\in\{1,\dots,n\}$. 
% D'après le (i) de la proposition \ref{prop1:t_structures_recollees} et le lemme ci-dessus
% (et le fait que $K\in\DP^{\leq a_0}(U)\cap\DP^{\geq a_0}(U)$),
% $$i_k^!L=i_k^!w^n_{\leq a_n}\dots w^1_{\leq a_1}Rj_*K=i_k^!w^k_{\leq a_k}\dots w^1_{\leq a_1}Rj_*K$$
% $$i_k^*L'=i_k^*w^n_{\geq a'_n}\dots w^1_{\geq a'_1}j_!K=i_k^*w^k_{\geq a'_k}\dots w^1_{\geq a'_1}j_!K.$$
% De plus, toujours d'après le lemme ci-dessus, 
% $$i_k^!w^{k-1}_{\leq a_{k-1}}\dots w^1_{\leq a_1}Rj_*K=i_k^!Rj_*K=0$$
% $$i_k^*w^{k-1}_{\geq a'_{k-1}}\dots w^1_{\geq a'_1}j_!K=i_k^*j_!K=0,$$
% donc 
% $$i_k^!L=i_k^!w^k_{\leq a_k}\dots w^1_{\leq a_1}Rj_*K=w_{\geq a_k+1}i_k^*Rj_*K[-1]\in\DP^{\geq a'_k}(S_k)$$
% $$i_k^*L=i_k^*w^k_{\geq a'_k}\dots w^1_{\geq a'_1}j_!K=w_{\leq a_k}i_k^!j_!K[1]\in\DP^{\leq a_k}(S_k).$$

% \end{preuve}

\section{Quelques lemmes techniques}

Ce paragraphe contient quelques lemmes utilisés dans les preuves des résultats des parties 4 et 5.

\begin{lemme}\label{fp_lemme1} Soient $X$ un schéma séparé de type fini sur $\Fi_q$ et $(S_\alpha)$ une partition finie de $X$ par des sous-schémas localement fermés. 
Pour tout $\alpha$, on note $i_\alpha:S_\alpha\fl X$ l'inclusion. Alors, pour tout $K\in D^b_c(X,\Q_\ell)$, on a
\begin{itemize}
\item[(a)] $K\in {}^{p}D_c^{\leq 0}$ \ssi pour tout $\alpha$, pour tout $i\geq 1$, on a $\Hp^i(i_\alpha^*K)=0$;
\item[(b)] $K\in {}^{p}D_c^{\geq 0}$ \ssi pour tout $\alpha$, pour tout $i\leq -1$, on a $\Hp^i(i_\alpha^!K)=0$.

\end{itemize}
De plus, si $U$ est un ouvert de $X$ réunion de strates, si $j:U\fl X$ est l'inclusion et si $K$ est un faisceau pervers sur $U$, alors $j_{!*}K$ est l'unique prolongement $L\in D^b_c(X,\Q_\ell)$ de $K$ tel que : pour tout $\alpha$ tel que $S_\alpha\subset X-U$, on a
$$\left\{\begin{array}{c}{}^{p}H^i(i_\alpha^*L)=0\mbox{ pour }i\geq 0 \\
{}^{p}H^i(i_\alpha^!L)=0\mbox{ pour }i\leq 0\end{array}\right..$$

\end{lemme}

\begin{preuve} Il suffit de montrer (a), car (b) en résulte par dualité. 
D'après [BBD] 2.2.5, les $i_\alpha^*$ sont $t$-exacts à droite, donc, si $K\in {}^{p}D_c^{\leq 0}$, on a bien $\Hp^i(i_\alpha^*K)=0$ pour tout $\alpha$ et pour $i\geq 1$.
 
Réciproquement, soit $K\in D^b_c(X,\Q_\ell)$ tel que pour tout $\alpha$, pour tout $i\geq 1$, $\Hp^i(i_\alpha^*K)=0$.
On sait par [BBD] 2.2.12 qu'un complexe $L$ est dans ${}^{p}D_c^{\leq 0}$ \ssi pour tout point $x$ de $X$, notant $i_x:x\fl X$ et $dim(x)=dim(\overline{\{x\}})$, on a $H^i(i_x^*L)=0$ pour $i<p(2 dim(x))$. On va utiliser cette caractérisation pour montrer que $K\in {}^{p}D_c^{\leq 0}$.
Soit $x$ un point de $X$, et soit $\alpha$ tel que $x\in S_\alpha$. Comme $S_\alpha$ est localement fermé dans $X$, $\overline{\{x\}}\cap S_\alpha$ est ouvert dans $\overline{\{x\}}$, donc $dim(\overline{\{x\}}\cap S_\alpha)=dim(\overline{\{x\}})$ ($\overline{\{x\}}$ est irréductible), et $dim(x)$ ne change pas si on considère $x$ comme un point de $S_\alpha$. 
Comme par hypothèse $i_\alpha^*K\in {}^{p}D_c^{\leq 0}(S_\alpha)$, on a bien $H^i(i_x^*K)=0$ si $i<p(2 dim(x))$. 

Montrons enfin la dernière assertion du lemme. $U$ est un ouvert de $X$ réunion de strates, $j:U\fl X$ est l'inclusion, $K$ est un faisceau pervers sur $U$. D'après [BBD] 1.4.24 (qui s'applique par [BBD] 2.2.3 et 2.2.11), on a $\Hp^0(i_\alpha^*j_{!*}K)=0$ pour tout $\alpha$ tel que $S_\alpha\subset X-U$. L'annulation des $\Hp^0(i_\alpha^!j_{!*}K)$ s'en déduit par dualité.

Soit $L\in D^b_c(X,\Q_\ell)$, muni d'un isomorphisme $j^*L\simeq K$, tel que, pour tout $\alpha$ tel que $S_\alpha\subset X-U$, on ait $\Hp^i(i_\alpha^*L)=0$ pour $i\geq 0$ et $\Hp^i(i_\alpha^!L)=0$ pour $i\leq 0$. D'après ce qui précède, on sait que $L$ est pervers.
En raisonnant par récurrence sur le cardinal de $\{\alpha\mbox{ tq }S_\alpha\subset X-U\}$, on se ramène au cas où $X-U=S_\alpha$ est une strate. Notons $i=i_\alpha$. 
On a un triangle distingué
$$i_*i^!L\fl L\fl Rj_*j^*L\simeq Rj_* K\stackrel{+1}{\fl},$$
d'où une suite exacte
$$\Hp^0(i_*i^! L)\fl \Hp^0(L)=L\fl \Hp^0(Rj_*K).$$
Comme $i_*$ est t-exact ([BBD] 2.2.6), $\Hp^0(i_*i^!L)=i_*\Hp^0(i^!L)=0$, et le morphisme $L\fl\Hp^0(Rj_*K)$ est injectif.
D'autre part, on a un triangle distingué
$$j_!j^*L\simeq j_!K\fl L\fl i_*i^*L\stackrel{+1}{\fl},$$
d'où une suite exacte
$$\Hp^0(j_!K)\fl L\fl \Hp^0(i_*i^* L)=i_*\Hp^0(i^* L)=0.$$
Le morphisme $\Hp^0(j_!K)\fl L$ est donc surjectif, ce qui finit la démonstration.

\end{preuve}

\begin{lemme}\label{fp_lemme2} Soient $X$ un schéma de type fini lisse purement de dimension $d$ sur un corps $k$ de caractéristique $0$ ou fini et $K\in D^b_c(X,\Q_\ell)$. On suppose que les $H^i(K)$ sont lisses. Alors, pour tout $i\in \Z$, $\Hp^i(K)=H^{i-d}(K)[d]$.

\end{lemme}

\begin{preuve} On montre le résultat par récurrence sur le cardinal $N(K)$ de\newline
 $\{i\in\Z\mbox{ tq }H^i(K)\not=0\}$.

Si $N(K)=1$, on a $K\simeq H^i(K)[-i]$ pour un $i\in\Z$, donc $K[i+d]$ est pervers, et
$$\Hp^j(K)=\left\{\begin{array}{ll}0 & \mbox{ si }j\not=i+d \\
H^i(K)[d] & \mbox{ si }j=i+d\end{array}\right..$$

Soit $K$ tel que $N(K)>1$, et supposons le résultat prouvé pour tous les $L$ tels que\newline
 $N(L)<N(K)$. Soit $i=max\{k\in\Z\mbox{ tq }H^k(K)\not=0\}$. Comme $H^i(K)[-i]$ est lisse, on a comme plus haut
$$\Hp^j(H^i(K)[-i])=\left\{\begin{array}{ll}0 & \mbox{ si }j\not=i+d \\
H^i(K)[d] & \mbox{ si }j=i+d\end{array}\right..$$
D'autre part, d'après l'hypothèse de récurrence, on a 
$$\Hp^j(\tau_{\leq i-1}K)=H^{j-d}(\tau_{\leq i-1}K)[d]=\left\{\begin{array}{ll}0 & \mbox{ si }j\geq i+d \\
H^{j-d}(K)[d] & \mbox{ si }j<i+d\end{array}\right..$$
On conclut en utilisant la suite exacte longue de cohomologie perverse du triangle distingué 
$$\tau_{\leq i-1}K\fl\tau_{\leq i}K\simeq K\fl H^i(K)[-i]\stackrel{+1}{\fl}.$$

\end{preuve}

\newpage

\chapter{Complexes pondérés sur les compactifications de Baily-Borel}

\section{Complexes pondérés}

Dans cette section, on définit à l'aide des foncteurs $w_{\leq\as}$ de 4.3
une famille de complexes sur la compactification de Baily-Borel,
qu'on appellera complexes pondérés,
et dont les complexes d'intersection sont des cas particuliers.
Ces complexes pondérés sont des analogues en caractéristique finie 
des complexes pondérés définis sur $M^{\K}(\G,\X)^*(\C)$ par
Goresky, Harder et MacPherson dans [GHM].

\begin{notation} Soient $\G$ un groupe algébrique sur $\Q$, $\lambda:\Gm\fl\G$ un cocaractère central, $\ell$ un nombre premier et $t\in\Z\cup\{\pm\infty\}$.
Pour tout $V\in Rep_{\G(\Q_\ell)}$, on note $w_{< t}V$ (resp. $w_{\geq t}V$) la plus grande sous-représentation de $V$
sur laquelle le poids relativement à $\lambda$ est $< t$ (resp. $\geq t$).

Les foncteurs exacts $w_{< t}$ et $w_{\geq t}$ s'étendent trivialement à la catégorie dérivée $D^b(Rep_{\G(\Q_\ell)})$.
Pour tout $V\in D^b(Rep_{\G(\Q_\ell)})$, on a
$$V=w_{< t}V\oplus w_{\geq t}V,$$
et, pour tout $i\in\Z$, $H^i(w_{< t}V)=w_{< t}H^i(V)$ est de poids $< t$ et $H^i(w_{\geq t}V)=w_{\geq t}H^i(V)$ de poids $\geq t$.

\end{notation}

\begin{lemme}\label{relation_tronques} Comme dans le paragraphe 2.1.4, on se place dans la situation suivante :
$\L$ est un groupe algébrique connexe
qui est produit quasi-direct de deux sous-groupes $\G_\lin$ et $\G$,
$(\G,\X)$ est une donnée de Shimura pure
telle que $\G_\lin(\Q)\cap\G(\Q)$ agisse trivialement sur $\X$,
$\K_L$ est un sous-groupe ouvert compact net de $\L(\Af)$,
et on note $\Hr=\K_L\cap\L(\Q)\G(\Af)$, $\Gamma_\ell=\K_L\cap Cent_{\L(\Q)}(\X)$.
On se place sur les réductions modulo $p$ des variétés, 
où le nombre premier $p$ est comme dans le corollaire \ref{poids_faisceaux_generalises}.
On note $d$ la dimension de $M^{\Hr/\Gamma_\lin}(\G,\X)$.
Alors, pour tous $V\in D^b(Rep_{\L(\Q_\ell)})$ et $a\in\Z\cup\{\pm\infty\}$,
$$w_{\leq a}\F^{\Hr/\Gamma_\lin}R\Gamma(\Gamma_\lin,V)=\F^{\Hr/\Gamma_\lin}R\Gamma(\Gamma_\lin,w_{\geq d-a}V)$$
$$w_{>a}\F^{\Hr/\Gamma_\lin}R\Gamma(\Gamma_\lin,V)=\F^{\Hr/\Gamma_\lin}R\Gamma(\Gamma_\lin,w_{<d-a}V).$$

\end{lemme}

\begin{preuve} Notons $t=d-a$, $K=\F^{\Hr/\Gamma_\lin}R\Gamma(\Gamma_\lin,V)$, $K_1=\F^{\Hr/\Gamma_\lin}R\Gamma(\Gamma_\lin,w_{\geq t}V)$ et $K_2=\F^{\Hr/\Gamma_\lin}R\Gamma(\Gamma_\lin,w_{<t}V)$.
Comme $V=w_{<t}V\oplus w_{\geq t}V$, on a $K=K_1\oplus K_2$.
Il suffit de montrer que $K_1$ est dans $\DP^{\leq a}(M^{\Hr/\Gamma_\lin}(\G,\X))$ et $K_2$ dans $\DP^{>a}(M^{\Hr/\Gamma_\lin}(\G,\X))$.
$M^{\Hr/\Gamma_\lin}(\G,\X)$ est lisse et pour tout $i\in\Z$, $H^i(K_1)=\F^{\Hr/\Gamma_\lin}H^i(\Gamma_\lin,w_{\geq t}V)$ et $H^i(K_2)=\F^{\Hr/\Gamma_\lin}H^i(\Gamma_\lin,w_{<t}V)$ sont lisses,
donc, d'après le lemme \ref{fp_lemme2}, pour tout $i\in\Z$, $\Hp^i(K_1)=H^{i-d}(K_1)[d]$ et $\Hp^i(K_2)=H^{i-d}(K_2)[d]$.
Soit $i\in\Z$.
D'après le corollaire \ref{poids_faisceaux_generalises}, $H^{i-d}(K_1)$ est de poids $\leq -t$ et $H^{i-d}(K_2)$ est de poids $>-t$.
On en déduit que $\Hp^i(K_1)=H^{i-d}(K_1)[d]$ est de poids $\leq -t+d=a$ et que $\Hp^i(K_2)=H^{i-d}(K_2)[d]$ est de poids $>-t+d=a$.

\end{preuve}

On se place maintenant dans la situation de la définition \ref{bon_modele_entier}, 
c'est-à-dire que $(\G,\X)$ est la donnée de Shimura de la section 1.2
et qu'on a inversé assez de nombres premiers pour avoir un modèle entier de la compactification
de Baily-Borel possédant toutes les propriétés souhaitables.
On travaille sur les réductions modulo $p$ des variétés.
On pose $M^*=M^{\K}(\G,\X)^*$, $M_0=M^{\K}(\G,\X)$ et, pour tout $r\in\{1,\dots,q\}$, on note $M_r$ l'union des strates de bord correspondant à $(\QP_r,\Y_r)$.
Pour tout $r\in\{0,\dots,q\}$, $dim(M_r)=(p-r)(q-r)$.
$(M_0,M_1,\dots,M_q)$ est une stratification de $M^*$ (au sens de la définition \ref{definition_stratification})
, et c'est toujours celle qu'on utilisera dans la suite.

\begin{definition}\label{def:complexes_ponderes} Soient $t_1,\dots,t_q\in\Z\cup\{\pm\infty\}$.
Pour tout $r\in\{1,\dots,q\}$, on pose $a_r=-t_r+(p-r)(q-r)$.
On définit un foncteur additif triangulé
$$W^{\geq t_1,\dots,\geq t_q}:D^b(Rep_{\G(\Q_\ell)})\fl D^b_m(M^{\K}(\G,\X)^*)$$
de la manière suivante : pour tout $m\in\Z$, si $V\in D^b(Rep_{\G(\Q_\ell)})$ est tel que $H^i(V)$ soit de poids $m$ pour tout $i\in\Z$, alors
$$W^{\geq t_1,\dots,\geq t_q}(V)=w_{\leq (-m+pq,-m+a_1,\dots,-m+a_q)}Rj_*\F^{\K}V.$$

\end{definition}

\begin{definition}\label{def:IC} Soit $V\in Rep_{\G(\Q_\ell)}$.
Comme $M^{\K}(\G,\X)$ est de dimension $pq$, $\F^{\K}V[pq]$ est un faisceau pervers sur $M^{\K}(\G,\X)$.
On pose
$$IC^{\K}V=(j_{!*}(\F^{\K}V[pq]))[-pq].$$

\end{definition}

\begin{proposition}\label{calcul_IC} 
\begin{itemize}
\item[(1)] Pour tous $t_1,\dots,t_q\in\Z\cup\{\pm\infty\}$ et pour tout $V\in D^b(Rep_{\G(\Q_\ell)})$, on a un isomorphisme canonique
$$D(W^{\geq t_1,\dots,\geq t_q}(V))\simeq W^{\geq s_1,\dots,\geq s_q}(V^*)[2pq](pq),$$
où $D$ est le foncteur dualisant,
$V^*=R\Hom(V,\Q_\ell)$ est la représentation duale de $V$ et $s_r=1-t_r+2r(r-n)$ pour tout $r\in\{1,\dots,q\}$.
\item[(2)] Notons, pour tout $r\in\{1,\dots,q\}$, $t_r=r(r-n)+1$ et $s_r=r(r-n)$.
Alors, pour tout $V\in Rep_{\G(\Q_\ell)}$, on a des isomorphismes canoniques
$$IC^{\K}V\simeq W^{\geq t_1,\dots,\geq t_q}(V)\iso W^{\geq s_1,\dots,\geq s_q}(V).$$

\end{itemize}
\end{proposition}

\begin{preuve}
\begin{itemize}
\item[(1)] On peut supposer que $V\in Rep_{\G(\Q_\ell)}$ et que $V$ est pure. 
Notons $m$ le poids de $V$.
$V^*$ est alors une représentation de $\G(\Q_\ell)$ de poids $-m$.
On pose $a_0=-m+pq$ et, pour tout $r\in\{1,\dots,q\}$, $a_r=-(t_r+m)+(p-r)(q-r)$.
Alors
$$W^{\geq t_1,\dots,\geq t_q}(V)=w_{\leq (a_0,\dots,a_q)}Rj_*\F^{\K}V,$$
donc
$$\begin{array}{rcl}D(W^{\geq t_1,\dots,\geq t_q}(V)) & = & w_{\geq (-a_0,\dots,-a_q)}(j_!\F^{\K}V^*[2pq](pq)) \\
& = & (w_{\geq (-a_0+2pq,\dots,-a_q+2pq)}j_!\F^{\K}V^*)[2pq](pq).\end{array}$$
D'après la proposition \ref{prop3:t_structures_recollees},
$$w_{\geq (-a_0+2pq,\dots,-a_q+2pq)}j_!\F^{\K}V^*=w_{\leq (-a_0+2pq,-a_1+2pq-1,\dots,-a_q+2pq-1)}Rj_*\F^{\K}V^*.$$
Notons $b_0=m+pq$ et, pour tout $r\in\{1,\dots,q\}$, $b_r=-(s_r-m)+(p-r)(q-r)$.
Alors $b_0=-a_0+2pq$ et, pour tout $r\in\{1,\dots,q\}$, $b_r=-a_r-1+2pq$.
Donc
$$D(W^{\geq t_1,\dots,\geq t_q}(V))= (w_{\leq (b_0,\dots,b_q)}Rj_*\F^{\K}V^*)[2pq](pq)=W^{\geq s_1,\dots,\geq s_q}(V^*)[2pq](pq).$$
\item[(2)] On peut supposer que $V$ est pure. Soit $m$ son poids.
Le faisceau $\F^{\K}V$ est lisse pur de poids $-m$ sur $M_0$, qui est lisse de dimension $pq$, 
donc le seul faisceau de cohomologie perverse non nul de $\F^{\K}V$
est $\Hp^{pq}\F^{\K}V=\F^{\K}V[pq]$, qui est de poids $-m+pq$.
D'après la proposition \ref{prop3:t_structures_recollees},
$$\begin{array}{rcl}IC^{\K}V[pq] = j_{!*}(\F^{\K}V[pq]) & = & w_{\leq (-m+pq,\dots,-m+pq)}Rj_*\F^{\K}V[pq] \\
& = & w_{\leq (-m+pq,-m+pq-1,\dots,-m+pq-1)}Rj_*\F^{\K}V[pq],\end{array}$$
donc 
$$IC^{\K}V=w_{\leq (-m+pq,\dots,-m+pq)}Rj_*\F^{\K}V=w_{\leq (-m+pq,-m+pq-1,\dots,-m+pq-1)}Rj_*\F^{\K}V.$$
Pour conclure, il suffit de remarquer que, pour tout $r\in\{1,\dots,q\}$,\newline
 $-t_r-m+(p-r)(q-r)=-m+pq-1$ et $-s_r-m+(p-r)(q-r)=-m+pq$.

\end{itemize}
\end{preuve}

\section{Restrictions des complexes pondérés aux strates}

On note $\S=\Gm.I_{p+q}$ le centre déployé de $\G$ et, pour tout $r\in\{1,\dots,q\}$, 
$$\S_r=\left\{\left(\begin{array}{ccc}\lambda^2 I_r & 0 & 0 \\
0 & \lambda I_{p-r} & 0 \\
0 & 0 & I_r\end{array}\right),\lambda\in\Gm\right\}$$
$\S_r$ est le centre déployé de $\G_r$.
On note $\chi_r$ le caractère de $\S_r$ défini par
$$\chi_r\left(\left(\begin{array}{ccc}\lambda^2I_r & 0 & 0 \\
0 & \lambda I_{p-r} & 0 \\
0 & 0 & I_r\end{array}\right)\right)=\lambda.$$
On identifie $X^*(\S_r)$ à $\Z$ en envoyant $\chi_r$ sur $1$, ce qui donne un ordre sur $X^*(\S_r)$.

\begin{notation} Soient $I\subset\{1,\dots,q\}$ et $V\in Rep_{\L_I(\Q_\ell)}$, avec $\L_I=\P_I/\N_I$ (cf 1.5).

Si $(t_i)_{i\in I}\in (\Z\cup\{\pm\infty\})^I$, on note
$$V_{<t_i,i\in I}$$
le sous-espace vectoriel de $V$ sur lequel, pour tout $i\in I$, $\S_i$ agit par des caractères $<t_i$.
Comme les $\S_i$ sont centraux dans $\L_I$, $V_{<t_i,i\in I}$ est stable par $\L_I(\Q_\ell)$.

La définition ci-dessus s'étend trivialement aux complexes et donne un foncteur exact
$$\left\{\begin{array}{rcl}D^b(Rep_{\L_I(\Q_\ell)}) & \fl & D^b(Rep_{\L_I(\Q_\ell)}) \\
V & \fle & V_{<t_i,i\in I}\end{array}\right..$$ 

On définit de même $V_{\geq t_i,i\in I}$.

\end{notation}

On se place dans la situation de la définition \ref{def:complexes_ponderes}.

On utilisera les notations suivantes (qui sont celles de 1.5) :
Soit $S\subset\{1,\dots,q\}$ non vide.
On pose
$$\P_S=\bigcap_{s\in S}\P_s.$$
C'est un sous-groupe parabolique standard de $\G$,
dont on note $\N_S$ le radical unipotent et $\L_S=\P_S/\N_S$ le quotient de Levi.
Soit $r=max(S)$.
Rappelons (voir la section 1.4) qu'on a trois sous-groupes distingués $\L'_{\lin,r}$, $\L_{\lin,r}$ et $\G_r=\QP_r/\N_r$ de $\L_r$
tels que $\G_r$ et $\L'_{\lin,r}$ soient d'intersection triviale,
que $\L_r$ soit produit quasi-direct de $\L_{\lin,r}$ et $\G_r$,
que $\L_{\lin,r}$ et $\L'_{\lin,r}$ soient égaux si $1\leq r<q$,
et que $\L_{\lin,q}=\L'_{\lin,q}\times\SU(p-q)$.

On a $\QP_r\subset\P_S\subset\P_r$, donc
$\P_S/\N_r$ est produit quasi-direct de $\G_r$ et d'un sous-groupe parabolique $\P_{\lin,S}$
de $\L_{\lin,r}$. 
On note $\P'_{\lin,S}$ le sous-groupe parabolique correspondant de $\L'_{\lin,r}$,
et $\L_{\lin,S}$ et $\L'_{\lin,S}$ les quotients de Levi respectifs de $\P_{\lin,S}$ et $\P'_{\lin,S}$.
On a $\L'_{\lin,S}=\L_{\lin,S}$ si $q\not\in S$, 
et $\L_{\lin,S}=\L'_{\lin,S}\times\SU(p-q)$ si $q\in S$.

Pour tout $g\in\G(\Af)$, on pose
$$\Hr_{g,S}=g\K g^{-1}\cap\P_S(\Q)\QP_r(\Af)=g\K g^{-1}\cap\P'_{\lin,S}(\Q)\QP_r(\Af)$$
$$\Hr_{g,\lin,S}=g\K g^{-1}\cap Cent_{\P_S(\Q)}(\X_r)\N_r(\Af)=g\K g^{-1}\cap\P'_{\lin,S}(\Q)\N_r(\Af)$$
$$\K_{g,S}=\Hr_{g,S}/\Hr_{g,\lin,S}$$
$$\Gamma_{g,S}=\Hr_{g,\lin,S}/(g\K g^{-1}\cap\N_S(\Q)\N_r(\Af))=(g\K g^{-1}\cap\P'_{\lin,S}(\Q)\N_S(\Af))/(g\K g^{-1}\cap\N_S(\Af)).$$
$\K_{g,S}$ s'identifie à un sous-groupe ouvert compact net de $\G_r(\Af)$,
et $\Gamma_{g,S}$ à un sous-groupe arithmétique net de $\L'_{\lin,S}(\Q)$ (ou de $\L_{\lin,S}(\Q)$).

\begin{theoreme}\label{th:restriction_W_bord} Soient $t_1,\dots,t_q\in\Z$ et $V\in D^b(Rep_{\G(\Q_\ell)})$. On suppose que tous les $H^i(V)$ sont purs de même poids $m\in\Z$.
On fixe $r\in\{1,\dots,q\}$ et $g\in\G(\Af)$, 
et on note $i$ l'inclusion de la strate $M^{\K_{g,\{r\}}}(\G_r,\X_r)$ dans $M^{\K}(\G,\X)$.
Alors
\vspace{1cm}
\begin{flushleft}$\displaystyle{[i^*W^{\geq t_1,\dots,\geq t_q}(V)]=}$\end{flushleft}
\begin{flushright}$\displaystyle{\F^{\K_{g,\{r\}}}\left(\sum_{S}\sum_{i\in I_S}\left[Ind_{\K_{p_ig,S\cup\{r\}}}^{\K_{g,\{r\}}}(-1)^{card(S)}R\Gamma\left(\Gamma_{p_ig,S\cup\{r\}},R\Gamma(Lie(\N_{S\cup\{r\}}),V)_{\geq t_r+m,< t_s+m,s\in S}\right)\right]\right),}$\end{flushright}
où $S$ parcourt l'ensemble des sous-ensembles de $\{1,\dots,r-1\}$, et,
pour tout $\protect{S\subset\{1,\dots,r-1\}}$,
on choisit un système de représentants $(p_i)_{i\in I_S}$ du double quotient\newline
$\P_{S\cup\{r\}}(\Q)\QP_r(\Af)\sous\P_r(\Q)\QP_r(\Af)/\Hr_{g,\{r\}}$.

Supposons que $V$ est concentré en degré $0$. Alors,
si $t_s=s(s-n)$ pour tout $s\in\{1,\dots,q\}$, ou
si $t_s=s(s-n)+1$ pour tout $s\in\{1,\dots,q\}$, on a
$IC^{\K}V=W^{\geq t_1,\dots,\geq t_q}(V)$
(proposition \ref{calcul_IC} (2)).
On obtient donc en particulier une formule pour $[i^*IC^{\K}V]$.

\end{theoreme}

\begin{preuve} Le théorème résulte du théorème \ref{th:simplification_formule_traces},
de la proposition ci-dessous,
et du (ii) de la proposition \ref{prop:Shapiro}
(qui permet de remplacer les images directes par des induites).

\end{preuve}

Rappelons qu'on a noté, pour tout $r\in\{1,\dots,q\}$, $M_r$ l'union des strates de bord
de $M^{\K}(\G,\X)^*$ associées à $(\QP_r,\Y_r)$. On note $i_r$ l'inclusion de $M_r$
dans $M^{\K}(\G,\X)^*$.

\begin{proposition}\label{restr_tronques_bord} Soient $r_1,\dots,r_c\in\{1,\dots,q\}$ tels que
$r_1<\dots <r_c$,
$a_1,\dots,a_c\in\Z\cup\{\pm\infty\}$, $V\in D^b(Rep_{\G(\Q_\ell)})$, $r\in\{r_c,\dots,q\}$ et $g\in\G(\Af)$.
Pour tout $i\in\{1,\dots,c\}$, on pose $t_i=-a_i+(p-r_i)(q-r_i)$.
On note 
$$L=Ri_{r_c*}w_{>a_c}i_{r_c}^*\dots Ri_{r_1*}w_{>a_1}i_{r_1}^*Rj_*\F^{\K}V$$
et $i$ l'inclusion de la strate $M^{\K_{g,\{r\}}}(\G_r,\X_r)$ dans $M^{\K}(\G,\X)^*$.
Soit
$$S=\{r_1,\dots,r_c,r\}.$$
Alors on a un isomorphisme canonique
$$i^*L\simeq\bigoplus_C L_C$$
où $C$ parcourt l'ensemble des doubles classes de 
$\P_S(\Q)\QP_r(\Af)\sous\P_r(\Q)\QP_r(\Af)/\Hr_{g,\{r\}}$,
et de plus, si $b\in\P_r(\Q)\QP_r(\Af)$ est un représentant de la double classe $C$,
on a un isomorphisme
$$L_C\simeq T_{b*}\F^{\K_{bg,S}}R\Gamma(\Gamma_{bg,S},R\Gamma(Lie(\N_S),V)_{<t_1,\dots,<t_c}).$$

\end{proposition}

\begin{preuve} Si $r>r_c$, on pose $d=c+1$, $r_{d}=r$, $a_{d}=-\infty$ et $t_{d}=+\infty$;
si $r=r_c$, on pose $d=c$.

On raisonne par récurrence sur $d$. 

Si $d=1$, le résultat cherché est une conséquence immédiate du théorème de Pink et du lemme \ref{relation_tronques}.

Supposons donc $d\geq 2$, et supposons le résultat prouvé pour tout $1\leq d'<d$.
On note
$$S'=\{r_1,\dots,r_{d-1}\}=S-\{r\}.$$
Calculons d'abord
$$M=i^*Ri_{r_{d-1}*}w_{>a_{d-1}}i_{r_{d-1}}^*\dots Ri_{r_1*}w_{>a_1}i_{r_1}^*Rj_*\F^{\K}(V).$$
Pour tout $h\in\G(\Af)$, soit $i_h$ l'inclusion de $M^{\K_{h,\{r_{d-1}\}}}(\G_{r_{d-1}},\X_{r_{d-1}})$ dans $M^{\K}(\G,\X)^*$. 
Le complexe $N=Ri_{r_{d-1}*}w_{>a_{d-1}}i_{r_{d-1}}^*\dots Ri_{r_1*}w_{>a_1}i_{r_1}^*Rj_*\F^{\K}(V)$ est égal à
$$\bigoplus_{h\in\P_{r_{d-1}}(\Q)\QP_{r_{d-1}}(\Af)\sous\G(\Af)/\K}Ri_{h*}i_h^*N,$$
donc on a
$$M=\bigoplus_{h\in\P_{r_{d-1}}(\Q)\QP_{r_{d-1}}(\Af)\sous\G(\Af)/K}i^*Ri_{h*}i_h^*N.$$
Soit $h\in\G(\Af)$. D'après l'hypothèse de récurrence, 
on a un isomorphisme
$$i_h^*N\simeq\bigoplus_{C'}N_{C'},$$
où $C'$ parcourt l'ensemble des doubles classes dans 
$\P_{S'}(\Q)\QP_{r_{d-1}}(\Af)\sous\G(\Af)/\K$, et, si $b$ est un représentant de $C'$,
$$N_{C'}\simeq T_{b*}\F^{\K_{bh,S'}}R\Gamma(\Gamma_{bh,S'},R\Gamma(Lie(\N_{S'}),V)_{<t_1,\dots,<t_{d-1}}).$$
Fixons $b\in\P_{r_{d-1}}(\Q)\QP_{r_{d-1}}(\Af)$ et calculons
$$i^*Ri_{h*}T_{b*}\F^{\K_{bh,S'}}R\Gamma(\Gamma_{bh,S'},R\Gamma(Lie(\N_{S'}),V)_{<t_1,\dots,<t_{d-1}}).$$
On définit $X$ par le carré cartésien suivant :
$$\xymatrix{X\ar[r]\ar[d] & M^{\K_{bh,S'}}(\G_{r_{d-1}},\X_{r_{d-1}})^*\ar[d]^{\overline{T}_b} \\
M^{\K_{g,\{r\}}}(\G_r,\X_r)\ar[r] & \overline{M^{\K_{h,\{r_{d-1}\}}}(\G_{r_{d-1}},\X_{r_{d-1}})}}$$
où $\overline{M^{\K_{h,\{r_{d-1}\}}}(\G_{r_{d-1}},\X_{r_{d-1}})}$ est l'adhérence de $M^{\K_{h,\{r_{d-1}\}}}(\G_{r_{d-1}},\X_{r_{d-1}})$ dans $M^{\K}(\G,\X)^*$.
Alors on a
$$X\simeq\coprod_{j\in J}M^{\K_{q_jbh,S}}(\G_r,\X_r),$$
et le morphisme $X\fl M^{\K_{g,\{r\}}}(\G_r,\X_r)$ est $\sqcup T_{q_j}$,
avec $(q_j)_{j\in J}$ un système fini (éventuellement vide) de représentants pour une certaine relation d'équivalence sur $\P_r(\Q)\QP_r(\Af)$, qu'on précisera plus tard.
Soit $q\in\P_r(\Q)\QP_r(\Af)$, et calculons la restriction à la strate $M^{\K_{qbh,S}}(\G_r,\X_r)$ de $M^{\K_{bh,S'}}(\G_{r_{d-1}},\X_{r_{d-1}})^*$ du complexe 
$$L'=\F^{\K_{bh,S'}}R\Gamma(\Gamma_{bh,S'},R\Gamma(Lie(\N_{S'}),V)_{<t_1,\dots,<t_{d-1}}).$$
Notons $\RP$ le sous-groupe parabolique de $\G_{r_{d-1}}$ correspondant à $\P_r$, c'est-à-dire\newline
 $\RP=(\QP_{r_{d-1}}\cap\P_r)/\N_{r_{d-1}}$.
Soient $\N_R=\N_S/\N_{S'}$ le radical unipotent de $\RP$, $(\QP_R,\Y)$ la composante rationnelle de bord de $(\G_{r_{d-1}},\X_{r_{d-1}})$ associée à $\RP$.
On a $(\QP_R,\Y)/\N_R=(\G_r,\X_r)$, et la partie linéaire $\RP_\lin$ de $\RP/\N_R$ est telle que $\L_{\lin,S}=\L_{\lin,S'}\times\RP_\lin$.
de plus, on vérifie facilement que
$$(q\K_{bh,S'}q^{-1}\cap\RP(\Q)\QP_R(\Af))/(q\K_{bh,S'}q^{-1}\cap\RP_\lin(\Q)\N_R(\Af))=\K_{qbh,S},$$
et que le sous-groupe arithmétique net 
$$\Gamma=(q\K_{bh,S'}q^{-1}\cap\RP_\lin(\Q)\N_R(\Af))/(q\K_{bh,S'}q^{-1}\cap\N_R(\Af))$$
de $\RP_\lin(\Q)$ s'identifie à 
$$\Gamma_{qbh,S}/\Gamma_{bh,S'}.$$
D'autre part, comme l'action de $\S_{r_1},\dots,\S_{r_{d-1}}$ sur $\N_R$ est triviale, on a
$$R\Gamma(Lie(\N_R),R\Gamma(Lie(\N_{S'},V)_{<t_1,\dots,<t_{d-1}}))\simeq R\Gamma(Lie(\N_S),V)_{<t_1,\dots,<t_{d-1}}.$$ 
Finalement, en utilisant le théorème de Pink, on voit que la restriction de $L'$ à $M^{\K_{qbh,S}}(\G_r,\X_r)$ est isomorphe à
$$\F^{\K_{qbh,S}}R\Gamma(\Gamma_{qbh,S},R\Gamma(Lie(\N_S),V)_{<t_1,\dots,<t_{d-1}}).$$
Il reste à compter les diagrammes
$$\xymatrix{M^{\K_{qbh,S}}(\G_r,\X_r)\ar[r]\ar[d] & M^{\K_{bh,S'}}(\G_{r_{d-1}},\X_{r})^*\ar[d] \\
M^{\K_{g,\{r\}}}(\G_r,\X_r)\ar[r] & \overline{M^{\K_{h,\{r_{d-1}\}}}(\G_{r_{d-1}},\X_{r_{d-1}})}}$$
modulo une relation d'équivalence convenable.
C'est ce qui a été fait dans la proposition \ref{combi_strates}.
La conclusion de cette proposition, combinée avec les calculs ci-dessus, montre que, si $(b_i)_{i\in I}$ est un système de représentants du double quotient $\P_S(\Q)\QP_r(\Af)\sous\P_r(\Q)\QP_r(\Af)/\Hr_{g,\{r\}}$, on a un isomorphisme
$$i^*M\simeq\bigoplus_{i\in I}T_{b_i*}\F^{\K_{b_ig,S}}R\Gamma(\Gamma_{b_ig,S},R\Gamma(Lie(\N_S),V)_{<t_1,\dots,<t_{d-1}})).$$
Le résultat cherché résulte alors  
du lemme \ref{relation_tronques}.

\end{preuve}

\begin{remarque} Nous pouvons maintenant rendre plus explicite le rapport entre les complexes pondérés
définis ici
et ceux de [GHM].
Soient $t_1,\dots,t_q\in\Z\cup\{\pm\infty\}$
et $V$ une représentation algébrique de $\G$,
qu'on suppose pure de poids $0$ pour simplifier.
Pour tout $r\in\{1,\dots,q\}$, on note
 $a_r=-t_r+(p-r)(q-r)$.
Alors, d'après le théorème \ref{th:simplification_formule_traces},
on a une égalité dans le groupe de Grothendieck de $D^b_m(M^{\K}(\G,\X)^*,\Q_\ell)$
$$[W^{\geq t_1,\dots,\geq t_q}V(\Q_\ell)]=\sum_{1\leq r_1<\dots<r_c\leq q}(-1)^c[Ri_{r_c*}w_{>a_{r_c}}i_{r_c}^*\dots Ri_{r_1*}w_{>a_1}i_{r_1}^*Rj_*\F^{\K}V(\Q_\ell)].$$
Or, d'après le calcul explicite des
$$L_{r_1,\dots,r_c}=Ri_{r_c*}w_{>a_{r_c}}i_{r_c}^*\dots Ri_{r_1*}w_{>a_{r_1}}i_{r_1}^*Rj_*\F^{\K}V(\Q_\ell)$$
qui a été fait dans la proposition ci-dessus,
il existe une manière naturelle de relever ces complexes en des complexes
sur un modèle entier $\Mod^{\K}_\Sigma(\G,\X)^*$, qu'on notera encore $L_{r_1,\dots,r_c}$.
Notons $L_{r_1,\dots,r_c}(\C)$ le complexe de faisceaux de $\Q_\ell$-espaces vectoriels
sur $M^{\K}(\G,\X)^*(\C)$ déduit du complexe
$L_{r_1,\dots,r_c}$ sur $\Mod^{\K}_\Sigma(\G,\X)^*$.
Alors la somme alternée
$$\sum_{1\leq r_1<\dots<r_c\leq q}(-1)^c[L_{r_1,\dots,r_c}(\C)]$$
est égale à la classe (dans le groupe de Grothendieck
de la catégorie dérivée de la catégorie des faisceaux de $\Q_\ell$-espaces vectoriels
sur $M^{\K}(\G,\X)^*(\C)$)
de l'image directe par le morphisme canonique de la compactification de Borel-Serre réductive
de $M^{\K}(\G,\X)(\C)$ sur $M^{\K}(\G,\X)^*(\C)$
du complexe pondéré de [GHM] associé au profil de poids $(t_1,\dots,t_q)$
et à coefficients dans $V$.

\end{remarque}

\section{La formule des traces}

\subsection{Rappels sur un théorème de Kottwitz}

Commençons par fixer quelques notations.

\begin{notation}Soient $p\not=\ell$ deux nombres premiers, $q=p^k$ une puissance de $p$, $\Fi_p\subset\Fi_q$ les corps finis à $p$ et $q$ éléments, $\Fi$ une clôture algébrique de $\Fi_q$ et $X$ un schéma de type fini sur $\Fi_q$.
On note $F$ l'endomorphisme de Frobenius de $X$, qui est l'identité sur l'espace topologique sous-jacent et l'élévation à la puissance $q$ sur les sections du faisceau structural.
Sur $X(\Fi)$, $F$ agit comme la substitution de Frobenius $\varphi:x\fle x^q$, donc $X(\Fi)^F=X(\Fi_q)$.

Considérons maintenant $K\in D^b_c(X,\Q_\ell)$. On a alors une correspondance de Frobenius
$$F^*:F^*K\iso K,$$
d'où un endomorphisme
$$F^*:R\Gamma_c(X_{\Fi},K_{\Fi})\fl R\Gamma_c(X_\Fi,K_\Fi),$$
qui est l'inverse de l'endomorphisme que $\varphi$ induit par transport de structure.
On note
$$Tr(F^*,R\Gamma_c(X_{\Fi},K_{\Fi}))=\sum_{i\in\Z}(-1)^i Tr(F^*,H^i_c(X_{\Fi},K_{\Fi})).$$

Soient $x$ un point fermé de $X$ et $n=[k(x):\Fi_q]$. Si $\overline{x}$ est un point géométrique de $X(\Fi)$ localisé en $x$, c'est un point fixe de $F^n$, et on note $F_x^*$ l'endomorphisme de $K_{\overline{x}}$ induit par $F^n$. 
La classe d'isomorphisme du couple $(K_{\overline{x}},F_x^*)$ ne dépend pas du choix de $\overline{x}$, donc la trace 
$$Tr(F_x^*,K_{\overline{x}})=\sum_{i\in\Z}(-1)^*Tr(F_x^*,H^i(K)_{\overline{x}})$$
ne dépend que de $x$; on la note $Tr(F_x^*,K_x)$.
\end{notation}

Nous allons énoncer le résultat principal de l'article [K2] pour la donnée de Shimura $(\G,\X)$ de la section 1.2, 
dans le cas particulier où la correspondance de Hecke considérée est l'identité (c'est-à-dire où, avec les notations de [K2], $g=1$).

\begin{theoreme}\label{th_Kottwitz} On se place dans la situation de la section 3.1 : 
$(\K,\Sigma)$ est comme dans la définition \ref{def:prob_mod}, 
$p$ est un nombre premier tel que $p\not\in \Sigma$, et on s'intéresse à la réduction modulo $p$ de $\Mod^{\K}_\Sigma(\G,\X)$, qu'on notera $M^{\K}(\G,\X)$ ou simplement $M$.
En particulier, $\K=\K_{0,p}\K^p$, avec $\K^p\subset\G(\Af^p)$.
Soient $\ell\in \Sigma$ et $V\in D^b(Rep_{\G(\Q_\ell)})$.
On note $K=\F^{\K}V$. 
Alors, pour tout $j\in\Nat^*$,
$$\sum_{x\in M(\Fi)^{F^j}}Tr(F_x^*,K_x)=\sum_{(\gamma_0;\gamma,\delta)}c(\gamma_0;\gamma,\delta)O_\gamma(f^p)TO_\delta(\phi)Tr(\gamma_0,V).$$

\nopreuve

\end{theoreme}

Il faut encore expliquer toutes les notations.
Rappelons qu'on a fixé une injection $E\fl\overline{\Q}_p$. 
On note $\wp$ l'idéal premier de $\O_E$ au-dessus de $p$ déterminé par cette injection, $L$ l'extension non ramifiée de degré $j$ de $E_\wp$, $\Q_p^{nr}$ l'extension non ramifiée maximale de $\Q_p$ et $\sigma\in Gal(\Q_p^{nr}/\Q_p)$ l'élément correspondant au Frobenius arithmétique $x\fle x^p$ de $\Fi$.

Dans la somme, le triplet $(\gamma_0;\gamma,\delta)$ parcourt un système de représentants des classes d'équivalence de triplets formés de $\gamma_0\in\G(\Q)$ semi-simple elliptique dans $\G(\R)$, de $\gamma=(\gamma_v)\in\G(\Af^p)$ et de $\delta\in\G(L)$ vérifiant les conditions de [K1] 2
(en particulier, pour toute place $v\not=p,\infty$ de $\Q$, $\gamma_0$ et $\gamma_v$ sont conjugués sous $\G(\overline{\Q}_v)$, et $\gamma_0$ et $N\delta$ sont conjugués sous $\G(\overline{\Q}_p)$, où $N\delta\in\G(\Q_p)$ est la norme de $\delta$, définie par
$N\delta=\delta.\sigma(\delta)\dots\sigma^{d-1}(\delta)$, avec $d=[L:\Q_p]$),
et tels que 
$$\alpha(\gamma_0;\gamma,\delta)=1,$$
où $\alpha$ est défini dans [K1] 2,
pour la relation d'équivalence suivante : deux triplets $(\gamma_0;\gamma,\delta)$ et $(\gamma_0';\gamma',\delta')$ sont équivalents \ssi
$\gamma_0$ et $\gamma_0'$ sont conjugués sous $\G(\overline{\Q})$, $\gamma$ et $\gamma'$ sont conjugués sous $\G(\Af^p)$ et $\delta$ et $\delta'$ sont $\sigma$-conjugués sous $\G(L)$ (c'est-à-dire qu'il existe $x\in\G(L)$ tel que $\delta'=x\delta\sigma(x)^{-1}$).

Rappelons qu'on a fixé un $\O_E$-réseau autodual $\Lambda=\O_E^{p+q}$ de $E^{p+q}$.
On note $\K_0$ le stabilisateur dans $\G(L)$ du réseau $\Lambda\otimes_\Z\O_L$. 

Soit $(\gamma_0;\gamma,\delta)$ un triplet comme ci-dessus.

On note $I_0$ le centralisateur de $\gamma_0$ dans $\G$.
Dans [K1] 3, Kottwitz montre qu'il existe une forme intérieure $I$ de $I_0$ telle que $I(\Af^p)$ soit le centralisateur de $\gamma$ dans $\G(\Af^p)$
et $I(\Q_p)$ le centralisateur tordu de $\delta$ dans $\G(L)$ 
(c'est-à-dire l'ensemble des $x\in\G(L)$ tels que $x\delta=\delta\sigma(x)$).

On pose
$$c(\gamma_0;\gamma,\delta)=vol(I(\Q)\sous I(\Af)).$$
(La définition de Kottwitz est $c(\gamma_0;\gamma,\delta)=vol(I(\Q)\sous I(\Af)).card(ker(ker^1(\Q,I_0)\rightarrow ker^1(\Q,\G)))$,
mais ici, $ker^1(\Q,\G)=\{1\}$ d'après [Sh] 5.8).

On choisit des mesures de Haar $dy$ et $dx$ sur $I(\Af^p)\sous\G(\Af^p)$ et $I(\Q_p)\sous\G(L)$.
Si $g$ est $\psi$ sont des fonctions localement constantes à support compact sur $I(\Af^p)\sous\G(\Af^p)$ et $I(\Q_p)\sous\G(L)$, on pose
$$O_\gamma(g)=\int_{I(\Af^p)\sous\G(\Af^p)}g(y^{-1}\gamma y)dy$$
$$TO_\delta(\psi)=\int_{I(\Q_p)\sous\G(L)}\psi(x^{-1}\delta\sigma(x))dx.$$
On note
$$f^p=\frac{\ungras_{\K^p}}{vol(\K^p)}$$
$$\phi=\frac{\ungras_{\K_0a\K_0}}{vol(\K_0)}.$$
Il reste à  expliquer qui est $a$ (cf la fin de [K1] 3).
On a défini dans 2.1.1 un morphisme
$$h_0r:{\Gm}_{,\C}\fl\G_\C$$
(noté $\mu_h$ dans [K1], où $L$ est noté $F$ et $F$ est noté $E$). 
Il est conjugué par $\G(\C)$ au morphisme 
$$\mu:{\Gm}_{,\C}\fl\G_\C,\quad z\fle\left(\begin{array}{cc}\left(1\otimes\frac{z+1}{2}+i\otimes\frac{z-1}{2i}\right)I_p & 0 \\
0 & \left(1\otimes\frac{z+1}{2}-i\otimes\frac{z-1}{2i}\right)I_q \end{array}\right).$$
L'image de $\mu$ arrive dans $T_\C$, où
$T$ est le sous-tore déployé sur $\O_L$ maximal du tore diagonal de $\GU(p,q)$, 
donc on peut prendre
$$a=\mu(\varpi_L^{-1}),$$
où $\varpi_L$ est une uniformisante de $\O_L$.

\begin{remarque} Le sous-groupe compact ouvert $\K$ de $\G(\Af)$ n'intervient que dans l'intégrale orbitale $O_\gamma(f^p)$.

\end{remarque}

\subsection{Formule des traces pour certains tronqués}

On se place dans la situation du paragraphe précédent.
Le but est maintenant d'écrire la formule des traces pour les complexes
$Ri_{r_c*}w_{>a_c}i_{r_c}^*\dots Ri_{r_1*}w_{>a_1}i_{r_1}^*Rj_*\F^{\K}V$.

La proposition suivante va nous permettre d'utiliser le théorème \ref{th_Kottwitz}.

\begin{proposition}\label{prop:formule_traces} Soient $S\subset\{1,\dots,r\}$, $r=max(S)$, $g\in\G(\Af)$, 
$a\in\P_r(\Q)\QP_r(\Af)$ et $W\in D^b(Rep_{\L_S(\Q_\ell)})$.
On note comme ci-dessus
$$\Gamma_{ag,S}=((ag)\K (ag)^{-1}\cap\P_S(\Q)\QP_r(\Af))/((ag)\K (ag)^{-1}\cap (\N_S\QP_r)(\Af)),$$
et soit
$$\begin{array}{rcl}\K^\circ_{ag,S} & = & ((ag)\K (ag)^{-1}\cap\N_S(\Q)\QP_r(\Af))/((ag)\K (ag)^{-1}\cap\N_S(\Q)\N_r(\Af)) \\
 & = & ((ag)\K (ag)^{-1}\cap (\N_S\QP_r)(\Af))/((ag)\K (ag)^{-1}\cap\N_S(\Af)).\end{array}$$
On a un revêtement étale fini
$$T_a:M'=M^{\K_{ag,S}}(\G_r,\X_r)\fl M=M^{\K_{g,\{r\}}}(\G_r,\X_r),$$
et on pose
$$M^\circ=M^{\K^\circ_{ag,S}}(\G_r,\X_r);$$
on a un revêtement étale fini
$$T_1:M^\circ\fl M'.$$
Enfin, on note
$$K=\F^{\K_{ag,S}}R\Gamma(\Gamma_{ag,S},W);$$
$T_1^*K$ est canoniquement isomorphe à $\F^{\K^\circ_{ag,S}}R\Gamma(\Gamma_{ag,S},W)$ (avec $\K^\circ_{ag,S}$ agissant sur $R\Gamma(\Gamma_{ag,S},W)$ via la représentation de $\G_r(\Q_\ell)$ dans ce complexe).

Alors, pour tout $j\in\Nat^*$, on a
$$\sum_{x\in M(\Fi)^{F^j}}Tr(F_x^*,(T_{a*}K)_x)=\frac{1}{[\K_{ag,S}:\K^\circ_{ag,S}]}\sum_{x\in M^\circ(\Fi)^{F^j}}Tr(F_x^*,(T_1^*K)_x).$$

\end{proposition}

\begin{preuve} Comme $T_a:M'\fl M$ est propre, on a pour tout $j\in\Nat^*$
$$\begin{array}{rcl}\sum_{x\in M(\Fi)^{F^j}}Tr(F_x^*,(T_{p*}K)_x) & = & Tr(F^{j*},R\Gamma_c(M_\Fi,(T_{a*}K)_\Fi)) \\
 & = & Tr(F^{j*},R\Gamma_c(M'_\Fi,K_\Fi)) \\
& = & \sum_{x\in M'(\Fi)^{F^j}}Tr(F_x^*,K_x)\end{array}.$$
La proposition résulte alors du lemme suivant, appliqué à $\G_r\subset\P_S/\N_S.$

\end{preuve}

\begin{lemme} On se place dans la situation de la définition \ref{def_faisceaux}, avec $\N=0$ : on a donc $\G,\G_\lin\subset\L$, une donnée de Shimura $(\G,\X)$, un sous-groupe compact ouvert net $\K_L\subset\L(\Af)$ et $W\in D^b(Rep_{\L(\Q_\ell)})$.
Notons 
$$\K=\K_L\cap\G(\Af)$$
$$\K'=\Hr/\Gamma_\lin=(\K_L\cap\L(\Q)\G(\Af))/(\K_L\cap\G_\lin(\Q)).$$
$\K$ s'envoie injectivement dans $\K'$, et son image est distinguée d'indice fini
égal au cardinal de $\Hr/\K\Gamma_\lin$.
On notera aussi $[\K':\K]$ cet indice.
$\K'$ agit sur $R\Gamma(\Gamma_\lin,W)$ via l'action de $\Hr$, comme le quotient $\Hr/\Gamma_\lin$, et $\K$ agit sur $R\Gamma(\Gamma_\lin,W)$ comme sous-groupe de $\G(\Af)$, qui agit lui-même sur $R\Gamma(\Gamma_\lin,W)$ parce qu'il commute avec $\Gamma_\lin$.
On se place toujours sur les réductions modulo un nombre premier $p$ assez grand.
Alors, pour tout $j\in\Nat^*$,
$$\sum_{x\in M^{\Hr/\Gamma_\lin}(\G,\X)(\Fi)^{F^j}}Tr(F_x^*,\F^{\K'}R\Gamma(\Gamma_\lin,W))=\frac{1}{[\K':\K]}\sum_{x\in M^{\K}(\G,\X)(\Fi)^{F^j}}Tr(F_x^*,\F^{\K}R\Gamma(\Gamma_\lin,W)).$$

\end{lemme}

\begin{preuve} Fixons $j\in\Nat^*$, et notons $\til{X}=M(\G,\X)$, $X'=M^{\Hr/\Gamma_\lin}(\G,\X)$, $X=M^{\K}(\G,\X)$. 
On a un diagramme commutatif
$$\xymatrix{\til{X}\ar[r]^{\varphi}\ar[dr]_{\varphi'} & X\ar[d]^{T_1} \\
& X'}$$
où $\varphi$ et $\varphi'$ sont les revêtements évidents.
$\varphi$ est un revêtement étale galoisien profini de groupe $\K$, et $\varphi'$ 
est un revêtement étale galoisien profini de groupe $\K'$.
Soit $x'\in X'(\Fi)^{F^j}$, et soit $y\in\til{X}(\Fi)$ tel que $x'=\varphi'(y)$. Comme $x'$ est stable par $F^j$, on peut écrire $F^j(y)=y.k$, avec $k\in\K'$ (uniquement déterminé). 
On a deux possibilités :
\begin{itemize}
\item[(1)] $x'$ est l'image par $T_1$ d'un élément de $X(\Fi)^{F^j}$, ce qui revient à demander que $k$ soit dans $\K$.
Alors les $[\K':\K]$ pré-images de $x'$ par $T_1$ sont dans $X(\Fi)^{F^j}$ (ce sont les $\varphi(y.l)$, pour $l\in\K'/\K$). 
De plus, comme $\F^{\K}R\Gamma(\Gamma_\lin,W)\simeq T_1^*\F^{\K'}R\Gamma(\Gamma_\lin,W)$, pour tout $x\in X(\Fi)$ tel que $T_1(x)=x'$, on a
$$Tr(F_x^*,\F^{\K}R\Gamma(\Gamma_\lin,W))=Tr(F_{x'}^*,\F^{\K'}R\Gamma(\Gamma_\lin,W)).$$
On en déduit que
$$Tr(F_{x'}^*,\F^{\K'}R\Gamma(\Gamma_\lin,W))=\frac{1}{[\K':\K]}\sum_{x\in X(\Fi)^{F^j},T_1(x)=x'}Tr(F_x^*,\F^{\K}R\Gamma(\Gamma_\lin,W)).$$
\item[(2)] $x'$ n'est pas l'image d'un élément de $X(\Fi)^{F^j}$, ce qui revient à demander que $k\not\in\K$.
Nous allons montrer que dans ce cas
$$Tr(F_{x'}^*,\F^{\K'}R\Gamma(\Gamma_\lin,W))=0,$$
ce qui finira la preuve du lemme.
On choisit $h\in\Hr=\K_L\cap\L(\Q)\G(\Af)$ tel que $h\Gamma_\lin=k$.
Comme $k\not\in\K$, $h$ n'est pas dans $\Gamma_\ell\K$.
D'après la proposition \ref{fibres_faisceaux}, on a
$$Tr(F_{x'}^*,\F^{\K'}R\Gamma(\Gamma_\lin,W))=Tr(k,R\Gamma(\Gamma_\lin,W))=Tr(h,R\Gamma(\Gamma_\lin,W)).$$
Il suffit d'appliquer le lemme ci-dessous.
\end{itemize}
\end{preuve}

\begin{lemme}\label{lemme:annulation_trace} Soient $\L$ un groupe réductif connexe, $\G,\G_\lin$ deux sous-groupes réductifs de $\L$
tels que $\L$ soit produit quasi-direct de $\G$ et $\G_\lin$
et $\K_L$ un sous-groupe ouvert compact net de $\L(\Af)$.
On suppose qu'il existe un sous-groupe distingué $\G'_\lin\subset\G_\lin$
et un sous-groupe $\G'\subset\L$ qui contient $\G$ comme un sous-groupe distingué
tels que $\G_\lin/\G'_\lin$ et $\G'/\G$ soient de type compact 
et que $\L$ soit produit direct de $\G'_\lin$ et $\G'$.
On note $\Hr=\K_L\cap\L(\Q)\G(\Af)$, $\Gamma_\ell=\K_L\cap\G_\lin(\Q)$ et 
$\K=\K_L\cap\G(\Af)$.
$\Gamma_\ell$, $\K$ et $\Gamma_\ell\K$ sont des sous-groupes distingués de $\Hr$.
Pour tout $W\in D^b(Rep_{\L(\Q_\ell)})$,
pour tout $h\in\Hr$, si $h\not\in \Gamma_\ell\K$, alors
$$Tr(h,R\Gamma(\Gamma_\ell,W))=0.$$

\end{lemme}

\begin{preuve} Comme $\K_L$ est net,
$\Hr=\K_L\cap\G'_\lin(\Q)\G(\Af)=\K_L\cap\G_\lin(\Q)\G(\Af)$ et 
$\Gamma_\lin=\G'_\lin(\Q)\cap\K_L$.
Notons $\Gamma'_\lin$ le projeté de $\Hr$ sur $\G'_\lin(\Q)$.
C'est un sous-groupe arithmétique net de $\G'_\lin(\Q)$ qui contient $\Gamma_\lin$ comme sous-groupe distingué.

Soient $h\in\Hr-\Gamma_\lin\K$ et $W\in D^b(Rep_{\L(\Q_\ell)})$.
On peut supposer que $W\in Rep_{\L(\Q_\ell)}$.
On écrit $h=\gamma.g$, avec $\gamma\in\G'_\lin(\Q)$ et $g\in\G(\Af)$.
Comme $h\not\in\Gamma_\lin\K$, on a $\gamma\in\Gamma'_\lin-\Gamma_\lin$.

Comme $\L$ est produit direct de $\G'_\lin$ et $\G'$, 
la représentation $W$ de $\L(\Q_\ell)$ est
somme directe de représentations de la forme $W_1\otimes W_2$, avec $W_1$ une représentation de $\G'_\lin(\Q_\ell)$ et $W_2$ une représentation de $\G'(\Q_\ell)$; 
pour une telle représentation $W_1\otimes W_2$, on a
$$R\Gamma(\Gamma_\lin,W_1\otimes W_2)=R\Gamma(\Gamma_\lin,W_1)\otimes W_2,$$
et
$$Tr(h,R\Gamma(\Gamma_\lin,W_1\otimes W_2))=Tr(\gamma,R\Gamma(\Gamma_\lin,W_1)).Tr(g,W_2).$$
Il suffit donc de prouver l'énoncé suivant : pour toute représentation $W_1$ de $\G'_\lin(\Q_\ell)$, on a
$$Tr(\gamma,R\Gamma(\Gamma_\lin,W_1))=0$$

Soit $Y$ l'espace symétrique du groupe $\G'_\lin$. 
Il a une compactification de Borel-Serre partielle $\overline{Y}^{BS}$, 
sur laquelle tout sous-groupe arithmétique de $\G'_\lin(\Q)$ agit proprement, 
et telle que $\Gamma_\lin\sous\overline{Y}^{BS}$ soit la compactification de Borel-Serre de l'espace localement symétrique $\Gamma_\lin\sous Y$ ([BS] 9.3).
Notons $j:\Gamma_\lin\sous Y\fl\Gamma_\lin\sous\overline{Y}^{BS}$ l'immersion ouverte.
Comme $\gamma$ normalise $\Gamma_\lin$, on a une correspondance de Hecke 
$$c_{\gamma}=(\overline{T}_{{\gamma}^{-1}},\overline{T}_1):\Gamma_\lin\sous\overline{Y}^{BS}\times\Gamma_\lin\sous\overline{Y}^{BS}\fl\Gamma_\lin\sous\overline{Y}^{BS}.$$
Soit $W_1$ une représentation de $\G_\lin(\Q_\ell)$. 
On peut lui associer un faisceau de $\Q_\ell$-espaces vectoriels $\F=\F^{\Gamma_\lin}W_1$ sur $\Gamma_\lin\sous Y$ vérifiant
$$R\Gamma(\Gamma_\lin\sous Y,\F)=R\Gamma(\Gamma_\lin\sous\overline{Y}^{BS},Rj_*\F)=R\Gamma(\Gamma_\lin,W_1),$$
et la correspondance $c_{\gamma}$ se relève en une correspondance cohomologique
$$u_{\gamma}:\overline{T}_{\gamma^{-1}}^*Rj_*\F\fl \overline{T}_1^!Rj_*\F=Rj_*\F$$
telle que
$$Tr(u_{\gamma},R\Gamma(\Gamma_\lin\sous\overline{Y}^{BS},Rj_*\F))=Tr(\gamma,R\Gamma(\Gamma_\lin,W_1)).$$
Pour montrer que $Tr(\gamma,R\Gamma(\Gamma_\lin,W_1))=0$, il suffit donc, d'après la formule des traces de Lefschetz, 
de montrer que la correspondance de Hecke $c_{\gamma}$ n'a pas de points fixes dans $\Gamma_\lin\sous\overline{Y}^{BS}$.
Supposons que $c_{\gamma}$ a un point fixe.
Il existe alors $y\in\overline{Y}^{BS}$ et $\delta\in\Gamma_\lin$ tels que $\gamma.y=\delta.y$, c'est-à-dire $(\delta^{-1}\gamma).y=y$.
Comme $\G_\lin(\Q)$ agit proprement sur $\overline{Y}^{BS}$, $\delta^{-1}\gamma$ est d'ordre fini; 
or $\delta^{-1}\gamma$ est un élément du sous-groupe arithmétique net $\Gamma'_\lin$ de $\G_\lin(\Q)$, donc $\delta^{-1}\gamma=1$, et $\gamma=\delta\in\Gamma_\lin$, 
ce qui contredit l'hypothèse.

\end{preuve}

\begin{corollaire}\label{formule_traces_sigma} Soient $r_1,\dots,r_c\in\{1,\dots,q\}$
tels que $r_1<\dots<r_c$,
$a_1,\dots,a_c\in\Z\cup\{\pm\infty\}$ et $V\in D^b(Rep_{\G(\Q_\ell)})$.
On note
$$t_i=-a_i+(p-r_i)(q-r_i)$$
pour tout $i\in\{1,\dots,c\}$ et
$$L=Ri_{r_c*}w_{>a_c}i_{r_c}^*\dots Ri_{r_1*}w_{>a_1}i_{r_1}^*Rj_*\F^{\K}V.$$
Enfin, soit $j\in\Nat^*$.
Alors
$$Tr(F^{j*},R\Gamma(M^{\K}(\G,\X)^*_\Fi,L_\Fi))=\sum_{s=r_c}^qT_s,$$
avec, pour tout $s\in\{r_c,\dots,q\}$, $T_s$ défini par
$$T_s=\sum_{(\gamma_0;\gamma,\delta)\in C_{s,j}}c(\gamma_0;\gamma,\delta)\chi(\L'_{\lin,S})TO_\delta(\phi^{(s)}_j)Tr(\gamma_0,R\Gamma(Lie(\N_S),V)_{<t_1,\dots,<t_c})$$
$$\sum_{i\in I}O_\gamma(f_{g_i})vol\left((g_i\K g_i^{-1}\cap\P_S(\Af))/(g_i\K g_i^{-1}\cap (\G'_s\N_S)(\Af))\right)^{-1},$$
où
\begin{itemize}
\item[$\bullet$] $S=\{r_1,\dots,r_c,s\}$.
\item[$\bullet$] $C_{s,j}$ est l'ensemble d'indices du théorème \ref{th_Kottwitz} pour le groupe $\G_s$; si $(\gamma_0;\gamma,\delta)\in C_{s,j}$, $c(\gamma_0;\gamma,\delta)$, $O_\gamma$ et $TO_\delta(\phi_j^{(s)})$ ont la même signification que dans l'énoncé du théorème \ref{th_Kottwitz};
\item[$\bullet$] $(g_i)_{i\in I}$ est un système de représentants de $\P_S(\Af)\sous\G(\Af)/\K$ dans $\G(\Af^p)$;
\item[$\bullet$] pour tout $i\in I$,
$$f_{g_i}=\frac{\ungras_{(g_i\K^p g_i^{-1}\cap\N_S(\Q)\QP_s(\Af^p))/(g_i\K^p g_i^{-1}\cap\N_S(\Q)\N_s(\Af^p))}}{vol((g_i\K^p g_i^{-1}\cap\N_S(\Q)\QP_s(\Af^p))/(g_i\K^p g_i^{-1}\cap\N_S(\Q)\N_s(\Af^p)))};$$
\item[$\bullet$] pour un groupe réductif $\H$, $\chi(\H)$ est défini dans [GKM] 7.10 : c'est un réel qui dépend du choix d'une mesure de Haar sur $\H(\Af)$ et qui est tel que, si $X_H$ est l'espace symétrique de $\H$ et $\K_H$ est un sous-groupe ouvert compact de $\H(\Af)$,
$$\chi(\H(\Q)\sous(X_H\times\H(\Af)/\K_H))=\chi(\H).vol(\K_H)^{-1}.$$
\end{itemize}
(Les diverses mesure de Haar utilisées plus haut doivent bien sûr être choisies de manière compatible.)
% Expliquons ce que nous entendons pas cela. 
% On sait que $\L_S=\L_{\lin,S}\times\G_s$; comme $\gamma\in\G_s(\Af^p)$, $I_L(\Af^p)=\L_{\lin,S}(\Af^p)\times I(\Af^p)$, où $I(\Af^p)$ est comme dans l'énoncé du théorème \ref{th_Kottwitz} le centralisateur de $\gamma$ dans $\G_s(\Af^p)$.
% La mesure de Haar $dy$ sur $I_L(\Af^p)\sous\L_S(\Af^p)$ est le quotient d'une mesure de Haar sur $\L_S(\Af^p)$ et d'une mesure de Haar sur $I_L(\Af^p)$.
% La première de ces mesures est le produit d'une mesure de Haar sur $\L_{\lin,S}(\Af^p)$ et d'une mesure de Haar sur $\G_s(\Af)$, et la deuxième est le produit de la même mesure sur $\L_{\lin,S}(\Af^p)$ et d'une mesure de Haar sur $I(\Af^p)$.
% Enfin, on utilise pour définir $\chi(\L_{\lin,S})$ une mesure de Haar sur $\L_{\lin,S}(\Af)$ compatible à la mesure déjà introduite sur $\L_{\lin,S}(\Af^p)$.

\end{corollaire}

\begin{preuve} Rappelons que, pour tout $s\in\{0,\dots,q\}$, $M_s$ est 
l'union des strates de bord de $M^{\K}(\G,\X)^*$ associées à $(\QP_s,\Y_s)$
et $i_s$ est l'inclusion de $M_s$ dans $M^{\K}(\G,\X)^*$.
On sait que $L$ est à support dans l'union des $M_s$, avec $r_c\leq s\leq q$, donc, d'après SGA 4 1/2 Rapport 3.2, on a
$$Tr(F^{j*},R\Gamma(M^{\K}(\G,\X)^*_\Fi,L_\Fi))=\sum_{s=r_c}^q T_s,$$
où
$$T_s=\sum_{x\in M_s(\Fi)^{F^j}}Tr(F_x^*,L_x).$$
Fixons $s\in\{r_c,\dots,q\}$.
Les propositions \ref{restr_tronques_bord} et \ref{prop:formule_traces} impliquent le résultat suivant : 
si $(g_i)_{i\in I}$ est un système de représentants de $\P_s(\Q)\QP_s(\Af)\sous\G(\Af)/\K$, 
et si, pour tout $i\in I$, $(p_{ij})_{j\in J_i}$ est un système de représentants de $\P_S(\Q)\QP_s(\Af)\sous\P_s(\Q)\QP_s(\Af)/(g_i\K g_i^{-1}\cap\P_s(\Q)\QP_s(\Af))$, alors
$$T_s=\sum_{i\in I}\sum_{j\in J_i}\frac{1}{[\K_{p_{ij}g_i,S}:\K^\circ_{p_{ij}g_i,S}]}\sum_{x\in M^{\K^0_{p_{ij}g_i,S}}(\G_s,\X_s)(\Fi)^{F^j}}Tr(F_x^*,L^{ij}_x),$$
où
$$L^{ij}=\F^{\K^\circ_{p_{ij}g_i,S}}R\Gamma(\Gamma_{p_{ij}g_i,S},R\Gamma(Lie(\N_S),V)_{<t_1,\dots,<t_c}).$$
Rappelons que, pour tout $g\in\G(\Af)$,
$$\K_{g,S}=(g\K g^{-1}\cap\P_S(\Q)\QP_s(\Af))/(g\K g^{-1}\cap\P_{\lin,S}(\Q)\N_s(\Af))$$
$$\K^\circ_{g,S}=(g\K g^{-1}\cap\N_S(\Q)\QP_s(\Af))/(g\K g^{-1}\cap\N_S(\Q)\N_s(\Af))$$
$$\Gamma_{g,S}=(g\K g^{-1}\cap\P_{\lin,S}(\Q)\N_s(\Af))/(g\K g^{-1}\cap\N_S(\Q)\N_s(\Af)).$$
On voit facilement que $(p_{ij}g_i)_{i\in I,j\in J_i}$ est un système de représentants de\newline
 $\P_S(\Q)\QP_s(\Af)\sous\G(\Af)/\K$, 
donc la formule ci-dessus se réécrit, en changeant les notations :
Soit $(g_i)_{i\in I}$ un système de représentants de $\P_S(\Q)\QP_s(\Af)\sous\G(\Af)/\K$.
Alors
$$T_s=\sum_{i\in I}\frac{1}{[\K_{g_i,S}:\K^\circ_{g_i,S}]}\sum_{x\in M^{\K^\circ_{g_i,S}}(\G_s,\X_s)(\Fi)^{F^j}}Tr(F_x^*,L^i_x),$$
où
$$L^i=\F^{\K^\circ_{g_i,S}}R\Gamma(\Gamma_{g,S},R\Gamma(Lie(\N_S),V)_{<t_1,\dots,<t_c}).$$
On utilise ensuite le théorème de Kottwitz (\ref{th_Kottwitz}).
Il implique que, pour tout $i\in I$,
\begin{flushleft}$\displaystyle{\sum_{x\in M^{\K^\circ_{g_i,S}}(\Fi)^{F^j}}Tr(F_x^*,L_x^i)=}$\end{flushleft}
\begin{flushright}$\displaystyle{\sum_{(\gamma_0;\gamma,\delta)\in C_{j,s}s}c(\gamma_0;\gamma,\delta)O_\gamma(f_{g_i})TO_\delta(\phi_j^{(s)})Tr(\gamma_0,R\Gamma(\Gamma_{g_i,S},R\Gamma(Lie(\N_S),V)_{<t_1,\dots,<t_c})).}$\end{flushright}
Comme
$$Tr(\gamma_0,R\Gamma(\Gamma_{g_i,S},R\Gamma(Lie(\N_S),V)_{<t_1,\dots,<t_c}))=\chi(\Gamma_{g_i,S})Tr(R\Gamma(Lie(\N_S),V)_{<t_1,\dots,<t_c})$$
et que ni $c(\gamma_0;\gamma,\delta)$ ni $TO_\delta(\phi_j^{(s)})$ ne dépendent de $g_i$, on obtient
\begin{flushleft}$\displaystyle{T_s=\sum_{(\gamma_0;\gamma,\delta)}c(\gamma_0;\gamma,\delta)TO_\delta(\phi_j^{(s)})\times}$\end{flushleft}
\begin{flushright}$\displaystyle{Tr(\gamma_0,R\Gamma(Lie(\N_S),V)_{<t_1,\dots,<t_c})\sum_{i\in I}\frac{1}{[\K_{g_i,S}:\K^\circ_{g_i,S}]}\chi(\Gamma_{g_i,S})O_\gamma(f_{g_i}).}$\end{flushright}
Notons, pour tout $g\in\G(\Af)$,
$$\Gamma^\circ_{g,S}=(g\K g^{-1}\cap\P_S(\Q)\QP_s(\Af))/(g\K g^{-1}\cap\N_S(\Q)\QP_s(\Af)).$$
Alors, pour tout $g\in\G(\Af)$,
$$[\K_{g,S}:\K^\circ_{g,S}]=[\Gamma^\circ_{g,S}:\Gamma_{g,S}],$$
donc
$$\frac{1}{[\K_{g,S}:\K^\circ_{g,S}]}\chi(\Gamma_{g,S})=\chi(\Gamma^\circ_{g,S}).$$

Changeons de notation, et notons $(g_i)_{i\in I}$ un système de représentants dans $\G(\Af^p)$ de $\P_S(\Af)\sous\G(\Af)/\K$.
Pour tout $i\in I$, soit $(p_{ij})_{j\in J_i}$ un système de représentants dans $\P_{\lin,S}(\Af)$ de 
$$\P_S(\Q)\QP_s(\Af)\sous\P_S(\Af)/(g_i\K g_i^{-1}\cap\P_S(\Af))=\L'_{\lin,S}(\Q)\sous\L'_{\lin,S}(\Af)/\K_{g_i,\lin,S},$$
où
$$\K_{g_i,\lin,S}=(g_i\K g_i^{-1}\cap\P_S(\Af))/(g_i\K g_i^{-1}\cap(\G'_s\N_S)(\Af)).$$
Alors $(p_{ij}g_i)_{i\in I,j\in J_i}$ est un système de représentants de $\P_S(\Q)\QP_s(\Af)\sous\G(\Af)/\K$, donc, d'après les calculs ci-dessus,
$$T_s=\sum_{(\gamma_0;\gamma,\delta)}c(\gamma_0;\gamma,\delta)TO_\delta(\phi_j^{(s)})Tr(\gamma_0,R\Gamma(Lie(\N_S),V)_{<t_1,\dots,<t_c})\sum_{i\in I}\sum_{j\in J_i}\chi(\Gamma^\circ_{p_{ij}g_i,S})O_\gamma(f_{p_{ij}g_i}).$$
Pour tous $g\in\G(\Af^p)$ et $p\in\P_{\lin,S}(\Af)$, $\K^\circ_{g,S}=\K^\circ_{pg,S}$, donc $O_\gamma(f_{g})=O_\gamma(f_{pg})$.
On en déduit que
$$\sum_{i\in I}\sum_{j\in J_i}\chi(\Gamma^\circ_{p_{ij}g_i})O_\gamma(f_{p_{ij}g_i})=\sum_{i\in I}O_\gamma(f_{g_i})\sum_{j\in J_i}\chi(\Gamma^\circ_{p_{ij}g_i}).$$
Soit $X_{\lin,S}$ l'espace symétrique de $\L'_{\lin,S}$.
Fixons $i\in I$. 
Alors
$$\L'_{\lin,S}(\Q)\sous(X_{\lin,S}\times\L'_{\lin,S}(\Af)/\K_{g_i,\lin,S})=\coprod_{j\in J_i}\Gamma^\circ_{p_{ij}g_i,S}\sous X_{\lin,S},$$
donc
$$\sum_{j\in J_i}\chi(\Gamma^\circ_{p_{ij}g_i,S})=\chi(\L_{\lin,S}(\Q)\sous(X_{\lin,S}\times\L_{\lin,S}(\Af)/\K_{g_i,\lin,S}))=\chi(\L'_{\lin,S})vol(\K_{g_i,\lin,S})^{-1}.$$

\end{preuve}

\subsection{Formule des traces pour un complexe pondéré}

\begin{theoreme}\label{th:formule_traces} On utilise les notations du corollaire \ref{formule_traces_sigma}.

Soient $t_1,\dots,t_q\in\Z\cup\{\pm\infty\}$ et $V\in D^b(Rep_{\G(\Q_\ell)})$.
Alors, pour tout $j\in\Nat^*$,
$$Tr(F^{j*},R\Gamma(M^{\K}(\G,\X)^*_\Fi,W^{\geq t_1,\dots,\geq t_q}(V)_\Fi))=$$
$$\sum_{S\subset\{1,\dots,q\}}\sum_{(\gamma_0;\gamma,\delta)\in C_{max(S),j}}c(\gamma_0;\gamma,\delta)\chi(\L_{\lin,S})TO_\delta(\phi^{(max(S))}_j)Tr(\gamma_0,R\Gamma(Lie(\N_S),V)_{\geq t_s,s\in S})$$
$$\sum_{i\in I_S}O_\gamma(f_{g_i})vol\left((g_i\K g_i^{-1}\cap\P_S(\Af))/(g_i\K g_i^{-1}\cap (\G'_s\N_S)(\Af))\right)^{-1},$$
avec la convention $max(\emptyset)=0$,
et où $(g_i)_{i\in I_S}$ est un système de représentants du double quotient $\P_S(\Af)\sous\G(\Af)/\K$ dans $\G(\Af^p)$.

\end{theoreme}

\begin{remarque} La somme sur les $S\subset\{1,\dots,q\}$ est en fait une somme sur les sous-groupes paraboliques standard de $\G$.
Le terme pour $S=\emptyset$, c'est-à-dire le terme correspondant à $\G$, est celui qui apparaît dans le théorème \ref{th_Kottwitz}.

\end{remarque}

\begin{preuvet} Le théorème résulte directement de la définition des complexes pondérés, du corollaire \ref{formule_traces_sigma}, du théorème \ref{th:simplification_formule_traces} et du fait que, pour tout $S\subset\{1,\dots,q\}$ et pour tout $\gamma_0\in\G_s(\Q)$,
$$Tr(\gamma_0,R\Gamma(Lie(\N_S),V)_{\geq t_s,s\in S})=\sum_{J\subset S}(-1)^{card(J)}Tr(\gamma_0,R\Gamma(Lie(\N_S),V)_{<t_j,j\in J}).$$

\end{preuvet}

\end{document}